%% file: cartesian_gray_double_cats.tex
\documentclass[letterpaper]{article}

\usepackage{import} 

\import{setup/}{setup}

\import{macros/}{macros}

\import{macros/}{local_macros}

\providecommand{\keywords}[1]{\textbf{Keywords: \,} #1}
\providecommand{\msc}[1]{\textbf{Mathematics Subject Classification: \,} #1}

\hypersetup
{
	pdfauthor = {Edward Morehouse} ,
	pdftitle = {Cartesian Gray-Monoidal Double Categories} ,
}

\title{Cartesian Gray-Monoidal Double Categories}

\author
{
	\small
	Edward Morehouse%
	\thanks
	{
	This research was supported by
	the ESF funded Estonian IT Academy research measure
	(project 2014-2020.4.05.19-0001).
	}
}
\date{}

\begin{document}
	\maketitle
	\subimport{content/}{outline}

\end{document}

%% file: setup/setup.tex

\usepackage{iftex} 
\usepackage{xparse} 
\usepackage{xifthen} 

\import{./}{docstyle} 
\import{./}{fontstyle} 
\import{./}{mathstyle} 
\import{./}{theoremstyle} 
\import{./}{codestyle} 
\import{./}{tikzstyle} 
\import{./}{bibstyle} 
\import{./}{hyperstyle} 

%% file: setup/docstyle.tex
\usepackage{newsec} 
\usepackage{parskip} 
\usepackage{verbatim} 
\usepackage{todonotes} 
\usepackage{comment} 
\usepackage{enumerate} 
\usepackage{multicol} 
\usepackage{authblk} 
\usepackage[hmargin=36mm , vmargin=36mm]{geometry} 

%% file: setup/fontstyle.tex
\ifpdftex
	\relax
\else
	\usepackage{fontspec}
	\defaultfontfeatures{Mapping = tex-text , Ligatures = NoCommon , Scale = MatchUppercase}
	\defaultfontfeatures[\ttfamily]{Scale = MatchLowercase}
\fi

%% file: setup/mathstyle.tex
\usepackage{amssymb}
\usepackage{amsmath}
	\numberwithin{equation}{section}
	\setlength{\multlinegap}{0pt}
\usepackage{mathtools}
\usepackage{stackengine} 
\usepackage{proof}
	\inferLineSkip=3pt 
	\inferLabelSkip=8pt 

\ifpdftex
	\usepackage{latexsym} 
	\usepackage{lmodern} 
	\usepackage{frenchmath} 
	\usepackage{upgreek} 
	\usepackage{mathabx} 
	\usepackage{mathbbol} 
\else
	\usepackage{unicode-math} 
	\unimathsetup
	{
		math-style = french , 
		partial = upright ,
		nabla = upright
	}
	\setmathfont{latinmodern-math.otf} 
\fi

\ifpdftex
	\usepackage{newunicodechar} 
	\newunicodechar{′}{{^{\prime}}}
	\newunicodechar{¨}{\text{\textasciidieresis}}
	\newunicodechar{≠}{\neq}
	\newunicodechar{−}{-}
	\newunicodechar{×}{\times}
	\newunicodechar{∧}{\wedge}
	\newunicodechar{∨}{\vee}
	\newunicodechar{⊕}{\oplus}
	\newunicodechar{⊖}{\ominus}
	\newunicodechar{⊗}{\otimes}
	\newunicodechar{⊙}{\odot}
	\newunicodechar{⊛}{\oasterisk}
	\newunicodechar{⋅}{\cdot}
	\newunicodechar{∘}{\circ}
	\newunicodechar{∙}{\bullet}
	\newunicodechar{∗}{\ast}
	\newunicodechar{°}{{}^\circ}
	\newunicodechar{⋆}{\star}
	\newunicodechar{∞}{\infty}
	\newunicodechar{∅}{\varnothing}
	\newunicodechar{⋯}{\cdots}
	\newunicodechar{…}{\ldots}
	\newunicodechar{→}{\rightarrow}
	\newunicodechar{←}{\leftarrow}
	\newunicodechar{⇸}{\stackinset{c}{-0.125ex}{c}{0.25ex}{$\shortmid$}{$\rightarrow$}}
	\newunicodechar{⇷}{\stackinset{c}{0.125ex}{c}{0.25ex}{$\shortmid$}{$\leftarrow$}}
	\newunicodechar{⤍}{\dashrightarrow}
	\newunicodechar{⤌}{\dashleftarrow}
	\newunicodechar{⇒}{\Rightarrow}
	\newunicodechar{⇐}{\Leftarrow}
	\newunicodechar{↓}{\downarrow}
	\newunicodechar{↑}{\uparrow}
	\newunicodechar{⇑}{\Uparrow}
	\newunicodechar{⇓}{\Downarrow}
	\newunicodechar{↔}{\leftrightarrow}
	\newunicodechar{↦}{\mapsto}
	\newunicodechar{◇}{\Diamond}
	\newunicodechar{□}{\square}
	\newunicodechar{⟨}{\langle}
	\newunicodechar{⟩}{\rangle}
	\newunicodechar{⌜}{\ulcorner}
	\newunicodechar{⌝}{\urcorner}
	\newunicodechar{⌞}{\llcorner}
	\newunicodechar{⌟}{\lrcorner}
	\newunicodechar{∫}{\int}
	\newunicodechar{∏}{\prod}
	\newunicodechar{∐}{\coprod}
	\newunicodechar{∑}{\sum}
	\newunicodechar{≅}{\cong}
	\newunicodechar{≃}{\simeq}
	\newunicodechar{∼}{\sim}
	\newunicodechar{≔}{\coloneqq}
	\newunicodechar{≕}{\eqqcolon}
	\newunicodechar{∈}{\in}
	\newunicodechar{∋}{\ni}
	\newunicodechar{∀}{\forall}
	\newunicodechar{∃}{\exists}
	\newunicodechar{⊤}{\top}
	\newunicodechar{⊥}{\bot}
	\newunicodechar{⊢}{\vdash}
	\newunicodechar{⊣}{\dashv}
	\newunicodechar{∆}{\Delta} 
	\newunicodechar{∇}{\nabla}
	\newunicodechar{∂}{\partial}
	\newunicodechar{∞}{\infty}
	\newunicodechar{α}{\alpha}
	\newunicodechar{β}{\beta}
	\newunicodechar{γ}{\gamma}
	\newunicodechar{δ}{\delta}
	\newunicodechar{ε}{\varepsilon}
	\newunicodechar{η}{\eta}
	\newunicodechar{ι}{\iota}
	\newunicodechar{κ}{\kappa}
	\newunicodechar{λ}{\lambda}
	\newunicodechar{ρ}{\rho}
	\newunicodechar{μ}{\mu}
	\newunicodechar{ν}{\nu}
	\newunicodechar{π}{\pi}
	\newunicodechar{φ}{\varphi}
	\newunicodechar{ψ}{\psi}
	\newunicodechar{θ}{\theta}
	\newunicodechar{σ}{\sigma}
	\newunicodechar{τ}{\tau}
	\newunicodechar{χ}{\chi}
	\newunicodechar{ξ}{\xi}
	\newunicodechar{υ}{\upsilon}
	\newunicodechar{ζ}{\zeta}
	\newunicodechar{ω}{\omega}
	\newunicodechar{Γ}{\Gamma}
	\newunicodechar{Δ}{\Delta}
	\newunicodechar{Λ}{\Lambda}
	\newunicodechar{Φ}{\Phi}
	\newunicodechar{Ψ}{\Psi}
	\newunicodechar{Ω}{\Omega}
\else
	\setmathfont[range = {bb,"02190,"02192,"021F7,"021F8,"0290C,"0290D,"0219C,"0219D,"029C4,"029C5,"025C7}]{STIXTwoMath-Regular.otf}
\fi

%% file: setup/theoremstyle.tex
\usepackage{amsthm}
	\numberwithin{equation}{section}

\newtheoremstyle
	{plainer} 
	{2ex} 
	{2ex} 
	{} 
	{} 
	{\bfseries} 
	{} 
	{\newline} 
	{} 

\theoremstyle{plainer} 
\newtheorem{theorem}{Theorem}
\newtheorem{lemma}[theorem]{Lemma}

\newtheorem{proposition}[theorem]{Proposition}

\newtheorem{definition}[theorem]{Definition}

\newtheorem{remark}[theorem]{Remark}

\numberwithin{theorem}{section}
\numberwithin{example}{section}

%% file: setup/codestyle.tex
\usepackage{xcolor} 
\usepackage{listings} 

\lstset
{
	tabsize = 2 ,
	basicstyle = \ttfamily ,
	showstringspaces = false ,
	breaklines=true ,
	breakatwhitespace=true ,
	language = ,
}

\lstset
{
	rangebeginprefix = \{--\ begin:\ ,
	rangebeginsuffix = \ --\},
	rangeendprefix = \{--\ end:\ ,
	rangeendsuffix = \ --\},
	includerangemarker = false
}

%% file: setup/tikzstyle.tex
\usepackage{tikz}

\usetikzlibrary
{
	calc ,
	matrix ,
	arrows ,
	angles ,
	quotes ,
	intersections ,
	positioning ,
	fit ,
	graphs ,
	shapes.geometric ,
	shapes.misc ,
	decorations ,
	decorations.markings ,
	decorations.pathmorphing ,
	decorations.pathreplacing ,
	decorations.text
}

\ifluatex
	\usetikzlibrary{graphdrawing}
	\usegdlibrary{trees , force}
\else
	\relax
\fi

\usepgflibrary
{
	arrows.meta
}

\usepackage{xcolor}

\tikzset
{
	diagram/.style =
	{
		line cap = butt ,
		line join = bevel ,
		arrows = -> ,
		> = angle 60 ,
		auto = left ,
		font = \small ,
		text height = 1.5ex , 
		text depth = 0.25ex , 
		inner sep = 0.25em ,
	} ,
	smaller/.append style =
	{
		font = \small
	} ,
	code node/.append style =
	{
		every node/.append style =
		{
			execute at begin node = \begin{texttt} ,
			execute at end node = \end{texttt}
		}
	} ,
	code diagram/.style =
	{
		diagram ,
		code node
	} ,
	math node/.append style =
	{
		every node/.append style =
		{
			execute at begin node = \begin{math} ,
			execute at end node = \end{math}
		}
	} ,
	math diagram/.style =
	{
		diagram ,
		math node
	} ,
	string diagram/.style =
	{
		math diagram ,
		arrows = -
	} ,
	online/.style =
	{
		shape = rectangle ,
		rounded corners = 2mm ,
		inner sep = 1.5pt ,
		fill = white ,
		opacity = 1.0 ,
		anchor = center	
	} ,
	double arrow/.style =
	{
		double distance between line centers = 2pt ,
		-implies 
	} ,
	equals/.style =
	{
		double distance between line centers = 2pt ,
		-implies , 
		arrows = -
	} ,
	mapto/.append style =
	{
		arrows = |-> ,
	} ,
	paph/.append style =
	{
		decorate ,
		decoration = {snake , segment length = 5mm , amplitude = 0.5mm} ,
	} ,
	includel/.append style =
	{
		arrows = Hooks[right]-> ,
	} ,
	includer/.append style =
	{
		arrows = Hooks[left]-> ,
	} ,
	monomorphism/.append style =
	{
		arrows = >-> ,
	} ,
	epimorphism/.append style =
	{
		arrows = ->> ,
	} ,
	cofibration/.append style =
	{
		arrows = >-> ,
	} ,
	fibration/.append style =
	{
		arrows = ->> ,
	} ,
	equivalence/.append style =
	{
		decoration = {markings , mark = at position 1/2 with {\node[online , transform shape] {∼};}} ,
		postaction = {decorate} ,
	} ,
	mark pos/.store in = \markpos ,
	mark pos = 1/2 ,
	proarrow/.append style =
	{
		postaction = {decorate} ,
		decoration =
		{
			markings ,
			mark = at position \markpos with {\draw [solid , -] (0 , -0.6ex) to (0 , 0.6ex);}
		} ,
	} ,
	proarrows/.append style =
	{
		postaction = {decorate} ,
		decoration =
		{
			markings ,
			mark = at position #1with {\draw [solid , -] (0 , -0.6ex) to (0 , 0.6ex);}
		} ,
	} ,
	heteromorphism/.append style =
	{
		decoration = {markings , mark = at position \markpos with {\draw [white , solid , semithick , -] (-1.0ex , 0) to (1.0ex , 0);}} ,
		postaction = {decorate} ,
	} ,
	arrow/.append style =
	{
		arrows = -> ,
	} ,
	string/.append style =
	{
		arrows = - ,
	} ,
	cofibration string/.append style =
	{
		decoration = {markings , mark = at position 1/4 with {\arrow{>}} , mark = at position 3/4 with {\arrow{>}}} ,
		postaction = {decorate} ,
	} ,
	fibration string/.append style =
	{
		decoration = {markings , mark = at position 4/8 with {\arrow{>}} , mark = at position 5/8 with {\arrow{>}}} ,
		postaction = {decorate} ,
	} ,
	overcross/.append style =
	{
		preaction = {draw = white , - , line width = 4pt , line cap = round}
	} ,
	label/.append style =
	{
		font = \scriptsize
	} ,
	incoming/.append style =
	{
		pos = 0 ,
		anchor = south
	} ,
	outgoing/.append style =
	{
		pos = 1 ,
		anchor = north
	} ,
	outcoming/.append style =
	{
		pos = 0 ,
		anchor = north
	} ,
	ingoing/.append style =
	{
		pos = 1 ,
		anchor = south
	} ,
	fromabove/.append style =
	{
		pos = 0 ,
		anchor = south
	} ,
	tobelow/.append style =
	{
		pos = 1 ,
		anchor = north
	} ,
	frombelow/.append style =
	{
		pos = 0 ,
		anchor = north
	} ,
	toabove/.append style =
	{
		pos = 1 ,
		anchor = south
	} ,
	fromleft/.append style =
	{
		pos = 0 ,
		anchor = east
	} ,
	toright/.append style =
	{
		pos = 1 ,
		anchor = west
	} ,
	fromright/.append style =
	{
		pos = 0 ,
		anchor = west
	} ,
	toleft/.append style =
	{
		pos = 1 ,
		anchor = east
	} ,
	forall/.append style =
	{
		line width = 0.5pt ,
	} ,
	exists/.append style =
	{
		densely dashed
	} ,
	structural/.append style =
	{
		opacity = 2/3 ,
	} ,
	0-cell/.style =
	{
		shape = rectangle ,
		rounded corners = 2mm ,
		inner sep = 2pt
	} ,
	1-cell/.style =
	{
		shape = rectangle ,
		rounded corners = 2mm ,
		inner sep = 1.5pt ,
		fill = white ,
		opacity = 1.0 ,
		anchor = center	
	} ,
	2-cell/.style =
	{
		draw = black!75 ,
		fill = white ,
		opacity = 0.9 ,
		inner sep = 0.75mm ,
		minimum size = 4.0mm
	} ,
	bead/.style =
	{
		2-cell ,
		shape = rectangle ,
		rounded corners = 2.0mm
	} ,
	diamond/.style =
	{
		2-cell ,
		shape = diamond
	} ,
	bra/.style =
	{
		2-cell ,
		shape = isosceles triangle ,
		isosceles triangle apex angle = 70 ,
		shape border rotate = -90 ,
		minimum size = 5mm
	} ,
	ket/.style =
	{
		2-cell ,
		shape = isosceles triangle ,
		isosceles triangle apex angle = 70 ,
		shape border rotate = 90 ,
		minimum size = 6mm
	} ,
	box/.style =
	{
		2-cell ,
		shape = rectangle
	} ,
	dot/.style =
	{
		2-cell ,
		shape = circle ,
		text height = 0ex ,
		text depth = 0ex ,
		inner sep = 0mm ,
		minimum size = 1.5mm
	} ,
	black dot/.style =
	{
		dot ,
		fill = gray
	} ,
	white dot/.style =
	{
		dot ,
		fill = white
	} ,
	gray dot/.style =
	{
		dot ,
		fill = gray!25
	} ,
	sheet/.style =
	{
		draw = black!50 ,
		fill = gray!10 ,
		opacity = 1.0 ,
		fill opacity = 0.5
	} ,
	torn/.append style =
	{
		decoration = {random steps , segment length = 1.25pt , amplitude = 0.75pt}
	} ,
	on sheet/.append style =
	{
		semithick
	} ,
	edge annotation/.append style =
	{
		anchor = west
	} ,
	compositor/.style =
	{
		draw = black!50 ,
		fill = gray!10 ,
		draw opacity = 0.5 ,
		fill opacity = 0.25 ,
		dashed
	} ,
	callout/.append style =
	{
		minimum height = 2mm ,
		minimum width = 4mm ,
		draw ,
		draw opacity = 0.25
	} ,
	callout-pre/.append style =
	{
		callout ,
		inner sep = 1.2mm ,
		solid
	} ,
	callout-post/.append style =
	{
		callout ,
		inner sep = 0.8mm ,
		dashed
	} ,
	pullback/.pic =
	{
		\draw [semithick , arrows = -] (-0.1 , 0.1) -- (0.1 , 0.1) -- (0.1 , -0.1) ;
	} ,
	pushout/.pic =
	{
		\draw [semithick , arrows = -] (-0.1 , 0.1) -- (-0.1 , -0.1) -- (0.1 , -0.1) ;
	} ,
}

\pgfdeclaredecoration{simple line}{start}
{
  \state{start}[width = +0pt,
                next state=step]{
    \pgfpathmoveto{\pgfpoint{0pt}{0pt}}
  }
  \state{step}[auto end on length    = 3pt,
               auto corner on length = 3pt,               
               width=+1pt]
  {
    \pgfpathlineto{\pgfpoint{1pt}{0pt}}
  }
  \state{final}
  {}
}

%% file: setup/bibstyle.tex
\usepackage
[
	backend = biber ,
	style = alphabetic ,
	maxnames = 6 ,
	url = false ,
]{biblatex} 

\ExecuteBibliographyOptions
{
	sorting = nyt ,
	backref = false ,
}

\DeclareBibliographyAlias{preprint}{online}

\bibliography{bibliography} 

%% file: setup/hyperstyle.tex
\usepackage{hyperref}
\hypersetup
{
	unicode = true ,  
	bookmarks = true , 
	bookmarksopen = true , 
	bookmarksopenlevel = 2 , 
	breaklinks = true , 
	colorlinks = false , 
	linkcolor = blue , 
	citecolor = blue , 
	urlcolor = blue , 
	pdfborder = {0 0 0} , 
	pdfborderstyle={/S/U/W 0} , 
}

%% file: macros/macros.tex
\NewDocumentCommand \define {s o m}
{%
	\hypertarget%
	{\IfValueTF{#2}{#2}{#3}}%
	{\IfBooleanTF{#1}{#3}{\emph{#3}}}%
}

\NewDocumentCommand \refer {s o m}
{%
	\hyperlink%
	{\IfValueTF {#2}{#2}{#3}}%
	{\IfBooleanTF{#1}{#3}{#3}}%
}


\NewDocumentCommand \ensuretext {m} {\textrm{#1}}

\NewDocumentCommand \centerdisplay {m} {\makebox[\linewidth][c]{$\displaystyle{#1}$}}


\NewDocumentCommand \bbb {m} {\mathbb{#1}}    
\NewDocumentCommand \scp {m} {\ensuretext{\textsc{#1}}}  


\RenewDocumentCommand \arg {o}
{
	\IfNoValueTF{#1}
	{{−}}
	{\overset{#1}{−}}
}

\NewDocumentCommand \infarg {o}
{
	\IfNoValueTF{#1}
	{{\_}}
	{\overset{#1}{\_}}
}



\NewDocumentCommand \bind {o m m} 
{
	\IfNoValueTF{#1}
	{} {\mathop{{#1}}}
	{#2} \mathrel{.} {#3}
}

\NewDocumentCommand \unaryminus {} {\scalebox{0.5}[1.0]{\( - \)}}
\NewDocumentCommand \inv {m} {{#1} ^{\unaryminus1}} 

\NewDocumentCommand \cat {m} {\bbb{#1}} 
\NewDocumentCommand \Cat {m} {\scp{#1}} 

\makeatletter
\newcommand*{\length}[1]{%
    \@tempcnta\z@
    \@for\@tempa:=#1\do{\advance\@tempcnta\@ne}%
    \the\@tempcnta%
}
\makeatother

\makeatletter
\newcount\my@repeat@count
\newcommand{\natrec}[2]{%
  \begingroup
  \my@repeat@count=\z@
  \@whilenum\my@repeat@count<#1\do{#2\advance\my@repeat@count\@ne}%
  \endgroup
}
\makeatother

\DeclarePairedDelimiter \hombrackets {(} {)}
\NewDocumentCommand \makehom {m} 
{
	\DeclareDocumentCommand \f {o m m} 
	{
		\IfNoValueTF{##1}
		{{##2} \mathrel{#1} {##3}}
		{{##1} \, \allowbreak \hombrackets{{##2} \mathrel{#1} {##3}}}
	}
}

\NewDocumentCommand \monoid {m m} 
{
	\DeclareDocumentCommand \f {m} 
	{
		\ifthenelse{\isempty{##1}}
		{
			{#2}
		}
		{
			\foreach \item [count=\index] in {##1}
			{
				\ifthenelse{\equal{\index}{1}}
				{} 
				{#1} 
				{\item}
			}
		}
	}
}

\DeclarePairedDelimiter \setbrackets {\{} {\}}
\NewDocumentCommand \set {m}
{
	\ifthenelse{\isempty{#1}}
	{∅}
	{
		\monoid{\mathpunct{,}}{\,}
		\setbrackets{\f{#1}}
	}
}

\DeclarePairedDelimiter \sequencebrackets {[} {]}
\NewDocumentCommand \sequence {m}
{
	\ifthenelse{\isempty{#1}}
	{∅}
	{
		\monoid{\mathpunct{,}}{\,}
		\sequencebrackets{\f{#1}}
	}
}


\RenewDocumentCommand \hom {o m m}
{
	\makehom{→}
	\f[#1]{#2}{#3}
}

\NewDocumentCommand \comp {O{1} m}
{
	\monoid
	{\mathbin{\natrec{#1}{\mathord{⋅}}}}
	{\mathrm{id} \ifthenelse{#1 > 1}{^{#1}}{}}
	\f{#2}
}

\NewDocumentCommand \idty {o m}
{
	\comp[#1]{} ({#2})
}

\NewDocumentCommand \pmoc {O{1} m}
{
	\monoid
	{\mathbin{\natrec{#1}{\mathord{∘}}}}
	{\mathrm{id} \ifthenelse{#1 > 1}{^{#1}}{}}
	\f{#2}
}

\NewDocumentCommand \id {m}
{
	\mathrm{id}
	\ifthenelse
	{\equal{#1}{}}
	{}
	{({#1})}
}

\NewDocumentCommand \adjoint {m}
{
	\monoid{\mathrel{⊣}}{\mathrm{id}}
	\f{#1}
}


\NewDocumentCommand \tensor {m}
{
	\monoid{⊗}{\mathrm{I}}
	\f{#1}
}






\NewDocumentCommand \product {m}
{
	\monoid{×}{1}
	\f{#1}
}

\DeclarePairedDelimiter \tuplebrackets {⟨} {⟩}
\DeclarePairedDelimiter \globaltuplebrackets {(} {)}
\NewDocumentCommand \tuple {s m}
{
	\IfBooleanTF{#1}
	{
		\monoid{\mathbin{,}}{\,}
		\globaltuplebrackets{\f{#2}}
	}
	{
		\ifthenelse{\isempty{#2}}
		{!}
		{
			\monoid{\mathbin{,}}{\,}
			\tuplebrackets{\f{#2}}
		}
	}
}


\DeclarePairedDelimiter \cotuplebrackets {[} {]}
\NewDocumentCommand \cotuple {m}
{
	\ifthenelse{\isempty{#1}}
	{¡}
	{
		\monoid{\mathbin{,}}{\,}
		\cotuplebrackets{\f{#1}}
	}
}


\NewDocumentCommand \boundary {o}
{
	∂
	{
		\IfNoValueTF{#1}
		{}
		{^{#1}}
	}
}





\NewDocumentCommand \rightclass {m m} 
{
	{#2}^{#1}
}
\NewDocumentCommand \leftclass {m m} 
{
	{}^{#1}{#2}
}

\NewDocumentCommand \since {m} {\ensuremath{[\text{#1}]}}
\NewDocumentEnvironment{relationalreasoning}{}{\begin{array}{cl}}{\end{array}}
\NewDocumentCommand \term {o o m}
{
	\IfNoValueTF{#2}
	{{#1} & {#3}\\}
	{{#1} & \since{#2} \\ & {#3} \\}
}


\NewDocumentEnvironment{maptable}{}{\begin{array}{rcl}}{\end{array}}


\NewDocumentCommand \sequent {o o m m}
{
	{#3} \mathrel{⊢\IfNoValueTF{#1}{}{^{#1}}\IfNoValueTF{#2}{}{_{#2}}}{#4}
}


\NewDocumentCommand \cancel {o m} 
{
	\IfNoValueTF{#1}
	{[{#2}]}
	{{[{#2}] ^{#1}}}
}

%% file: macros/local_macros.tex
\NewDocumentCommand \interchanger {} {χ}

\NewDocumentCommand \interchange {m}
{
	\interchanger_{\tuple*{#1}}
}

\NewDocumentCommand \swapper {} {S}

\NewDocumentCommand \swap {m m}
{
	\swapper_{\tuple*{{#1} , {#2}}}
}

\NewDocumentCommand \braiding {} {σ}

\NewDocumentCommand \braid {o m m}
{
	\braiding
	{\IfNoValueTF{#1}{}{#1}}
	\tuple*{{#2} , {#3}}
}

\NewDocumentCommand \yangbaxt {m m m}
{
	\braiding
	\tuple*{{#1} , {#2} , {#3}}
}

\NewDocumentCommand \syllepsis {} {υ}

\NewDocumentCommand \syllepsize {o m m}
{
	\syllepsis
	{\IfNoValueTF{#1}{}{#1}}
	\tuple*{{#2} , {#3}}
}

\NewDocumentCommand \comparitor {o}
{
	\IfNoValueTF{#1}
	{θ}
	{θ_{#1}}
}



\RenewDocumentCommand \lim {} {\underset{⟵}{\mathrm {lim}}}


\NewDocumentCommand \push {o m m} 
{
	\IfNoValueTF{#1}
	{
		{#2}_!
	}
	{
		Σ_{#1} ({#2})
	}
	\ifthenelse{\isempty{#3}}{}{({#3})}
}
\NewDocumentCommand \forward {o m} 
{
	\overrightarrow{#2}
}
\NewDocumentCommand \pull {o m m} 
{
	\IfNoValueTF{#1}
	{
		{#2}^*
	}
	{
		\inv{#1} ({#2})
	}
	\ifthenelse{\isempty{#3}}{}{({#3})}
}
\NewDocumentCommand \backward {o m} 
{
	\overleftarrow{#2}
}

\NewDocumentCommand \cartesianlift {o m m} 
{
	\overset{\hom{∙}{#3}}{#2}
}
\NewDocumentCommand \opcartesianlift {o m m} 
{
	\overset{\hom{#3}{∙}}{#2}
}



\NewDocumentCommand \proarrow {}{⇸}
\NewDocumentCommand \prohom {o m m}
{
	\makehom{\proarrow}
	\f[#1]{#2}{#3}
}

\NewDocumentCommand \procomp {m}
{
	\monoid{⊙}{\mathrm{U}}
	\f{#1}
}


\NewDocumentCommand \whisker {m}
{
	\monoid{⊛}{\mathrm{Id}}
	\f{#1}
}










\NewDocumentCommand \het {o m m}
{
	\makehom{⤍}
	\f[#1]{#2}{#3}
}


\NewDocumentCommand \pathhom {o m m}
{
	\makehom{⇝}
	\f[#1]{#2}{#3}
}







\NewDocumentCommand \lifts {}
{⧄}
\NewDocumentCommand \orthogonal {s m m}
{
	\IfBooleanTF{#1}
		{{#2} \mathrel{\overline{\lifts}} {#3}}
		{{#2} \mathrel{\lifts} {#3}}
}

\NewDocumentCommand \homotopy {o m m}
{
	\IfNoValueTF{#1}
	{
		\makehom{↝}
		\f{#2}{#3}
	}
	{
		\makehom{\overset{#1} {↝}}
		\f{#2}{#3}
	}
}

\NewDocumentCommand \hcomp {O{1} m}
{
	\monoid
	{\mathbin{\natrec{#1}{\mathord{+}}}}
	{0 \ifthenelse{#1 > 1}{^{#1}}{}}
	\f{#2}
}






\NewDocumentCommand \doublehom {o m m m m} 
{
	\IfNoValueTF{#1}
	{\prescript{#4\vphantom{#3}}{#2\vphantom{#5}}{}{◇}{}^{#3\vphantom{#4}}_{#5\vphantom{#2}}}
	{{#1} \, \allowbreak \hombrackets{\prescript{#4\vphantom{#3}}{#2\vphantom{#5}}{}{◇}{}^{#3\vphantom{#4}}_{#5\vphantom{#2}}}}
}




\NewDocumentCommand \fix {o m}
{
	\mathop{\mathrm{fix} {\IfValueTF {#1}{_{#1}}{}}} \ifthenelse{\isempty{#2}}{}{({#2})}
}




\NewDocumentCommand \boundrel {o m m}
{
	{#2} \mathbin{⪯} {#3}
	\IfNoValueTF{#1}{}{ : #1}
}


\NewDocumentCommand \costaction{o m m} 
{
	{#2} \mathbin{+_{\IfNoValueTF{#1}{⋅}{#1}}} {#3}
}



\NewDocumentCommand \evaluatesto {o m m}  
{
	{#2} \mathrel{\mathord{⇓}{\IfNoValueTF{#1}{}{^{#1}}}} {#3}
}

\NewDocumentCommand \diag {} {∆}


%% file: content/outline.tex
%

\begin{abstract} \pdfbookmark[1]{Abstract}{section: abstract}
	\subimport*{parts/}{abstract}

\end{abstract}

\keywords
{
	double categories,
	Gray-monoidal structure,
	surface diagrams
}

\msc
{
	18N10, 
	18N20, 
	18M30 
}

\begin{sect}{Introduction} \label{section: introduction}
	\subimport*{parts/}{introduction}

\end{sect}

\begin{sect}{Double Categories} \label{section: double categories}
	\subimport*{parts/}{double_categories}

\end{sect}

\begin{sect}{Locally Cubical Gray Categories} \label{section: locally cubical gray categories}
	\subimport*{parts/}{locally_cubical_gray_categories}

\end{sect}

\begin{sect}{Gray-Monoidal Double Categories} \label{section: gray-monoidal double categories}
	\subimport*{parts/}{gray_monoidal_double_categories}

\end{sect}

\begin{sect}{Cartesian Structure} \label{section: cartesian structure}
	\subimport*{parts/}{cartesian_structure}

\end{sect}

\begin{sect}{Conclusion} \label{section: conclusion}
	\subimport*{parts/}{conclusion}

\end{sect}


\begin{sect*}{Article Version Information} \pdfbookmark[1]{Article Version Information}{section: version_information}
	\subimport*{parts/}{version_information}

\end{sect*}

\begin{sect*}{References} \pdfbookmark[1]{References}{section: references}
	\printbibliography[heading=none]
\end{sect*}

%% file: content/parts/abstract.tex

In this paper we present cartesian structure for symmetric Gray-monoidal double categories.
To do this we first introduce locally cubical Gray categories,
which are three-dimensional categorical structures analogous to classical, locally globular, Gray categories.
The motivating example comprises double categories themselves,
together with their functors, transformations, and modifications.
A one-object locally cubical Gray category is a Gray-monoidal double category.
Braiding, syllepsis, and symmetry for these is introduced in a manner analogous to that for $2$-categories.
Adding cartesian structure requires the introduction of doubly-lax functors of double categories to manage the order of copies.
The resulting theory is algebraically rather complex,
largely due to the bureaucracy of linearizing higher-dimensional boundary constraints.
Fortunately, it has a relatively simple and compelling representation in the graphical calculus of surface diagrams,
which we present.

%% file: content/parts/introduction.tex

A monoidal extension of a categorical structure allows us to combine together multiple things as a single thing.
We can regard this as adding new structure to do the combining in a coherent way.
Alternatively, we can view a monoidal extension of a categorical structure
as a one-object instance of another categorical structure with an additional dimension,
where the monoidal product of things in the original structure
corresponds to composition in the loop space of endomorphisms of the new one.

In the classical case such monoidal structure is \emph{dimension-sup'ing},
in the sense that tensoring an $m$-cell with an $n$-cell yields an $(m \mathop{∨} n)$-cell.
However, the setting of premonoidal categories \cite{power-1997-premonoidal_categories}
gives a glimpse of another possibility.
There, the tensor product of $0$-cells yields a $0$-cell,
but to make a $1$-cell we tensor either a $1$-cell with a $0$-cell or a $0$-cell with a $1$-cell.
We can think of such a monoidal structure as \emph{dimension-summing},
in the sense that tensoring an $m$-cell with an $n$-cell yields an $(m \mathop{+} n)$-cell.
In a premonoidal category we don't obtain a piece of structure by tensoring two $1$-cells
because there are no $2$-cells.
However, we may chose to obtain a property, in the form of an \emph{interchange law},
which acts as a relation on the $1$-dimensional boundary of a forbidden $2$-cell.

One dimension up, there are more possibilities.
With a $2$-category we may obtain a $0$-cell by tensoring two $0$-cells ($0 = 0 + 0$),
a $1$-cell by tensoring a $1$-cell and $0$-cell or a $0$-cell and $1$-cell ($1 = 1 + 0 = 0 + 1$),
and a $2$-cell by tensoring a $2$-cell and $0$-cell or a $0$-cell and $2$-cell ($2 = 2 + 0 = 0 + 2$).
But we also have another possibility: a $2$-cell obtained by tensoring two $1$-cells ($2 = 1 + 1$).
This categorifies the property of interchange into the structure of an \emph{interchanger}.

Gray-monoidal structure for $2$-categories refines classical monoidal structure
by requiring a choice of ordering on the tensor of $2$-cells.
We can't tensor two $2$-cells
directly because the dimension would be too high.
But we can tensor the first $2$-cell with the $0$-cell domain of the second $2$-cell,
and tensor the $0$-cell codomain of first $2$-cell with the second $2$-cell,
or vice-versa.
This yields a horizontally consecutive pair of $2$-cells that we can compose.
A choice of ordering for the tensor of $2$-cells
determines an ordering on the $1$-cells in the boundary of the resulting composite,
which can then be manipulated by vertical composition with interchangers.

This dimension-summing monoidal structure on $2$-categories was studied by Gray, after whom it has been named.
The Gray-monoidal structure on $2$-categories can be used to define $3$-dimensional, locally globular, \emph{Gray categories}.
Gray categories are significant in that they are the algebraic structure
comprising $2$-categories themselves, together with a hierarchy of their morphisms \cite{gray-1974-formal_category_theory}.
Moreover, every fully weak tricategory is equivalent to one \cite{street-1995-coherence_for_tricategories, gurski-2013-tricategories}.
Once we have Gray categories we can recognize Gray-monoidal $2$-categories as their one-object instances.
Like ordinary monoidal $2$-categories, Gray-monoidal $2$-categories can be given a symmetric braiding structure,
which was developed by Kapranov and Voevodsky \cite{voevodsky-1994-2categories},
Baez and Neuchl \cite{baez-1996-braided_monoidal_2categories},
Day and Street \cite{street-1997-monoidal_bicategories}, and
Crans \cite{crans-1998-sylleptic_monoidal_2categories}.

Double categories are $2$-dimensional categories of cubical shape,
in the sense that there are two independent dimensions in which $1$-cells may extend,
and a $2$-cell is a square bounded by a pair of each sort of $1$-cell.
We may impose Gray-monoidal structure on these as well.
However, now we have four sorts of $1$-cell interchanger instead of just one.
Gray-monoidal structure on double categories has been investigated by Böhm \cite{bohm-2019-monoidal_product_of_double_categories}.
We use this structure to define a locally cubical analogue of Gray categories,
where the homs are double categories rather than $2$-categories.

These \emph{locally cubical Gray categories} are significant in that they are the algebraic structure
comprising double categories themselves, together with a hierarchy of their morphisms.
Moreover, they generalize classical, locally globular, Gray categories
in the sense that the the latter may be seen as instances of the former that are discrete in one dimension.
One-object locally cubical Gray categories are \emph{Gray-monoidal double categories},
which we equip with a symmetric braiding structure in a manner similar to that for $2$-categories.

A symmetric monoidal structure is \emph{cartesian} if it allows us to uniformly duplicate and delete things,
as was observed in the case of $1$-categories by Fox \cite{fox-1976-cartesian_categories}.
It may seem counterintuitive that we can have cartesian structure
in the context of the dimension-summing Gray-monoidal product.
But in fact, in this setting we can still duplicate things
provided that we impose and maintain an order on the copies.
Consequently, the relationship between one copy and two is no longer strictly functorial,
nor is it invariant under swapping the copies,
as in the case of the ordinary monoidal product.
We propose a notion of (vertical) cartesian structure for Gray-monoidal double categories
that is compatible with these constraints.

The plan of the paper is as follows.
In section \ref{section: double categories}
we review the standard definitions of double categories and their hierarchy of morphisms.
This serves to introduce the constructions from the double-categorical literature that we need
along with the graphical calculus of surface diagrams that we will use to represent them.
We observe that the hom determined by a given pair of double categories itself has the structure of a double category,
and moreover that elements of consecutive homs can be composed.

In section \ref{section: locally cubical gray categories}
we give an algebraic presentation of the structure formed by double categories, functors, transformations, and modifications,
which we call a \emph{locally cubical Gray category}.
In fact, there is a family of such structures parameterized by the variance of homogeneous interchangers.
Locally cubical Gray categories are $3$-dimensional categorical structures
that are cubical in two dimensions and globular in the third.
They constitute the cubical generalization of classical, locally globular, Gray categories.

In section \ref{section: gray-monoidal double categories}
we introduce one-object instances of locally cubical Gray categories,
which we call \emph{Gray-monoidal double categories},
and find them to be essentially Böhm's \emph{double categorical Gray-monoids}.
We extend the monoidal structure with braiding, syllepsis, and symmetry
in a manner similar to that in the globular case,
as developed by Kapranov and Voevodsky, Baez and Neuchl, Day and Street, and Crans.

In section \ref{section: cartesian structure}
we give \emph{symmetric Gray-monoidal double categories} Fox-cartesian structure
by equipping them with suitable duplication and deletion.
In order to make duplication compatible with composition
we need the notion of a \emph{doubly-lax functor} of double categories,
whose comparison cells collate copies.
We also find that we need multiple duplicator maps
in order to account for the different orders in which copies can be produced.

\paragraph{Acknowledgments}
I am grateful to Michael Shulman and two anonymous reviewers for their helpful comments.

%% file: content/parts/double_categories.tex

A double category is a $2$-dimensional categorical structure of cubical shape.
Strict double categories were introduced by Ehresmann \cite{ehresmann-1962-double_categories}.
The weak form considered here has been studied extensively by Grandis and Paré
\cite{grandis-1999-limits_in_double_categories, grandis-2004-adjoint_for_double_categories,
grandis-2008-kan_extensions_in_double_categories, grandis-2017-span_representations_of_double_categories}.
The $2$-dimensional representation of constructions in double categories using string diagrams
was studied by Myers \cite{myers-2016-string_diagrams_for_double_categories},
who gave geometric conditions sufficient to ensure that a string diagram is valid,
in the sense of representing a composable configuration of cells in a double category.

We can characterize a weak double category as a weak category internal to
the $2$-category of (suitably small) categories, functors, and natural transformations.
In the following definition, and throughout the paper, all composition is written in left-to-right order.
Our generic composition notation uses an infix dot ``$\comp{\arg , \arg}$'' and units ``$\comp{}$''.

\begin{definition}[double category]
	A (weak) \define{double category} $\cat{D}$ consists of ordinary categories $\cat{D}_0$ and $\cat{D}_1$
	together with functors:
	$$
		\begin{tikzpicture}[math diagram , x = {(24mm , 0mm)} , y = {(0mm , -8mm)} , baseline=(align.base)]
			\node [anchor = east] (D2) at (0 , 0) {\cat{D}_1 ×_{\cat{D}_0} \cat{D}_1} ;
			\node [anchor = center] (D1) at (1 , 0) {\cat{D}_1} ;
			\node [anchor = west] (D0) at (2 , 0) {\cat{D}_0} ;
			\draw (D2) to node [auto] {\procomp{\arg , \arg}} (D1) ;
			\draw (D0) to node [swap] {\procomp{}} (D1) ;
			\draw (D1.north east) to [bend left = 30] node [auto] {L} (D0.north west) ;
			\draw (D1.south east) to [bend right = 30] node [swap] {R} (D0.south west) ;
			\draw (D2.north east) to [bend left = 30] node [auto] {π_0} (D1.north west) ;
			\draw (D2.south east) to [bend right = 30] node [swap] {π_1} (D1.south west) ;
			\node (align) at (0 , 0) {} ;
		\end{tikzpicture}
	$$
	where $\tuple*{π_0 , π_1}$ is a pullback of $\tuple*{R , L}$,
	and such that
	\begin{description}
		\item[identity boundaries:]
			$\comp{\procomp{} , L}  =  \comp{} \, \cat{D}_0  =  \comp{\procomp{} , R}$
		\item[composition boundaries:]
			$\comp{\procomp{\arg , \arg} , L}  =  \comp{π_0 , L}$
			and
			$\comp{\procomp{\arg , \arg} , R}  =  \comp{π_1 , R}$
	\end{description}
	together with coherent natural isomorphisms with the following components
	\begin{description}
		\item[unitors:]
			$
				λ \tuple*{m} : \hom{\procomp{\procomp{} (L m) , m}}{m}
			$
			and
			$
				ρ \tuple*{m}: \hom{\procomp{m , \procomp{} (R m)}}{m}
			$
		\item[associator:]
			$
				κ \tuple*{m , n , p} : \hom{\procomp{(\procomp{m , n}) , p}}{\procomp{m , (\procomp{n , p})}}
			$
	\end{description}	
\end{definition}

We call objects of $\cat{D}_0$ \define[double category object]{objects} or \emph{$0$-cells}
of the double category $\cat{D}$,
morphisms of $\cat{D}_0$ its \define[double category arrow]{arrows} or \emph{vertical $1$-cells},
objects of $\cat{D}_1$ its \define[double category proarrow]{proarrows} or \emph{horizontal $1$-cells},
and morphisms of $\cat{D}_1$ its \define[double category square]{squares} or \emph{$2$-cells}.

The functors $L$ and  $R$ pick out the ``left'' and ``right''
boundary objects of a proarrow, and arrows of a square, respectively.
For a proarrow $m : \cat{D}_1$, we write ``$m : \prohom{A}{B}$'' to indicate that $L (m) = A$ and $R (m) = B$.
This parallels the boundary notation for an arrow ``$f : \hom{A}{C}$'' inherited from $\cat{D}_0$.
The functor $\procomp{}$ gives the \define[identity proarrow]{identity} proarrow on an object, and square on an arrow;
and $\procomp{\arg , \arg}$ gives the \define[proarrow composition]{composite} of consecutive proarrows,
and of squares in the proarrow dimension.
Composition of arrows, and of squares in the (strict) arrow dimension,
is inherited from $\cat{D}_0$ and $\cat{D}_1$.

The coherence of the associator and unitors can be characterized in the same way as for bicategories, namely by
Mac Lane's associator coherence ``pentagon equation'' relating terms of the form
$\hom{\procomp{(\procomp{(\procomp{l , m}) , n}) , p}}{}\procomp{l , (\procomp{m , (\procomp{n , p})})}$
and middle unit coherence ``triangle equation'' relating those of the form
$\hom{\procomp{(\procomp{m , \procomp{}}) , n}}{\procomp{m , n}}$
\cite{mac_lane-1998-categories, kelly-1964-maclane_coherence}.

If the unitor natural isomorphisms are identities
then the double category is called \define[unitary double category]{unitary}.
If the associator is an identity as well
then it is \define[strict double category]{strict}.
In the following we assume our double categories to be at least strict for identity proarrows, in the sense that
$\procomp{\procomp{} A , \procomp{} A} = \procomp{} A$ and
$λ \tuple*{\procomp{} A} = \comp{} (\procomp{} A) = ρ \tuple*{\procomp{} A}$
(by coherence the transitive equality is true in any double category).
By the triangle equation this implies
$κ \tuple*{\procomp{} A , \procomp{} A , \procomp{} A} = \comp{} (\procomp{} A)$ as well.
Such double categories are sometimes called \define[preunitary double category]{preunitary}
\cite{grandis-2019-multiple_categories}. 

A full boundary specification for a square in a double category $\cat{D}$ is written as
``$\doublehom[\doublehom[\cat{D}]{\hom{A}{C}}{\hom{B}{D}}{\prohom{A}{B}}{\prohom{C}{D}}]{f}{g}{m}{n}$'',
but any well-formed prefix may be omitted if it is clear from context or irrelevant.
A square $α$ with this boundary can be depicted as the following string diagram,
where our convention is to draw the proarrow dimension horizontally and the arrow dimension vertically.
The point representing $α$ has been ``fattened up'' to a bead to facilitate labeling.
\begin{equation} \label{walking square}

\end{equation}

Composition of squares is depicted as pasting along a compatible shared boundary morphism
in the appropriate dimension:
$$
	%
$$

By the functoriality of $\procomp{,}$ we have the equations
$$
	\procomp{(\comp{α , γ}) , (\comp{β , δ})} = \comp{(\procomp{α , β}) , (\procomp{γ , δ})}
	\quad \text{and} \quad
	\procomp{\comp{} \, m , \comp{} \, n}  =  \comp{} (\procomp{m , n})
	,
$$
the former of which is a $2$-dimensional associative law known as \define{middle-four exchange}.
These imply that each of the following diagrams has a unique interpretation.
$$
	%
$$
By the functoriality of $\procomp{}$ we have the equations
$$
	\procomp{} (\comp{f , g})  =  \comp{\procomp{} f , \procomp{} g}
	\quad \text{and} \quad
	\procomp{} (\comp{} \, A)  =  \comp{} (\procomp{} A)
	,
$$
the latter of which provides a well defined notion of a (double) \define{identity square} on an object, which we write as ``$\comp[2]{} A$''.
These imply that each of the following diagrams has a unique interpretation.
$$
	%
$$
Note the convention of suppressing the drawing of composition unit cells.

Unitor naturality implies that squares can ``slide past'' them,
in the sense that for any square $α : \doublehom{f}{g}{m}{m′}$ we have
$$
	\comp{(\procomp{\procomp{} \, f , α}) , λ m′} = \comp{λ m , α}
	\; \text{and} \;
	\comp{(\procomp{α , \procomp{} \, g}) , ρ m′} = \comp{ρ m , α}
	.
$$
$$
	%
$$
Similarly, squares can slide past an associator as
$$
	\comp{(\procomp{(\procomp{α , β}) , γ}) , κ  \tuple*{m′ , n′ , p′}}  =  \comp{κ \tuple*{m , n , p} , (\procomp{α , (\procomp{β , γ})})}
	.
$$
$$
	%
$$

A square of a double category is called \define[globular square]{globular}, or a \emph{disk},
if it has only identity boundary morphisms in some dimension.
A globular square whose boundary proarrows are identities is an \define{arrow disk},
and one whose boundary arrows are identities is a \define{proarrow disk}.
In \define{disk notation} we overload the hom notation for globular squares,
writing  ``$\prohom{f}{g}$'' for $\doublehom{f}{g}{\procomp{}}{\procomp{}}$
and ``$\hom{m}{n}$'' for $\doublehom{\comp{}}{\comp{}}{m}{n}$.
Composition of arrow- or proarrow disks in the arrow dimension is strictly unital and associative
because it is ordinary $1$-categorical composition in $\cat{D}_1$.
Composition of arrow disks in the proarrow dimension is also strictly unital and associative
by \refer[preunitary double category]{preunitarity}.

An \define{invertible disk} is one with a two-sided inverse in the dimension where it has not-necessarily-trivial boundary.
An \define{identity disk} is one of the form $\procomp{} \, f : \prohom{f}{f}$ or $\comp{} \, m : \hom{m}{m}$.
A (double) \refer{identity square} is an identity disk in both dimensions.
In order to declutter notation we may use \define{dimensional promotion}
to elide an ``$\comp{}$'' or ``$\procomp{}$'' from a subterm occurrence of an expression or diagram
when its dimension is evident from the context.

The sub-double category determined by arrow disks can be identified with a (strict) $2$-category.
Indeed, there is a functor $\Cat{Arr} : \hom{\Cat{DblCat}}{2\Cat{Cat}}$ that does this\footnote
{We haven't defined morphisms of double categories yet, but we will do so shortly.}.
It has a left adjoint that fully embeds a $2$-category as a double category where all proarrows are trivial,
making $2\Cat{Cat}$ a coreflective subcategory of $\Cat{DblCat}$
\cite{moser-2021-thesis}. 
We refer to both of these as the \define[arrow 2-category]{arrow $2$-category} of a double category.
Similarly, each double category has a \define{proarrow bicategory}.

\refer[strict double category]{Strict double categories} are congenial to diagrammatics
because they let us depict a compatible configuration of squares without explicit bracketing.
There is a coherence theorem for double categories \cite{grandis-2019-multiple_categories} 
which implies that given a diagram in a strict double category,
any elaboration of its proarrow boundary to terms of a weak double category
can be extended to the interior by elaboration with unitors and associators;
and moreover, that all ways of doing so result in diagrams representing the same composite square.
Except for the sake of emphasis, in diagrams we will omit explicit bracketing of proarrow boundaries,
as well as explicit coherator proarrow disks.
When presenting an equation that holds up to coherators, we will write ``$≅$'' rather than ``$=$''
as a reminder that coherators may be inserted to unify the proarrow boundaries.

A double category that plays an important theoretical role is the \define{walking square double category},
whose only non-identity cells of each dimension are depicted in diagram \eqref{walking square}.
There is also a \define{singleton double category} $\cat{\product{}}$
comprising just one object, one arrow, one proarrow, and one square.
Given a pair of double categories $\cat{C}$ and $\cat{D}$,
we define the \define{cartesian ordered pair double category} $\product{\cat{C} , \cat{D}}$
with $(\product{\cat{C} , \cat{D}})_i  =  \product{\cat{C}_i , \cat{D}_i}$ for $i ∈ \set{0 , 1}$
and the composition structure given factor-wise.


A functor between double categories is an internal functor between internal categories.

\begin{definition}[strict functor of double categories] \label{definition: strict double category functor}
	A (strict) \define[strict double category functor]{functor} of double categories, $F : \hom{\cat{C}}{\cat{D}}$,
	consists of a pair of functors of $1$-categories
	$F_0 : \hom{\cat{C}_0}{\cat{D}_0}$ and
	$F_1 : \hom{\cat{C}_1}{\cat{D}_1}$
	that are compatible with the structural boundary functors $L$ and $R$ in the sense that
	$\comp{F_1 , L^{\cat{D}}} = \comp{L^{\cat{C}} , F_0}$ and
	$\comp{F_1 , R^{\cat{D}}} = \comp{R^{\cat{C}} , F_0}$,
	$$

	$$
	which strictly preserve the proarrow-dimension composition structure
	in the sense that
	$\nobreak{\comp{F_0 , \procomp{}^{\cat{D}}}} = \comp{\procomp{}^{\cat{C}} , F_1}$ and
	$\comp{F_1 ×_{F_0} F_1 , \procomp{,}^{\cat{D}}} = \comp{\procomp{,}^{\cat{C}} , F_1}$,
	$$
		%
	$$
	and which strictly preserve the unitors and associator,
	in the sense that
	$$
		F_1 (λ^{\cat{C}} m)  =  λ^{\cat{D}} (F_1 m)
		, \quad
		F_1 (ρ^{\cat{C}} m)  =  ρ^{\cat{D}} (F_1 m)
		, \quad
		F_1 (κ^{\cat{C}} \tuple*{m , n , p})  =  κ^{\cat{D}} \tuple*{F_1 m , F_1 n , F_1 p}
		.
	$$
\end{definition}

Except for the sake of emphasis we will omit the subscripts on a functor's constituent maps.
For dimensional uniformity we sometimes refer to a functor image of something as a ``component''.

For a functor $F : \hom{\cat{C}}{\cat{D}}$,
the $F$-image of a $\cat{C}$-square
$φ : \doublehom[\doublehom[]{\hom{A}{A′}}{\hom{B}{B′}}{\prohom{A}{B}}{\prohom{A′}{B′}}]{f}{g}{m}{n}$
looks like the same square in $\cat{D}$, but with an ``$F$'' applied to each cell label,
as shown in the string diagram on the right.
It will be useful to regard this string diagram for $F φ$
as the geometric projection, in a dimension orthogonal to both the arrow and proarrow dimensions,
of the surface diagram formed by juxtaposing a surface representing the functor $F$
with a diagram of global elements\footnote
{
	We have not yet defined the higher-dimensional cells involved, but will do so momentarily.
}
representing the square $φ$,
as shown on the left.
$$

$$

Functors of double categories preserve arrow-dimension composition
because functors of $1$-categories preserve morphism composition.
The requirement that they also preserve proarrow-dimension composition
implies that for consecutive $\cat{C}$-squares $φ$ and $ψ$
we have $F (\procomp{φ , ψ}) = \procomp{F φ , F ψ}$,
as depicted in the equation between string diagrams on the right.
We represent either/both of these terms with the surface diagram on the left.
$$
	%
$$

In this paper we treat surface diagrams \emph{algebraically},
that is, as graphical notation for a term, or an equivalence class of terms, within a categorical structure of interest.
In order to be \define[valid diagram]{valid} a diagram must admit a block decomposition into atomic subdiagrams,
representing a composable configuration of cells.
As in the diagram above, this decomposition may not be unique;
however, the composite cells obtained by all such decompositions must be equivalent.
This is in contrast to a \emph{geometric} treatment of diagrams,
which seeks to describe conditions on diagrams as geometrical (or topological) objects
sufficient to ensure that they represent well-defined elements within the categorical structure of interest.
Such a presentation for a globular surface diagram calculus is given by Hummon \cite{hummon-2012-surface_diagrams}.

Consecutive functors of double categories $F : \hom{\cat{C}}{\cat{D}}$ and $G : \hom{\cat{D}}{\cat{E}}$
compose with $(\comp{F , G})_i = \comp{F_i , G_i}$.
This composition is associative and unital because the composition of functors of $1$-categories is.
In surface diagrams we depict functor composition by sequential juxtaposition of their respective surfaces.
We usually don't bother drawing identity functors at all.

For each double category $\cat{C}$ there is a unique functor to the \refer{singleton double category}
$\tuple{}_{\cat{C}} : \hom{\cat{C}}{\cat{\product{}}}$.
There is also a \define{diagonal functor} to the \refer{cartesian ordered pair double category} of $\cat{C}$ with itself
$\diag _{\cat{C}} : \hom{\cat{C}}{\product{\cat{C} , \cat{C}}}$
sending each object, arrow, proarrow, and square to the ordered pair consisting of two copies of its argument cell.
We also have a \define{swap functor} $\swap{\cat{C}}{\cat{D}} : \hom{\product{\cat{C} , \cat{D}}}{\product{\cat{D} , \cat{C}}}$
that transposes the factors.
In the absence of ambiguity we may omit the subscripts.


Transformations between functors between double categories come in two different dimensions,
corresponding to the arrow and proarrow dimensions of double categories themselves.
Moreover, transformations in each of these dimensions can have either -- or both -- of two possible variances,
known as ``oplax'' and ``lax''.

\begin{definition}[arrow-dimension oplax transformation] \label{definition: arrow-dimension oplax transformation}
	An \define{arrow-dimension oplax transformation} of functors of double categories
	$α : \hom[\hom[]{\cat{C}}{\cat{D}}]{F}{G}$
	consists of the following data.
	\begin{description}
		\item[{\define*[arrow transformation object-component arrow]{object-component arrows}}:]
			for each object of the domain double category $A : \cat{C}$
			an arrow of the codomain double category $α A : \hom[\cat{D}]{F A}{G A}$,
		\item[{\define*[arrow transformation proarrow-component square]{proarrow-component squares}}:]
			for each proarrow of the domain double category $m : \prohom[\cat{C}]{A}{B}$
			a square of the codomain double category $α m : \doublehom[\cat{D}]{α A}{α B}{F m}{G m}$,
		\item[{\define*[arrow transformation arrow-component disk]{arrow-component disks}}:]
			for each arrow of the domain double category $f : \hom[\cat{C}]{A}{A′}$
			an \refer{arrow disk} of the codomain double category $α f : \prohom[\cat{D}]{\comp{F f , α A′}}{\comp{α A , G f}}$.
	\end{description}
	This data is required to satisfy the following relations.
	\begin{description}
		\item[preservation of proarrow composition:]
			for an object $A$ and
			consecutive proarrows $m : \prohom{A}{B}$ and $n :  \prohom{B}{C}$
			of the domain double category we have
			\begin{equation} \label{arrow transformation proarrow composition preservation}
				α (\procomp{} \, A) = \procomp{} (α A)
				\quad \text{and} \quad
				α (\procomp{m , n}) = \procomp{α m , α n}
				,
			\end{equation}
		\item[compatibility with arrow composition:]
			for an object $A$ and
			consecutive arrows $f : \hom{A}{A′}$ and $g :  \hom{A′}{A′′}$
			of the domain double category we have
			\begin{equation} \label{arrow transformation arrow composition compatibility}
				α (\comp{} \, A) = \procomp{} (α A)
				\quad \text{and} \quad
				α (\comp{f , g}) = \procomp{(\comp{\procomp{} (F f) , α g}) , ({\comp{α f , \procomp{} (G g)}})}
				,
			\end{equation}
		\item[naturality for squares:]
			for a square $φ : \doublehom{f}{g}{m}{n}$ of the domain double category
			we have
			\begin{equation} \label{arrow transformation square naturality}
				\procomp{(\comp{F φ , α n}) , α g} ≅ \procomp{α f , (\comp{α m , G φ})}
				.
			\end{equation}
	\end{description}
\end{definition}

In surface diagrams we draw an arrow-dimension transformation as a line
that vertically separates the surfaces representing its boundary functors.

A \refer[arrow transformation proarrow-component square]{proarrow-component square}
can be seen as the projection of the juxtaposition of the arrow-dimension transformation
with the global element corresponding to the proarrow.
$$

$$
An \refer[arrow transformation arrow-component disk]{arrow-component disk}
can be seen as the projection of the juxtaposition of the arrow-dimension transformation
with the global element corresponding to the arrow.
$$
	%
$$
The oplax variance of an arrow-dimension transformation corresponds to an ``upward'' slope of its line
relative to lines representing arrows in the domain double category.

In string diagrams we usually don't label the points depicting component squares,
and instead represent them as crossings of their boundary lines,
with the line representing the transformation crossing ``behind'' the one representing the arrow or proarrow
as a mnemonic for the fact that it comes ``later''\footnote
{
	We call this the \define[prospective perspective]{prospective} perspective,
	with the opposite convention representing the \define[retrospective perspective]{retrospective} perspective
}.

Preservation of binary proarrow composition \eqref{arrow transformation proarrow composition preservation}
says that the two possible ways to read the surface diagram on the left are equal,
giving the equation between their projection string diagrams on the right.
$$

$$

Compatibility with binary arrow composition \eqref{arrow transformation arrow composition compatibility}
says that the two possible ways to read the surface diagram on the left are equal,
giving the equation between their projection string diagrams on the right.
$$
	%
$$
Note that the proarrow boundaries are equal by \refer[preunitary double category]{preunitarity}.

Preservation of nullary proarrow composition \eqref{arrow transformation proarrow composition preservation} and
compatibility with nullary arrow composition \eqref{arrow transformation arrow composition compatibility}
say that the three possible ways to read the surface diagram on the left are equal,
giving the equation between their projection string diagrams on the right.
$$
	\centerdisplay{
	%
	}
$$

Naturality for squares \eqref{arrow transformation square naturality}
says that the terms represented by the two \refer[valid diagram]{valid} boundary-preserving perturbations
of the surface diagram on the left are equal,
giving the relation between projection string diagrams on the right.
$$
	%
$$
We can unify the proarrow boundaries of these string diagrams by conjugating them by unitors.
Note that the surface diagram on the left is not itself \refer[valid diagram]{valid},
in that it does not represent a structure in $\cat{D}$,
because the definition of \refer{arrow-dimension oplax transformation} does not specify a square-component \emph{anything}.
Instead, it represents a relation between the structures represented by its valid boundary-preserving perturbations.
Geometrically, this requires that
projections of lines intersect one another only pairwise and transversely (i.e., each intersection point has a neighborhood in which it forms a crossing),
projections of points don't intersect one another at all,
and projections of lines intersect only those of points on their own boundary \cite{morehouse-2022-gray_2categories}.

In the case that the square $φ$ has trivial proarrow boundary,
\eqref{arrow transformation square naturality} specializes to the following naturality relation for arrow disks.
\begin{equation} \label{arrow transformation naturality for arrow disks} 

\end{equation}
In the case that the square $φ$ has trivial arrow boundary,
\eqref{arrow transformation square naturality} specializes to the following naturality relation for proarrow disks,
which is a strict equality by unitor naturality.
\begin{equation} \label{arrow transformation naturality for proarrow disks} 
	%
\end{equation}

An \define{arrow-dimension lax transformation}
has arrow-component disks arising from a ``downward'' slope of the line representing the transformation
relative to lines representing arrows in the domain double category.
An \define{arrow-dimension pseudo transformation}
is one that is both lax and oplax, with \refer[invertible disk]{invertible} arrow-component disks.
For an arrow-dimension pseudo transformation $α : \hom{F}{G}$
we will regard the oplax variance the ``forward'' one,
and for arrow $f : \hom{A}{B}$ write
``$α f$'' for the component disk with oplax orientation and
``$\inv{(α f)}$'' for the one with lax orientation.
The invertibility of these disks gives us the following equations,
which are strict by preunitarity.
\begin{equation} \label{equation: pseudo interchanger laws}
	\procomp{α f , \inv{(α f)}}  =  \procomp{} (\comp{F f , α B})
	\qquad , \qquad
	\procomp{\inv{(α f)} , α f}  =  \procomp{} (\comp{α A , G f})
\end{equation}
They identify the projection string diagrams of the \refer[valid diagram]{valid} boundary-preserving perturbations
of each of the following surface diagrams, resolving the tangent points.
$$

$$
An arrow-dimension transformation is \define[arrow-dimension strict transformation]{strict}
if it has \refer[identity disk]{identity} arrow-component disks.

Arrow-dimension transformations compose as suggested by their diagrammatics.
For consecutive arrow-dimension oplax transformations $α : \hom{F}{G}$ and $β : \hom{G}{H}$
and proarrow $m : \prohom{A}{B}$ we have
\begin{equation} \label{arrow transformation proarrow component composition}
	(\comp{a , β}) m  =  \comp{α m , β m}
	\quad \text{and} \quad
	(\comp{} \, F) m  =  \comp{} \, (F m)
\end{equation}
or
$$
	%
$$
and for arrow $f : \hom{A}{A′}$ we have
\begin{equation} \label{arrow transformation arrow component composition}
	(\comp{a , β}) f  =  \procomp{(\comp{α f , \procomp{} \, (β A′)}) , (\comp{\procomp{} \, (α A) , β f})}
	\quad \text{and} \quad
	(\comp{} \, F) f  =  \procomp{} \, (F f)
\end{equation}
or
$$
	%
$$

Proarrow-dimension transformations are defined similarly, transposing the roles of the arrows and proarrows.
Diagrammatically, this is just a reflection about the $\tuple*{\procomp{,} , \comp{,}}$-diagonal.
However, because of the differing strictness of composition structure in the two dimensions,
the laws are not exactly algebraically dual.

\begin{definition}[proarrow-dimension oplax transformation]
	A \define{proarrow-dimension oplax transformation} of functors of double categories
	$γ : \prohom[\hom[]{\cat{C}}{\cat{D}}]{F}{G}$
	consists of the following data.
	\begin{description}
		\item[{\define*[proarrow transformation object-component proarrow]{object-component proarrows}}:]
			for each object of the domain double category $A : \cat{C}$
			a proarrow of the codomain double category $γ A : \prohom[\cat{D}]{F A}{G A}$,
		\item[{\define*[proarrow transformation arrow-component square]{arrow-component squares}}:]
			for each arrow of the domain double category $f : \hom[\cat{C}]{A}{A′}$
			a square of the codomain double category $γ f : \doublehom[\cat{D}]{F f}{G f}{γ A}{γ A′}$,
		\item[{\define*[proarrow transformation proarrow-component disk]{proarrow-component disks}}:]
			for each proarrow of the domain double category $m : \prohom[\cat{C}]{A}{B}$
			a \refer{proarrow disk} of the codomain double category $γ m : \hom[\cat{D}]{\procomp{F m , γ B}}{\procomp{γ A , G m}}$.
	\end{description}
	This data is required to satisfy the following relations.
	\begin{description}
		\item[preservation of arrow composition:]
			for an object $A$ and
			consecutive arrows $f : \hom{A}{A′}$ and $g :  \hom{A′}{A′′}$
			of the domain double category we have
			\begin{equation} \label{proarrow transformation arrow composition preservation}
				γ (\comp{} \, A) = \comp{} (γ A)
				\quad \text{and} \quad
				γ (\comp{f , g}) = \comp{γ f , γ g}
				,
			\end{equation}
		\item[compatibility with proarrow composition:]
			for an object $A$ and
			consecutive proarrows $m : \prohom{A}{B}$ and $n :  \prohom{B}{C}$
			of the domain double category we have
			\begin{equation} \label{proarrow transformation proarrow composition compatibility}
				\begin{array}{l}
					γ (\procomp{} \, A) ≅ \comp{} (γ A)
					\quad \text{and} \quad
					\\
					γ (\procomp{m , n}) ≅ \comp{(\procomp{\comp{} (F m) , γ n}) , \inv{κ} \tuple*{F m , γ B, G n} , ({\procomp{γ m , \comp{} (G n)}})}
					,
					\\
				\end{array}
			\end{equation}
		\item[naturality for squares:]
			for a square $φ : \doublehom{f}{g}{m}{n}$ of the domain double category
			we have
			\begin{equation} \label{proarrow transformation square naturality}
				\comp{(\procomp{F φ , γ g}) , γ n} = \comp{γ m , (\procomp{γ f , G φ})}
				.
			\end{equation}
	\end{description}
\end{definition}

In surface diagrams we draw a proarrow-dimension transformation as a line
that horizontally separates the surfaces representing its boundary functors.

An \refer[proarrow transformation arrow-component square]{arrow-component square}
can be seen as the projection of the juxtaposition of the proarrow-dimension transformation
with the global element corresponding to the arrow.
$$

$$

A \refer[proarrow transformation proarrow-component disk]{proarrow-component disk}
can be seen as the projection of the juxtaposition of the proarrow-dimension transformation
with the global element corresponding to the proarrow.
$$
	%
$$
The oplax variance of a proarrow-dimension transformation corresponds to a ``leftward'' slope of its line
relative to lines representing proarrows in the domain double category.

Preservation of binary arrow composition \eqref{proarrow transformation arrow composition preservation}
says that the two possible ways to read the surface diagram on the left are equal,
giving the equation between their projection string diagrams on the right.
$$
	%
$$

Compatibility with binary proarrow composition \eqref{proarrow transformation proarrow composition compatibility}
says that the two possible ways to read the surface diagram on the left are equal,
giving the relation between their projection string diagrams on the right.
$$
	%
$$
We can unify the proarrow boundaries of the these string diagrams
as $\hom{\procomp{(\procomp{F m , F n}) , γC}}{\procomp{γ A , (\procomp{G m , G n})}}$
by conjugating the latter by associators.
Moreover, that diagram implicitly contains an associator disk $\inv{κ} \tuple*{F m , γ B, G n}$.
Such bureaucracy is the price we must pay for weak composition structure.
Fortunately, the diagrammatics keeps it out of our way when we don't care about it;
and when we do, we need only elaborate our diagrams with coherator cells.

Preservation of nullary arrow composition \eqref{proarrow transformation arrow composition preservation} and
compatibility with nullary proarrow composition \eqref{proarrow transformation proarrow composition compatibility}
say that the three possible ways to read the surface diagram on the left are equal,
giving the relation between their projection string diagrams on the right.
$$

$$
We can unify the proarrow boundaries of these string diagrams by conjugating the last one by unitors.

Naturality for squares \eqref{proarrow transformation square naturality}
says that the terms represented by the two \refer[valid diagram]{valid} boundary-preserving perturbations
of the surface diagram on the left are equal,
giving the equation between projection string diagrams on the right.
$$
	%
$$
There are again special cases for globular squares corresponding to equations
\eqref{arrow transformation naturality for arrow disks} and \eqref{arrow transformation naturality for proarrow disks}.

\define[proarrow-dimension lax transformation]{Lax},
\define[proarrow-dimension pseudo transformation]{pseudo}, and
\define[proarrow-dimension strict transformation]{strict}
variants of proarrow-dimension transformations
are defined analogously with those of arrow-dimension transformations.
The composition structure of proarrow-dimension transformations is, modulo strictness,
the $\tuple*{\procomp{,} , \comp{,}}$-reflection of that for arrow-dimension transformations,
and is obvious from their diagrammatics.
We write composites of proarrow-dimension transformations using the same notation as for composites of proarrows,
namely $\procomp{\arg , \arg}$ and $\procomp{}$.

Compatible pairs of transformations in each dimension together determine square-shaped boundaries.
Cells that inhabit these boundaries are known as cubical modifications.

\begin{definition}[modification] \label{definition: modification}
	For parallel functors of double categories $F , G , F′ , G′ : \hom{\cat{C}}{\cat{D}}$,
	\refer[arrow-dimension oplax transformation]{arrow-dimension oplax transformations} $α : \hom{F}{F′}$ and $β : \hom{G}{G′}$,
	and \refer[proarrow-dimension oplax transformation]{proarrow-dimension oplax transformations} $γ : \prohom{F}{G}$ and $δ : \prohom{F′}{G′}$,
	a (cubical) \define{modification}
	$μ : \doublehom[\doublehom[\hom[]{\cat{C}}{\cat{D}}]{\hom{F}{F′}}{\hom{G}{G′}}{\prohom{F}{G}}{\prohom{F′}{G′}}]{α}{β}{γ}{δ}$
	consists of the following data.
	\begin{description}
		\item[{\define*[modification object-component square]{object-component squares}}:]
			for each object of the domain double category $A : \cat{C}$
			a square of the codomain double category $μ A : \doublehom[\cat{D}]{α A}{β A}{γ A}{δ A}$.
	\end{description}
	This data is required to satisfy the following relations.
	\begin{description}
		\item[naturality for arrows:]
			for an arrow $f : \hom{A}{A′}$ of the domain double category
			we have
			\begin{equation} \label{modification arrow naturality}
				\procomp{α f , (\comp{μ A , δ f})}  ≅  \procomp{(\comp{γ f , μ A′}) , β f}
				,
			\end{equation}
		\item[naturality for proarrows:]
			for a proarrow $m : \prohom{A}{B}$ of the domain double category
			we have
			\begin{equation} \label{modification proarrow naturality}
				\comp{γ m , (\procomp{μ A , β m})} = \comp{(\procomp{α m , μ B}) , δ m}
				.
			\end{equation}
	\end{description}
\end{definition}

In surface diagrams we draw a modification as a point
that horizontally separates the lines representing its boundary arrow-dimension transformations
and vertically separates the lines representing its boundary proarrow-dimension transformations.

An \refer[modification object-component square]{object-component square}
can be seen as the projection of the juxtaposition of the modification
with the global element corresponding to the object.
$$

$$

Naturality for arrows \eqref{modification arrow naturality}
says that the terms represented by the two \refer[valid diagram]{valid} boundary-preserving perturbations
of the surface diagram on the left are equal,
giving the relation between projection string diagrams on the right.
$$
	%
$$
We can unify the proarrow boundaries of these string diagrams by conjugating them by unitors.

Naturality for proarrows \eqref{modification proarrow naturality}
says that the terms represented by the two \refer[valid diagram]{valid} boundary-preserving perturbations
of the surface diagram on the left are equal,
giving the equation between projection string diagrams on the right.
$$
	%
$$

We can similarly define modifications for any of the other combinations
of transformation variance,
so long as the variance in each dimension is consistent.

A modification is \define[globular modification]{globular}
if its \refer[modification object-component square]{object-component squares} are \refer[globular square]{globular}.
A globular modification is \define[invertible modification]{invertible}
if its object-component disks are \refer[invertible disk]{invertible},
and is an \define[identity modification]{identity}
if they are (double) \refer[identity square]{identity squares}.

For globular modification $μ : \prohom{α}{β}$,
\eqref{modification arrow naturality} specializes to the following naturality relation for arrows.
\begin{equation} \label{arrow disk modification naturality for arrows} 

\end{equation}
Similarly, \eqref{modification proarrow naturality} specializes to the following naturality relation for proarrows.
\begin{equation} \label{arrow disk modification naturality for proarrows} 
	%
\end{equation}
There are, of course, analogous relations for globular modifications whose components are proarrow disks.

Modifications compose component-wise.
For a modification $μ$ and a horizontally consecutive modification $ν$ we have $(\procomp{μ , ν}) A = \procomp{(μ A) , (ν A)}$,
and for a vertically consecutive modification $μ′$ we have $(\comp{μ , μ′}) A = \comp{(μ A) , (μ′ A)}$,
providing unique interpretations for the following diagrams.
$$
	%
$$
In the nullary case, for an arrow-dimension transformation $α$ and a proarrow-dimension transformation $γ$
we have $(\procomp{} \, α) A = \procomp{} (α A)$ and $(\comp{} \, γ) A = \comp{} (γ A)$,
providing unique interpretations for the following diagrams.
$$
	%
$$

For each pair of double categories $\cat{C}$ and $\cat{D}$
the functors, arrow- and proarrow-dimension transformations, and modifications bounded by them
comprise the objects, arrows, proarrows, and squares
of a hom double category \cite{grandis-2019-multiple_categories}. 
Moreover, the components of these structures provide a form of composition for elements of consecutive homs.

\begin{lemma}[transformation whiskering] \label{theorem: transformation whiskering}
	Let $C : \hom{\cat{C}′}{\cat{C}}$ and $D : \hom{\cat{D}}{\cat{D}′}$ be functors of double categories.
	For transformations (of some variance)
	$α : \hom[\hom[]{\cat{C}}{\cat{D}}]{F}{F′}$ and $γ : \prohom[\hom[]{\cat{C}}{\cat{D}}]{F}{G}$
	we can define \define[transformation whiskering]{whiskered transformations} (of the same variance) as follows.
	$$
		\whisker{C , α} : \hom[\hom[]{\cat{C}′}{\cat{D}}]{\comp{C , F}}{\comp{C , F′}}
		\quad , \quad
		\whisker{α , D} : \hom[\hom[]{\cat{C}}{\cat{D}′}]{\comp{F , D}}{\comp{F′ , D}}
		\quad ,
	$$
	$$
		\whisker{C , γ} : \prohom[\hom[]{\cat{C}′}{\cat{D}}]{\comp{C , F}}{\comp{C , G}}
		\quad , \quad
		\whisker{γ , D} : \prohom[\hom[]{\cat{C}}{\cat{D}′}]{\comp{F , D}}{\comp{G , D}}
		\quad \phantom{,}
	$$
\end{lemma}

\begin{proof}
	In the first case with oplax variance, the
	\refer[arrow transformation object-component arrow]{object-component arrows},
	\refer[arrow transformation proarrow-component square]{proarrow-component squares}, and
	\refer[arrow transformation arrow-component disk]{arrow-component disks}
	of $\whisker{C , α}$ are the respective $α$-components of the $C$-images:
	$$
		(\whisker{C , α}) A ≔ α (C A)
		\quad , \quad
		(\whisker{C , α}) m ≔ α (C m)
		\quad , \quad
		(\whisker{C , α}) f ≔ α (C f)
	$$
	Proarrow composition preservation \eqref{arrow transformation proarrow composition preservation},
	arrow composition compatibility \eqref{arrow transformation arrow composition compatibility}, and
	square naturality \eqref{arrow transformation square naturality}
	for $\whisker{C , α}$ follow from the same property for $α$.
	The other cases are similar.
\end{proof}

In surface diagrams, we use surface juxtaposition to represent whiskering.
For example, associating the following diagram to the right represents the definiendum $(\whisker{C , α}) m$,
while associating it to the left represents its definiens.
$$
	\begin{tikzpicture}[string diagram , x = {(16mm , -8mm)} , y = {(0mm , -16mm)} , z = {(24mm , 0mm)} , baseline=(align.base)]
		\coordinate (volume) at (1/2 , 1/2 , 5/2) ;
		\node at ($ (volume) + (0 , 0 , 0) $) {\cat{D}} ;
		\coordinate (sheet) at ($ (volume) + (0 , 0 , -1/2) $) ;
		\draw [sheet]
			($ (sheet) + (-1/2 , -1/2 , 0) $) to
			($ (sheet) + (1/2 , -1/2 , 0) $) to coordinate [pos = 1/2] (right)
			($ (sheet) + (1/2 , 1/2 , 0) $) to
			($ (sheet) + (-1/2 , 1/2 , 0) $) to coordinate [pos = 1/2] (left)
			cycle ;
		\draw [on sheet] (left) to  node [fromleft] {α} node [toright] {α} (right) ;
		\node at ($ (sheet) + (0 , -1/3 , 0) $) {F} ;
		\node at ($ (sheet) + (0 , 1/3 , 0) $) {F′} ;
		\coordinate (volume) at ($ (volume) + (0 , 0 , -1) $) ;
		\node at ($ (volume) + (0 , -0 , 0) $) {\cat{C}} ;
		\coordinate (sheet) at ($ (sheet) + (0 , 0 , -1) $) ;
		\draw [sheet]
			($ (sheet) + (-1/2 , -1/2 , 0) $) to
			($ (sheet) + (1/2 , -1/2 , 0) $) to
			($ (sheet) + (1/2 , 1/2 , 0) $) to
			($ (sheet) + (-1/2 , 1/2 , 0) $) to
			cycle ;
		\node at (sheet) {C} ;
		\coordinate (volume) at ($ (volume) + (0 , 0 , -1) $) ;
		\node at ($ (volume) + (0 , 0 , 0) $) {\cat{C}′} ;
		\coordinate (sheet) at ($ (sheet) + (0 , 0 , -1) $) ;
		\draw [sheet]
			($ (sheet) + (-1/2 , -1/2 , 0) $) to coordinate [pos = 1/2] (top)
			($ (sheet) + (1/2 , -1/2 , 0) $) to
			($ (sheet) + (1/2 , 1/2 , 0) $) to coordinate [pos = 1/2] (bot)
			($ (sheet) + (-1/2 , 1/2 , 0) $) to
			cycle ;
		\draw [on sheet] (top) to node [fromabove] {m} node [tobelow] {m} (bot) ;
		\node at ($ (sheet) + (-1/3 , 0 , 0) $) {A} ;
		\node at ($ (sheet) + (1/3 , 0 , 0) $) {B} ;
		\coordinate (volume) at ($ (volume) + (0 , 0 , -1) $) ;
		\node at ($ (volume) + (0 , 0 , 0) $) {\cat{\product{}}} ;
		\node (align) at (volume) {} ;
	\end{tikzpicture}
$$

\begin{lemma}[modification whiskering] \label{theorem: modification whiskering}
	Let $C : \hom{\cat{C}′}{\cat{C}}$ and $D : \hom{\cat{D}}{\cat{D}′}$ again be functors of double categories.
	For a \refer{modification}
	$μ : \doublehom[\doublehom[\hom[]{\cat{C}}{\cat{D}}]{\hom{F}{F′}}{\hom{G}{G′}}{\prohom{F}{G}}{\prohom{F′}{G′}}]{α}{β}{γ}{δ}$
	we can define \define[modification whiskering]{whiskered modifications} as follows.
	$$
		\whisker{C , μ} :
		\doublehom[\doublehom[]
			{\hom{\comp{C , F}}{\comp{C , F′}}}
			{\hom{\comp{C , G}}{\comp{C , G′}}}
			{\prohom{\comp{C , F}}{\comp{C , G}}}
			{\prohom{\comp{C , F′}}{\comp{C , G′}}}]
		{\whisker{C , α}}{\whisker{C , β}}{\whisker{C , γ}}{\whisker{C , δ}}
		\quad , \quad
				\whisker{μ , D} :
		\doublehom[\doublehom[]
			{\hom{\comp{F , D}}{\comp{F′ , D}}}
			{\hom{\comp{G , D}}{\comp{G′ , D}}}
			{\prohom{\comp{F , D}}{\comp{G , D}}}
			{\prohom{\comp{F′ , D}}{\comp{G′ , D}}}]
		{\whisker{α , D}}{\whisker{β , D}}{\whisker{γ , D}}{\whisker{δ , D}}
	$$
	$$

	$$
\end{lemma}

\begin{proof}
	In the first case, the \refer[modification object-component square]{object-component squares} of $\whisker{C , μ}$
	are the $μ$-components of the $C$-images: $(\whisker{C , μ}) A ≔ μ (C A)$.
	Naturality for arrows \eqref{modification arrow naturality}, and for proarrows \eqref{modification proarrow naturality}, for $\whisker{C , μ}$
	follows from the same property for $μ$.
	The other case is similar.
\end{proof}

\begin{lemma}[modification whiskering functoriality] \label{theorem: modification whiskering functoriality}
	Let $C : \hom{\cat{C}′}{\cat{C}}$ and $D : \hom{\cat{D}}{\cat{D}′}$ again be functors of double categories.
	For a modification $μ$ with double category boundary $\hom{\cat{C}}{\cat{D}}$,
	with a horizontally consecutive modification $ν$
	and a vertically consecutive modification $μ′$,
	we have the following equations of modifications.
	$$
		\whisker{C , (\procomp{μ , ν})} = \procomp{(\whisker{C , μ}) , (\whisker{C ,  ν})}
		\quad , \quad
		\whisker{(\procomp{μ , ν}) , D} = \procomp{(\whisker{μ , D}) , (\whisker{ν , D})}
		\quad ,
	$$
	$$
		\whisker{C , (\comp{μ , μ′})} = \comp{(\whisker{C , μ}) , (\whisker{C ,  μ′})}
		\quad , \quad
		\whisker{(\comp{μ , μ′}) , D} = \comp{(\whisker{μ , D}) , (\whisker{μ′ , D})}
		\quad \phantom{,}
	$$
	And for transformations $α : \hom{F}{F′}$ and $γ : \prohom{F}{G}$
	with double category boundary $\hom{\cat{C}}{\cat{D}}$
	we have the following equations of modifications.
	$$
		\whisker{C , \procomp{} \, α} = \procomp{} (\whisker{C , α})
		\quad , \quad
		\whisker{\procomp{} \, α , D} = \procomp{} (\whisker{α , D})
		\quad ,
	$$
	$$
		\whisker{C , \comp{} \, γ} = \comp{} (\whisker{C , γ})
		\quad , \quad
		\whisker{\comp{} \, γ , D} = \comp{} (\whisker{γ , D})
		\quad \phantom{,}
	$$
\end{lemma}

\begin{proof}
	This is because both modification composition and whiskering are defined component-wise.
	For example,
	$$
		(\whisker{C , (\procomp{μ , ν})}) A
		\, = \,
		(\procomp{μ , ν}) (C A)
		\, = \,
		\procomp{μ (C A) , ν (C A)}
		\, = \,
		\procomp{(\whisker{C , μ}) A , (\whisker{C , ν}) A}
	$$
	and
	$$
		(\whisker{C , \procomp{} \, α}) A
		\, = \,
		(\procomp{} \, α) (C A)
		\, = \,
		\procomp{} (α (C A))
		\, = \,
		\procomp{} ((\whisker{C , α}) A)
		\, = \,
		(\procomp{} (\whisker{C , α})) A
	$$
\end{proof}

\begin{lemma}[transformation interchangers] \label{theorem: transformation interchangers}
	For \refer{arrow-dimension oplax transformation}s $α : \hom[\hom[]{\cat{C}}{\cat{D}}]{F}{F′}$ and $β : \hom[\hom[]{\cat{D}}{\cat{E}}]{H}{H′}$
	and \refer{proarrow-dimension oplax transformation}s $γ : \prohom[\hom[]{\cat{C}}{\cat{D}}]{F}{G}$ and $δ : \prohom[\hom[]{\cat{D}}{\cat{E}}]{H}{I}$
	we can define \define{transformation interchanger} \refer{modification}s as follows (using \refer{disk notation} for the globular ones):
	$$
		\interchange{α , δ} :
		\doublehom[\doublehom[]
			{\hom{\comp{F , H}}{\comp{F′ , H}}}
			{\hom{\comp{F , I}}{\comp{F′ , I}}}
			{\prohom{\comp{F , H}}{\comp{F , I}}}
			{\prohom{\comp{F′ , H}}{\comp{F′ , I}}}]
		{\whisker{α , H}}{\whisker{α , I}}{\whisker{F , δ}}{\whisker{F′ , δ}}
		\quad , \quad
		\interchange{γ , β} :
		\doublehom[\doublehom[]
			{\hom{\comp{F , G}}{\comp{F′ , G}}}
			{\hom{\comp{F , G′}}{\comp{F′ , G′}}}
			{\prohom{\comp{F , G}}{\comp{F , G′}}}
			{\prohom{\comp{F′ , G}}{\comp{F′ , G′}}}]
		{\whisker{F , β}}{\whisker{G , β}}{\whisker{γ , H}}{\whisker{γ , H′}}
		\quad ,
	$$
	$$
		\interchange{α , β} :
		\prohom[\hom[]{\comp{F , H}}{\comp{F′ , H′}}]
		{\comp{(\whisker{α , H}) , (\whisker{F′ , β})}}{\comp{(\whisker{F , β}) , (\whisker{α , H′})}}
		\quad ,
	$$
	$$
		\interchange{γ , δ} :
		\hom[\prohom[]{\comp{F , H}}{\comp{G , I}}]
		{\procomp{(\whisker{γ , H}) , (\whisker{G , δ})}}{\procomp{(\whisker{F , δ}) , (\whisker{γ , I})}}		%
		\quad \phantom{,}
	$$
\end{lemma}

\begin{proof}
	In each case, we define an \refer[modification object-component square]{object-component square} of $\interchange{φ , ψ}$
	to be the arrow/proarrow-component square of $ψ$ for the object-component arrow/proarrow of $φ$.
	For example, $\interchange{α , δ} A ≔ δ (α A)$.
	Naturality for arrows \eqref{modification arrow naturality} and for proarrows \eqref{modification proarrow naturality} of $\interchange{φ , ψ}$
	both follow from the naturality for squares \eqref{arrow transformation square naturality} or \eqref{proarrow transformation square naturality} of $ψ$.
\end{proof}

In surface diagrams we use surface juxtaposition to represent transformation interchangers.
$$

	\quad \phantom{,}
$$

\begin{remark}[transformation interchanger functoriality] \label{theorem: transformation interchanger functoriality}
	Interchangers of composite transformations behave according to their diagrammatics
	because transformation interchangers are defined component-wise.
	For example, the diagram
	$$
		%
	$$
	suggests the equation
	$$
		\interchange{\comp{α , α′} , β}
		\; = \;
		\procomp{(\comp{\procomp{} (\whisker{α , H}) , \interchange{α′ , β}}) , (\comp{\interchange{α , β} , \procomp{} (\whisker{α′ , H′})})}
	$$
	This is indeed the case.  For $A : \cat{C}$ we have:
	$$
		\interchange{\comp{α , α′} , β} A
		\, = \,
		β ((\comp{α , α′}) A)
		\, = \,
		β (\comp{α A , α′ A})
		\, \overset{\eqref{arrow transformation arrow composition compatibility}}{=} \,
		\procomp{(\comp{\procomp{} (H (α A)) , β (α′ A)}) , ({\comp{β (α A) , \procomp{} (H′ (α′ A))}})}
		\, =
	$$
	$$
		\procomp{(\comp{\procomp{} ((\whisker{α , H}) A) , \interchange{α′ , β} A}) , ({\comp{\interchange{α , β}A , \procomp{} (\whisker{α′ , H′})A}})}
		\, = \,
		[\procomp{(\comp{\procomp{} (\whisker{α , H}) , \interchange{α′ , β}}) , (\comp{\interchange{α , β} , \procomp{} (\whisker{α′ , H′})})}] A
	$$
\end{remark}
In the next section we will give an algebraic presentation of the ``transverse'' composition structure of morphisms of double categories.

%% file: content/parts/locally_cubical_gray_categories.tex

A classical Gray category may be thought of as a category
enriched in $2$-dimensional \emph{globular} categories under the Gray tensor product. 
The motivating example is the structure
formed by $2$-categories, together with their functors, (oplax and/or lax) transformations, and (globular) modifications.
This arises from the fact that the category of $2$-categories is closed monoidal under the Gray tensor product, 
as was shown by Gray \cite{gray-1974-formal_category_theory}.

Analogously, a locally cubical Gray category may be thought of as
a category enriched in $2$-dimensional \emph{cubical} categories under a suitable Gray tensor product.
Our motivating example is the structure
formed by double categories, together with their functors, (oplax and/or lax) transformations, and (cubical) modifications.
It has been shown by Böhm that the category of double categories is also closed monoidal
under a cubical version of the Gray tensor product
\cite{bohm-2019-monoidal_product_of_double_categories}.

In the globular setting, $2$-dimensional categories may have transformations of several possible variances, but in only one dimension.
In contrast, double categories have transformations in two independent dimensions,
reflecting the two types of morphisms within double categories themselves.
In locally cubical Gray categories this gives rise to four distinct types of interchanger,
which we can think of as vertical--vertical, vertical--horizontal, horizontal--vertical, and horizontal--horizontal,
and where we have a choice of variance in each dimension independently.
The homogeneous cubical interchangers behave like the globular interchangers of classical Gray-categories,
while the heterogeneous ones give rise to non-globular natural squares.

This determines a family of $3$-dimensional categorical structures
that are cubical in two dimensions and globular in the third.
Rather than trying to enumerate all of the structures of different variances and strictnesses that result,
we instead present the one corresponding to the choices made in section \ref{section: double categories};
namely, \emph{\refer[preunitary double category]{preunitary}} weak \refer[double category]{double categories},
with \emph{strict} \refer[strict double category functor]{functors},
and \emph{oplax} \refer[arrow-dimension oplax transformation]{arrow-dimension}
and \refer[proarrow-dimension oplax transformation]{proarrow-dimension} transformations.
Other variants may be constructed analogously.

We break the following definition into parts
in order to introduce notation and diagrammatics as we go,
and to try to explain the (relatively simple) geometric intuition
behind the (rather cumbersome) formalism.

\begin{definition}[locally cubical Gray category: $n$-cells]
	A \define{locally cubical Gray category} $\cat{C}$
	has the following cell data.
	\begin{description}
		
		\item[{\define*[lc gray category 0-cell]{$0$-cells}}:]
			a collection of objects known as $0$-cells,
		
		\item[{\define*[lc gray category 1-cell]{$1$-cells}}:]
			for each pair of $0$-cells $A$ and $B$
			a collection of $1$-cells,
			$\hom[\cat{C}]{A}{B}$,
		
		\item[{\define*[lc gray category vertical 2-cell]{vertical $2$-cells}}:]
			for each parallel pair of $1$-cells $f , f′ : \hom[\cat{C}]{A}{B}$
			a collection of vertical $2$-cells,
			$\hom[\hom[\cat{C}]{A}{B}]{f}{f′}$,
		
		\item[{\define*[lc gray category horizontal 2-cell]{horizontal $2$-cells}}:]
			for each parallel pair of $1$-cells $f , g : \hom[\cat{C}]{A}{B}$
			a collection of horizontal $2$-cells,
			$\prohom[\hom[\cat{C}]{A}{B}]{f}{g}$,
		
		\item[{\define*[lc gray category 3-cell]{$3$-cells}}:]
			for each parallel quadruple of $1$-cells $f , g , f′ , g′ : \hom[\cat{C}]{A}{B}$,
			vertical $2$-cells $α : \hom[\hom[\cat{C}]{A}{B}]{f}{f′}$ and $β : \hom[\hom[\cat{C}]{A}{B}]{g}{g′}$,
			and horizontal $2$-cells $γ : \prohom[\hom[\cat{C}]{A}{B}]{f}{g}$ and $δ : \prohom[\hom[\cat{C}]{A}{B}]{f′}{g′}$,
			a collection of $3$-cells,
			$\doublehom[\doublehom[\hom[\cat{C}]{A}{B}]{\hom{f}{f′}}{\hom{g}{g′}}{\prohom{f}{g}}{\prohom{f′}{g′}}]{α}{β}{γ}{δ}$.
		
	\end{description}
\end{definition}

For brevity, we may omit any well-formed prefix of a boundary specification that can be inferred from context or is irrelevant.
This lets us describe a $3$-cell $φ : \doublehom{α}{β}{γ}{δ}$, leaving the lower-dimensional structure implicit.

In surface diagrams we represent
$0$-cells as volumes, 
$1$-cells as planes separating the volumes representing their boundary $0$-cells in the ``principal'' or ``transverse'' dimension,
vertical $2$-cells as lines vertically separating the planes representing their boundary $1$-cells,
horizontal $2$-cells as lines horizontally separating the planes representing their boundary $1$-cells,
and $3$-cells as points horizontally separating the lines representing their boundary vertical $2$-cells
and vertically separating the lines representing their boundary horizontal $2$-cells.
We typically ``fatten up'' these points into beads to facilitate labeling.
Thus, we may depict the $3$-cell $φ$ above as follows.
$$
	\begin{tikzpicture}[string diagram , x = {(16mm , -8mm)} , y = {(0mm , -16mm)} , z = {(32mm , 0mm)} , baseline=(align.base)]
		\coordinate (volume) at (1/2 , 1/2 , 1/2) ;
		\node at ($ (volume) + (0 , 0 , 0) $) {B} ;
		\coordinate (sheet) at ($ (volume) + (0 , 0 , -1/2) $) ;
		\draw [sheet]
			($ (sheet) + (-1/2 , -1/2 , 0) $) to coordinate [pos = 1/2] (top)
			($ (sheet) + (1/2 , -1/2 , 0) $) to coordinate [pos = 1/2] (right)
			($ (sheet) + (1/2 , 1/2 , 0) $) to coordinate [pos = 1/2] (bot)
			($ (sheet) + (-1/2 , 1/2 , 0) $) to coordinate [pos = 1/2] (left)
			cycle ;
		\draw [on sheet , name path = arrow] (left) to node [fromleft] {α} node [toright] {β} (right) ;
		\draw [on sheet , name path = proarrow] (top) to node [fromabove] {γ} node [tobelow] {δ} (bot) ;
		\path [name intersections = {of = arrow and proarrow , by = {square}}] ;
		\node [bead] at (square) {φ} ;
		\node at ($ (square) + (-1/4 , -1/4 , 0) $) {f} ;
		\node at ($ (square) + (1/4 , -1/4 , 0) $) {g} ;
		\node at ($ (square) + (-1/4 , 1/4 , 0) $) {f′} ;
		\node at ($ (square) + (1/4 , 1/4 , 0) $) {g′} ;
		\coordinate (volume) at ($ (volume) + (0 , 0 , -1) $) ;
		\node at ($ (volume) + (0 , 0 , 0) $) {A} ;
		\node (align) at (volume) {} ;
	\end{tikzpicture}
$$

\begin{definition}[locally cubical Gray category: local structure]
	For a given pair of $0$-cells $A$ and $B$,
	the $1$-cells, vertical $2$-cells, horizontal $2$-cells, and $3$-cells with $0$-cell boundary $\hom{A}{B}$
	constitute the objects, arrows, proarrows, and squares
	of a(n, in our case, \refer[preunitary double category]{preunitary} weak) \refer{double category}.
\end{definition}

We use the notation for double categories, both linear and graphical, established in section \ref{section: double categories}
to describe this local structure.

\begin{definition}[locally cubical Gray category: whiskerings]
	We may compose $1$-cells as follows.
	\begin{description}
		
		\item[{\define*[lc gray category nullary 1-cell composition]{nullary $1$-cell composition}}:]
			for each $0$-cell $A$ we have an identity $1$-cell
			$\whisker{} \, A : \hom{A}{A}$,
		
		\item[{\define*[lc gray category binary 1-cell composition]{binary $1$-cell composition}}:]
			for consecutive $1$-cells $f : \hom{A}{B}$ and $g : \hom{B}{C}$ we have a composite $1$-cell
			$\whisker{f , g} : \hom{A}{C}$.
		
	\end{description}

	We may compose a $2$-cell or $3$-cell having $0$-cell boundary $\hom{A}{B}$
	with a $1$-cell $a : \hom{A′}{A}$ or
	with a $1$-cell $b : \hom{B}{B′}$ as follows.
	\begin{description}
		
		\item[{\define*[lc gray category vertical 2-cell whiskering]{vertical $2$-cell whiskering}}:]
			for a vertical $2$-cell $α : \hom{f}{f′}$ we have vertical $2$-cells
			$\whisker{a , α} : \hom{\whisker{a , f}}{\whisker{a , f′}}$ and
			$\whisker{α , b} : \hom{\whisker{f , b}}{\whisker{f′ , b}}$,
		
		\item[{\define*[lc gray category horizontal 2-cell whiskering]{horizontal $2$-cell whiskering}}:]
			for a horizontal $2$-cell $γ : \prohom{f}{g}$ we have horizontal $2$-cells
			$\whisker{a , γ} : \prohom{\whisker{a , f}}{\whisker{a , g}}$ and
			$\whisker{γ , b} : \prohom{\whisker{f , b}}{\whisker{g , b}}$,
		
		\item[{\define*[lc gray category 3-cell whiskering]{$3$-cell whiskering}}:]
			for a $3$-cell $φ : \doublehom{α}{β}{γ}{δ}$ we have $3$-cells
			$\whisker{a , φ} : \doublehom{\whisker{a , α}}{\whisker{a , β}}{\whisker{a , γ}}{\whisker{a , δ}}$ and
			$\whisker{φ , b} : \doublehom{\whisker{α , b}}{\whisker{β , b}}{\whisker{γ , b}}{\whisker{δ , b}}$.
		
	\end{description}
\end{definition}

For the sake of dimensional uniformity, we may regard the composition of $1$-cells $\whisker{f , g}$
either as the whiskering of $g$ by $f$ on the left, or as the whiskering of $f$ by $g$ on the right.
In surface diagrams we depict the $3$-cell whiskerings $\whisker{a , φ}$ and $\whisker{φ , b}$  as follows,
where the boundaries depict the whiskerings of the respective lower-dimensional cells.
$$

$$

\begin{definition}[locally cubical Gray category: interchangers]
	For a vertical or horizontal $2$-cell having $0$-cell boundary $\hom{A}{B}$
	and a vertical or horizontal $2$-cell having $0$-cell boundary $\hom{B}{C}$
	we have $3$-cells having $0$-cell boundary $\hom{A}{C}$
	as follows.
	\begin{description}
		
		\item[{\define*[lc gray category vertical-horizontal interchanger]{vertical-horizontal interchanger}}:]
			for vertical $2$-cell $α : \hom[\hom[]{A}{B}]{f}{f′}$
			and horizontal $2$-cell $δ : \prohom[\hom[]{B}{C}]{g}{g′}$
			we have a $3$-cell
			$
				\interchange{α , δ} : \doublehom
				{\whisker{α , g}}
				{\whisker{α , g′}}
				{\whisker{f , δ}}
				{\whisker{f′ , δ}}
			$,
		
		\item[{\define*[lc gray category horizontal-vertical interchanger]{horizontal-vertical interchanger}}:]
			for horizontal $2$-cell $γ : \prohom[\hom[]{A}{B}]{f}{f′}$
			and vertical $2$-cell $β : \hom[\hom[]{B}{C}]{g}{g′}$
			we have a $3$-cell
			$
				\interchange{γ , β} : \doublehom
				{\whisker{f , β}}
				{\whisker{f′ , β}}
				{\whisker{γ , g}}
				{\whisker{γ , g′}}
			$,
		
		\item[{\define*[lc gray category vertical-vertical interchanger]{vertical-vertical interchanger}}:]
			for vertical $2$-cells $α : \hom[\hom[]{A}{B}]{f}{f′}$ and $β : \hom[\hom[]{B}{C}]{g}{g′}$
			we have a $3$-cell
			$
				\interchange{α , β} : \doublehom
					{\comp{(\whisker{α , g}) , (\whisker{f′ , β})}}
					{\comp{(\whisker{f , β}) , (\whisker{α , g′})}}
					{\procomp{} (\whisker{f , g})}
					{\procomp{} (\whisker{f′ , g′})}
			$,
		
		\item[{\define*[lc gray category horizontal-horizontal interchanger]{horizontal-horizontal interchanger}}:]
			for horizontal $2$-cells $γ : \prohom[\hom[]{A}{B}]{f}{f′}$ and $δ : \prohom[\hom[]{B}{C}]{g}{g′}$
			we have a $3$-cell
			$
				\interchange{γ , δ} : \doublehom
					{\comp{} (\whisker{f , g})}
					{\comp{} (\whisker{f′ , g′})}
					{\procomp{(\whisker{γ , g}) , (\whisker{f′ , δ})}}
					{\procomp{(\whisker{f , δ}) , (\whisker{γ , g′})}}
			$.
		
	\end{description}
\end{definition}

In surface diagrams we depict the \define{heterogeneous interchanger} $3$-cells like this:
$$

	\quad , \quad
$$
and the \define{homogeneous interchanger} $3$-cells like this:
$$
	%
	\quad . \quad
$$
It is the orientations of the homogeneous interchangers that determines the variance of a locally cubical Gray category.
Here, we have oriented them so that both vertical and horizontal $2$-cells
are ``eager'' relative to those that precede them in the transverse dimension.
This corresponds to the oplax variance.
If instead $2$-cells were ``lazy'' with respect to their predecessors
then the locally cubical Gray category would be lax in that dimension,
and if they were ``indifferent'' then it would be pseudo.

We overload the adjectives ``oplax'' and ``lax'' to refer
both to the global property of having all homogeneous interchangers of a given orientation in a given dimension,
as well as to the local property of the orientation of a particular homogeneous interchanger $3$-cell.
For example, we may say that in a locally cubical Gray category that is vertically pseudo,
each oplax interchanger of vertical $2$-cells $\interchange{α , β}$ necessarily has a lax horizontal inverse $\inv{\interchange{α , β}}$.
Note that while a locally cubical Gray category may be pseudo in either dimension,
each homogeneous interchanger $3$-cell within it must be either oplax or lax.

We may regard whiskerings and interchangers as aspects of a single dimension-summing composition operation along $0$-cells,
but with shifted indices in the sense that composing an $(m + 1)$-cell with an $(n + 1)$-cell yields an $(m + n + 1)$-cell.
When either $m = 0$ or $n = 0$ this is a whiskering,
and when $m = n = 1$ it is an interchanger.
This lets us represent both whiskerings and interchangers uniformly in surface diagrams
as juxtapositions of string diagrams embedded on consecutive surfaces.

\begin{definition}[locally cubical Gray category: whiskering laws]
	For a $3$-cell $φ : \doublehom{α}{β}{γ}{δ}$ with $0$-cell boundary $\hom{A}{B}$,
	\begin{description}
		
		\item[nullary composite whiskering:]
			identity $1$-cells are neutral for whiskering:
			\begin{equation} \label{equation: nullary composite whiskering}
				\whisker{\whisker{} \, A , φ}  =  φ  =  \whisker{φ , \whisker{} \, B}
			\end{equation}
		
		\item[binary composite whiskering:]
			whiskering by a composite $1$-cell is the same as consecutive whiskering:
			\begin{equation} \label{equation: binary composite whiskering}
				\whisker{(\whisker{a′ , a}) , φ}  =  \whisker{a′ , (\whisker{a , φ})}
				\quad \text{and} \quad
				\whisker{φ , (\whisker{b , b′})}  =  \whisker{(\whisker{φ , b}) , b′}
			\end{equation}
		
		\item[two-sided whiskering:]
			whiskering on both sides together is associative:
			\begin{equation} \label{equation: two-sided whiskering}
				\whisker{(\whisker{a , φ}) , b}  =  \whisker{a , (\whisker{φ , b})}
			\end{equation}
		
		\item[whiskering vertical functoriality:]
			whiskering is functorial with respect to vertical composition:
			\begin{equation} \label{equation: whiskering vertical functoriality}
				\begin{array}{c}
					\whisker{a , (\comp{} \, γ)}  =  \comp{} (\whisker{a , γ})
					\quad , \quad
					\whisker{a , (\comp{φ , φ′})}  =  \comp{(\whisker{a , φ}) , (\whisker{a , φ′})}
					\quad , \quad
					\\
					\whisker{(\comp{} \, γ) , b}  =  \comp{} (\whisker{γ , b})
					\quad , \quad
					\whisker{(\comp{φ , φ′}) , b}  =  \comp{(\whisker{φ , b}) , (\whisker{φ′ , b})}
					\quad \phantom{,} \quad
				\end{array}
			\end{equation}
		
		\item[whiskering horizontal functoriality:]
			whiskering is functorial with respect to horizontal composition:
			\begin{equation} \label{equation: whiskering horizontal functoriality}
				\begin{array}{c}
					\whisker{a , (\procomp{} \, α)}  =  \procomp{} (\whisker{a , α})
					\quad , \quad
					\whisker{a , (\procomp{φ , ψ})}  =  \procomp{(\whisker{a , φ}) , (\whisker{a , ψ})}
					\quad , \quad
					\\
					\whisker{(\procomp{} \, α) , b}  =  \procomp{} (\whisker{α , b})
					\quad , \quad
					\whisker{(\procomp{φ , ψ}) , b}  =  \procomp{(\whisker{φ , b}) , (\whisker{ψ , b})}
					\quad \phantom{,} \quad
				\end{array}
			\end{equation}
		
		\item[whiskering horizontal coherators:]
			whiskering preserves the horizontal composition coherators:
			\begin{equation} \label{equation: whiskering horizontal coherators}
				\begin{array}{c}
					\whisker{a , λ \tuple*{γ}}  =  λ \tuple*{\whisker{a , γ}}
					\quad , \quad
					\whisker{a , ρ \tuple*{γ}}  =  ρ \tuple*{\whisker{a , γ}}
					\quad , \quad
					\\
					\whisker{a , κ \tuple*{γ_0 , γ_1 ,  γ_2}}  =  κ \tuple*{\whisker{a , γ_0} , \whisker{a , γ_1} , \whisker{a , γ_2}}
					\quad , \quad
					\\
					\whisker{λ \tuple*{γ} , b}  =  λ \tuple*{\whisker{γ , b}}
					\quad , \quad
					\whisker{ρ \tuple*{γ} , b}  =  ρ \tuple*{\whisker{γ , b}}
					\quad , \quad
					\\
					\whisker{κ \tuple*{γ_0 , γ_1 ,  γ_2} , b}  =  κ \tuple*{\whisker{γ_0 , b} , \whisker{γ_1 , b} , \whisker{γ_2 , b}}
					\quad \phantom{,} \quad
				\end{array}
			\end{equation}
		
	\end{description}
\end{definition}

Note that these equations between $3$-cells imply equations between their respective lower-dimensional boundary cells.
The whiskering laws can be understood as saying that transverse composition with $1$-cells
is strictly associative and unital, as well as functorial for the local double categorical composition structure.

The composite whiskering laws
\eqref{equation: nullary composite whiskering}, \eqref{equation: two-sided whiskering}, and \eqref{equation: binary composite whiskering},
assert that each of following diagrams represents a unique $3$-cell.
$$

	\quad \phantom{,}
$$
Because whiskering a $2$-cell or a $3$-cell with any number of $1$-cells on either side is unambiguous
we adopt an unbracketed notation for this,
letting us write the unique $3$-cells in each of
equations \eqref{equation: binary composite whiskering} and \eqref{equation: two-sided whiskering}
as $\whisker{a′ , a , φ}$, $\whisker{φ , b , b′}$, and $\whisker{a , φ , b}$.

The whiskering functoriality laws for  binary composites in equations
\eqref{equation: whiskering vertical functoriality} and \eqref{equation: whiskering horizontal functoriality}
assert that each of following diagrams, along with its $\tuple*{\procomp{,} , \comp{,}}$-reflection, represents a unique $3$-cell.
$$

$$

\begin{definition}[locally cubical Gray category: interchanger laws]
	For consecutive $1$-cells $f : \hom{A}{B}$ and $g : \hom{B}{C}$,
	\begin{description}
		
		\item[interchanger extremal whiskering:]
			for a $2$-cell $α$ with $0$-cell boundary $\hom{A}{B}$
			and a $2$-cell $β$ with $0$-cell boundary $\hom{B}{C}$,
			each either horizontal or vertical,
			and $1$-cells $a : \hom{A′}{A}$ and $c : \hom{C}{C′}$,
			\begin{equation} \label{equation: interchanger extremal whiskering}
				\interchange{\whisker{a , α} , β}  =  \whisker{a , \interchange{α , β}}
				\quad \text{and} \quad
				\interchange{α , \whisker{β , c}}  =  \whisker{\interchange{α , β} , c}
			\end{equation}
		
		\item[interchanger medial whiskering:]
			for a $2$-cell $α$ with $0$-cell boundary $\hom{A}{B}$
			and a $2$-cell $β$ with $0$-cell boundary $\hom{B′}{C}$,
			each either horizontal or vertical,
			and $1$-cell $b : \hom{B}{B′}$,
			\begin{equation} \label{equation: interchanger medial whiskering}
				\interchange{\whisker{α , b} , β}  =  \interchange{α , \whisker{b , β}}
			\end{equation}
		
		\item[vertical-horizontal composite interchangers:]
			for vertical $2$-cells $α : \hom{f}{f′}$ and $β : \hom{f′}{f′′}$,
			and horizontal $2$-cells $γ : \prohom{g}{g′}$ and $δ : \prohom{g′}{g′′}$,
			\begin{equation} \label{equation: vertical-horizontal composite interchangers}
\begin{array}{c}
					\interchange{\comp{} \, f , γ}  =  \comp{} (\whisker{f , γ})
					\quad , \quad
					\interchange{\comp{α , β} , γ}  =  \comp{\interchange{α , γ} , \interchange{β , γ}}
					\quad , \quad
					\\
					\interchange{α , \procomp{} \, g}  =  \procomp{} (\whisker{α , g})
					\quad , \quad
					\interchange{α , \procomp{γ , δ}}  =  \procomp{\interchange{α , γ} , \interchange{α , δ}}
					\quad \phantom{,} \quad
					\\
				\end{array}
			\end{equation}
		
		\item[horizontal-vertical composite interchangers:]
			for horizontal $2$-cells $γ : \prohom{f}{f′}$ and $δ : \prohom{f′}{f′′}$,
			and vertical $2$-cells $α : \hom{g}{g′}$ and $β : \hom{g′}{g′′}$,
			\begin{equation} \label{equation: horizontal-vertical composite interchangers}
				\begin{array}{c}
					\interchange{\procomp{} \, f , α}  =  \procomp{} (\whisker{f , α})
					\quad , \quad
					\interchange{\procomp{γ , δ} , α}  =  \procomp{\interchange{γ , α} , \interchange{δ , α}}
					\quad , \quad
					\\
					\interchange{γ , \comp{} \, g}  =  \comp{} (\whisker{γ , g})
					\quad , \quad
					\interchange{γ , \comp{α , β}}  =  \comp{\interchange{γ , α} , \interchange{γ , β}}
					\quad \phantom{,} \quad
					\\
				\end{array}
			\end{equation}
		
		\item[vertical-vertical composite interchangers:]
			for vertical $2$-cells $α : \hom{f}{f′}$, $α′ : \hom{f′}{f′′}$,
			$β : \hom{g}{g′}$,  and $β′ : \hom{g′}{g′′}$,
			\begin{equation} \label{equation: vertical-vertical composite interchangers}
				\begin{array}{c}
					\interchange{\comp{} \, f , β}  =  \procomp{} (\whisker{f , β})
					\quad , \quad
					\interchange{\comp{α , α′} , β}  =  \procomp
					{
						(\comp{\procomp{} (\whisker{α , g}) , \interchange{α′ , β}}) ,
						(\comp{\interchange{α , β} , \procomp{} (\whisker{α′ , g′})})
					}
					\quad , \quad
					\\
					\interchange{α , \comp{} \, g}  =  \procomp{} (\whisker{α , g})
					\quad , \quad
					\interchange{α , \comp{β , β′}}  =  \procomp
					{
						(\comp{\interchange{α , β} , \procomp{} (\whisker{f ′ , β′})}) ,
						(\comp{\procomp{} (\whisker{f , β}) , \interchange{α , β′}})
					}
					\quad \phantom{,} \quad
					\\
				\end{array}
			\end{equation}
		
		\item[horizontal-horizontal composite interchangers:]
			for horizontal $2$-cells $γ : \hom{f}{f′}$, $γ′ : \hom{f′}{f′′}$,
			$δ : \hom{g}{g′}$,  and $δ′ : \hom{g′}{g′′}$,
			\begin{equation} \label{equation: horizontal-horizontal composite interchangers}
				\centerdisplay{
				\begin{array}{c}
					\interchange{\procomp{} f , δ}  =  \comp{λ , \inv{ρ}}
					\quad , \quad
					\interchange{\procomp{γ , γ′} , δ}  =  \comp
					{
						κ ,
						(\procomp{\comp{} (\whisker{γ , g}) , \interchange{γ′ , δ}}) ,
						\inv{κ} ,
						(\procomp{\interchange{γ , δ} , \comp{} (\whisker{γ′ , g′})}) ,
						κ
					}
					\quad , \quad
					\\
					\interchange{γ , \procomp{} g}  =  \comp{ρ , \inv{λ}}
					\quad , \quad
					\interchange{γ , \procomp{δ , δ′}}  =  \comp
					{
						\inv{κ} ,
						(\procomp{\interchange{γ , δ} , \comp{} (\whisker{f′ , δ′})}) ,
						κ ,
						(\procomp{\comp{} (\whisker{f , δ}) , \interchange{γ , δ′}}) ,
						\inv{κ}
					}
					\quad \phantom{,} \quad
					\\
				\end{array}
				}
			\end{equation}
			where we have suppressed the horizontal composition coherator indices for readability.
		
		\item[interchanger naturality:]
			for a $3$-cell $φ : \doublehom{α}{α′}{γ}{γ′}$ with $0$-cell boundary $\hom{A}{B}$,
			a horizontal $2$-cell $δ : \prohom[\hom[]{B}{C}]{g}{g′}$, and
			a vertical $2$-cell $β : \hom[\hom[]{B}{C}]{g}{g′}$,
			\begin{equation} \label{equation: interchanger naturality in first index}
				\begin{array}{c}
					\comp{(\procomp{(\whisker{φ , g}) , \interchange{α′ , δ}}) , \interchange{γ′ , δ}}
					=
					\comp{\interchange{γ , δ} , (\procomp{\interchange{α , δ} , (\whisker{φ , g′})})}
					\quad \text{and}
					\\
					\comp{\inv{ρ} , [\procomp{(\comp{(\whisker{φ , g}) , \interchange{γ′ , β}}) , \interchange{α′ , β}}] , ρ}
					=
					\comp{\inv{λ} , [\procomp{\interchange{α , β} , (\comp{\interchange{γ , β} , (\whisker{φ , g′})})}] , λ}
				\end{array}
			\end{equation}
			and for a $3$-cell $ψ : \doublehom{β}{β′}{δ}{δ′}$ with $0$-cell boundary $\hom{B}{C}$,
			a horizontal $2$-cell $γ : \prohom[\hom[]{A}{B}]{f}{f′}$, and
			a vertical $2$-cell $α : \hom[\hom[]{A}{B}]{f}{f′}$,
			\begin{equation} \label{equation: interchanger naturality in last index}
				\begin{array}{c}
					\comp{\interchange{γ , δ} , (\procomp{(\whisker{f , ψ}) , \interchange{γ , β′}})}
					=
					\comp{(\procomp{\interchange{γ , β} , (\whisker{f′ , ψ})}) , \interchange{γ , δ′}}
					\quad \text{and}
					\\
					\comp{\inv{λ} , [\procomp{\interchange{α , β} , (\comp{(\whisker{f , ψ}) , \interchange{α , δ′}})}] , λ}
					=
					\comp{\inv{ρ} , [\procomp{(\comp{\interchange{α , δ} , (\whisker{f′ , ψ})}) , \interchange{α , β′}}] , ρ}
				\end{array}
			\end{equation}
			where we have again suppressed the horizontal composition coherator indices for readability.
		
	\end{description}
\end{definition}

This concludes the definition of a locally cubical Gray category.

The interchanger-whiskering laws \eqref{equation: interchanger extremal whiskering} and \eqref{equation: interchanger medial whiskering}
can be understood as saying that transverse composition involving a pair of $2$-cells is associative.
The composite interchanger laws
\eqref{equation: vertical-horizontal composite interchangers}, \eqref{equation: horizontal-vertical composite interchangers},
\eqref{equation: vertical-vertical composite interchangers}, and \eqref{equation: horizontal-horizontal composite interchangers}
 specify the transverse composition of $2$-cells that are local composites.
Because we have only $3$ dimensions available we can't obtain a structure by composing a $3$-cell with a $2$-cell,
so we instead impose a property in the form of the interchanger naturality laws
\eqref{equation: interchanger naturality in first index} and \eqref{equation: interchanger naturality in last index}.

The algebraic presentation of these laws is quite cumbersome.
This is largely a consequence of using a $1$-dimensional notation to encode a $3$-dimensional theory.
In surface diagrams they become more perspicuous.

The interchanger whiskering laws \eqref{equation: interchanger extremal whiskering} and \eqref{equation: interchanger medial whiskering}
assert that each of following diagrams, along with its $\tuple*{\procomp{,} , \comp{,}}$-reflection, represents a unique $3$-cell.
$$

	\quad \phantom{,}
$$
Because the interchange of a pair of $2$-cells interspersed with any number of $1$-cells is unambiguous
we adopt an interchanger notation with arbitrarily many additional $1$-cell indices,
letting us write the unique $3$-cells in each of
equations \eqref{equation: interchanger extremal whiskering} and \eqref{equation: interchanger medial whiskering}
as $\interchange{a , α , β}$, $\interchange{α , β , c}$, and $\interchange{α , b , β}$.

The heterogeneous binary composite interchanger laws in equations
\eqref{equation: vertical-horizontal composite interchangers} and \eqref{equation: horizontal-vertical composite interchangers}
assert that each of following diagrams, along with its $\tuple*{\procomp{,} , \comp{,}}$-reflection, represents a unique $3$-cell.
$$

$$

The homogeneous binary composite interchanger laws in equations
\eqref{equation: vertical-vertical composite interchangers} and \eqref{equation: horizontal-horizontal composite interchangers}
assert that each of following diagrams, along with its $\tuple*{\procomp{,} , \comp{,}}$-reflection elaborated by associators,
represents a unique $3$-cell.
$$
	%
$$

The nullary composite interchanger laws in equations
\eqref{equation: vertical-horizontal composite interchangers}, \eqref{equation: horizontal-vertical composite interchangers},
\eqref{equation: vertical-vertical composite interchangers}, and \eqref{equation: horizontal-horizontal composite interchangers},
together with the whiskering functoriality laws for nullary composites in equations
\eqref{equation: whiskering vertical functoriality} and \eqref{equation: whiskering horizontal functoriality},
assert that each of following diagrams, along with its $\tuple*{\procomp{,} , \comp{,}}$-reflection, represents a unique $3$-cell.
$$

$$

Finally, the interchanger naturality laws
\eqref{equation: interchanger naturality in first index} and \eqref{equation: interchanger naturality in last index}
assert the equations obtained by boundary-preserving perturbation of each of the following surface diagrams.
$$
	%
	\quad \phantom{,}
$$

This structure was chosen so that double categories together with their morphisms constitute a model.

\begin{proposition} \label{theorem: locally cubical Gray category of double categories}
	There is a locally cubical Gray category, $\Cat{DblCat}_G$, having
	\refer[double category]{double categories} as \refer[lc gray category 0-cell]{$0$-cells},
	\refer[strict double category functor]{strict functors} as \refer[lc gray category 1-cell]{$1$-cells},
	\refer[arrow-dimension oplax transformation]{arrow-dimension oplax transformations} as \refer[lc gray category vertical 2-cell]{vertical $2$-cells},
	\refer[proarrow-dimension oplax transformation]{proarrow-dimension oplax transformations} as \refer[lc gray category horizontal 2-cell]{horizontal $2$-cells}, and
	\refer[modification]{cubical modifications} as \refer[lc gray category 3-cell]{$3$-cells}.
\end{proposition}

\begin{proof}
	\refer[lc gray category nullary 1-cell composition]{Nullary} and
	\refer[lc gray category binary 1-cell composition]{binary composite} $1$-cells
	are identity and composite \refer[strict double category functor]{functors}.
	\refer[lc gray category vertical 2-cell whiskering]{Whiskering of a $2$-cell by a $1$-cell} and
	\refer[lc gray category 3-cell whiskering]{whiskering of a $3$-cell by a $1$-cell}
	are transformation whiskering (lemma \ref{theorem: transformation whiskering}),
	and modification whiskering (lemma \ref{theorem: modification whiskering}), respectively.
	\refer[lc gray category vertical-horizontal interchanger]{$2$-cell interchangers}
	are transformation interchangers (lemma \ref{theorem: transformation interchangers}).
	
	The whiskering of a modification by a functor is unital \eqref{equation: nullary composite whiskering}
	and associative \eqref{equation: binary composite whiskering}, \eqref{equation: two-sided whiskering}
	because it is defined component-wise.
	It is compatible with local composition \eqref{equation: whiskering vertical functoriality}, \eqref{equation: whiskering horizontal functoriality}
	by modification whiskering functoriality (lemma \ref{theorem: modification whiskering functoriality}).
	It preserves the weak composition coherators \eqref{equation: whiskering horizontal coherators} on one side by definition,
	and on the other side because \refer[strict double category functor]{(strict) functors} do.
	
	The interchanger--whiskering compatibility laws \eqref{equation: interchanger extremal whiskering}, \eqref{equation: interchanger medial whiskering}
	hold because both are defined component-wise.
	The composite interchanger laws
	\eqref{equation: vertical-horizontal composite interchangers}, \eqref{equation: horizontal-vertical composite interchangers},
	\eqref{equation: vertical-vertical composite interchangers}, \eqref{equation: horizontal-horizontal composite interchangers}
	hold by transformation interchanger functoriality (remark \ref{theorem: transformation interchanger functoriality}).
	Finally, interchanger naturality \eqref{equation: interchanger naturality in first index}, \eqref{equation: interchanger naturality in last index}
	holds by transformation naturality for squares
	\eqref{arrow transformation square naturality}, \eqref{proarrow transformation square naturality}
	if the rightmost cell is an arrow/proarrow-dimension transformation,
	and by modification naturality for arrows/proarrows
	\eqref{modification arrow naturality}, \eqref{modification proarrow naturality}
	if it is a modification.
\end{proof}

Locally cubical Gray categories generalize classical Gray categories in the following sense.

\begin{proposition}
	There is an isomorphism between (locally globular) Gray categories
	and locally cubical Gray categories having only trivial horizontal $2$-cells.
\end{proposition}

\begin{proof}
	Apply the coreflection $\Cat{Arr} : \hom{\Cat{DblCat}}{2\Cat{Cat}}$ to the homs
	and observe that the whiskering and interchanger laws for Gray categories
	are the restrictions of those for locally cubical Gray categories
	involving only vertical $2$-cells and globular $3$-cells.
\end{proof}

\paragraph{Related Constructions}
The structure of a locally cubical Gray category
is a cubical generalization of that of a classical Gray category \cite{gray-1974-formal_category_theory},
with the increased complexity arising mostly from the proliferation of interchangers.
Regarding whiskerings and interchangers
as aspects of a single transverse-dimension composition operation,
together with the graphical calculus of surface diagrams,
makes this added complexity quite manageable.

Another construction closely related to locally cubical Gray categories
is that of Garner and Gurski's \define{locally cubical bicategories} 
\cite{garner-2009-tricategory_structures}.
These can be understood as bicategories enriched in double categories;
that is, using the weak monoidal product, rather than the Gray one.
Locally cubical bicategories are weaker than locally cubical Gray categories
in the sense that transverse composition is weakly, rather than strictly, associative and unital.
But they are also stricter in the sense that the transverse composition of horizontal or vertical $2$-cells
yields another $2$-cell of the same sort rather than an interchanger $3$-cell,
and $3$-cells can be composed with one another to yield another $3$-cell.

%% file: content/parts/gray_monoidal_double_categories.tex

Whenever we have some sort of enriched categorical structure we can consider its one-object instances,
giving us some sort of monoidal structure on the hom.
In the case of a classical, locally globular, Gray category
such a structure is known as a ``Gray monoid'', a ``semistrict monoidal $2$-category'', or as we will call it, a ``Gray-monoidal $2$-category''.
In the case of a locally cubical Gray category we have something similar,
except that the carrier of the Gray-monoidal structure is a double category instead of a $2$-category.

\begin{definition}[Gray-monoidal double category]
	A \define{Gray-monoidal double category} is a one-object \refer{locally cubical Gray category},
	regarded as the hom double category of endomorphisms
	endowed with Gray-monoidal structure by the transverse composition.
\end{definition}

This has the effect of shifting everything down by a dimension,
moving the local structure of the hom double category to center stage.
The only $0$-cell, $⋆$, goes away.
Its identity $1$-cell, $\whisker{} \, ⋆$, becomes the tensor unit object, written ``$\tensor{}$''.
Other $1$-cells, $A , B : \hom{⋆}{⋆}$, become objects as well,
and their composition, $\whisker{A , B}$, becomes their tensor product, written ``$\tensor{A , B}$''.
Whiskering by a $1$-cell, $\whisker{\arg , A}$ or $\whisker{A , \arg}$,
becomes tensoring with an object, $\tensor{\arg , A}$ or $\tensor{A , \arg}$.
Interchanger $3$-cells for vertical and horizontal $2$-cells become interchanger squares for arrows and proarrows.
In short, the dimension-summing transverse combining operation has its indices un-shifted
so that combining an $m$-cell with an $n$-cell yields an $(m + n)$-cell,
manifesting as a tensor when one of $m$ or $n$ is $0$, and as an interchanger when both are $1$.
The surface diagram calculus remains unchanged, except that we never need to label the volumes.

We can think of a Gray-monoidal double category as a double category with additional structure.
In order to do this we introduce the Gray tensor product of double categories.
In the globular setting, Gray constructed a tensor product functor for $2$-categories
$\tensor{\arg , \arg} : \hom{\product{2\Cat{Cat} , 2\Cat{Cat}}}{2\Cat{Cat}}$
and showed it to be left adjoint to the internal hom in the sense that $\adjoint{\tensor{\arg , \cat{A}} , \hom{\cat{A}}{\arg}}$
\cite{gray-1974-formal_category_theory}. 
In the cubical setting, Böhm characterized the corresponding functor for double categories
$\tensor{\arg , \arg} : \hom{\product{\Cat{DblCat} , \Cat{DblCat}}}{\Cat{DblCat}}$
\cite{bohm-2019-monoidal_product_of_double_categories} 
by applying results about representability \cite{gurski-2017-gray_tensor_product}
to the fact that $\Cat{DblCat}$ is locally presentable \cite{pronk-2007-model_double_categories}
with the \refer{walking square double category} as a strong generator.

Recall that for double categories $\cat{C}$ and $\cat{D}$ we can form the \refer{cartesian ordered pair double category}
$\product{\cat{C} , \cat{D}}$.
This gives the object map of a cartesian product functor for double categories
$\product{\arg , \arg} : \hom{\product{\Cat{DblCat} , \Cat{DblCat}}}{\Cat{DblCat}}$.
If we instead use Böhm's tensor product functor
$\tensor{\arg , \arg} : \hom{\product{\Cat{DblCat} , \Cat{DblCat}}}{\Cat{DblCat}}$
then we get a different ordered pair double category.

\begin{definition}[Gray ordered pair double category] \label{definition: Gray ordered pair double category}
	For double categories $\cat{C}$ and $\cat{D}$,
	let $\cat{G}$ be the \refer{locally cubical Gray category} freely generated by them as consecutive hom double categories,
	say $\cat{C} = \hom{X}{Y}$ and $\cat{D} = \hom{Y}{Z}$.
	The \define{Gray ordered pair double category} $\tensor{\cat{C} , \cat{D}}$
	is then the transitive hom double category $\hom{X}{Z}$ in $\cat{G}$.
\end{definition}

Explicitly, a generating $(\tensor{\cat{C} , \cat{D}})$-object (of type $\tuple*{0 , 0}$, which is the only type for objects)
is the ordered pair of a $\cat{C}$-object and a $\cat{D}$-object.
This corresponds to the transverse composition of $1$-cells in the locally cubical Gray category $\cat{G}$.
A generating $(\tensor{\cat{C} , \cat{D}})$-arrow of type $\tuple*{1 , 0}$ is the ordered pair of a $\cat{C}$-arrow and a $\cat{D}$-object,
and one of type $\tuple*{0 , 1}$ is the ordered pair of a $\cat{C}$-object and a $\cat{D}$-arrow.
These correspond to the two types of whiskerings of vertical $2$-cells by $1$-cells in $\cat{G}$.
Similarly, we have generating $(\tensor{\cat{C} , \cat{D}})$-proarrows of types $\tuple*{1 , 0}$ and $\tuple*{0 , 1}$,
corresponding to horizontal $2$-cell whiskerings in $\cat{G}$.
A generating $(\tensor{\cat{C} , \cat{D}})$-square is one of the following possible types.
Generating type $\tuple*{2 , 0}$ and $\tuple*{0 , 2}$ squares are ordered pairs of a square and an object.
These correspond to whiskerings of $3$-cells by $1$-cells in $\cat{G}$.
Then there are four subtypes of generating type $\tuple*{1 , 1}$ squares,
which we call $\tuple*{v , v}$, $\tuple*{v , h}$, $\tuple*{h , v}$, and $\tuple*{h , h}$,
consisting of ordered pairs of arrows or proarrows.
These correspond to the four types of cubical interchangers,
$\product{\set{\text{vertical} , \text{horizontal}} , \set{\text{vertical} , \text{horizontal}}}$ in $\cat{G}$. 

Note that the variances of the homogeneous type $\tuple*{1 , 1}$ $(\tensor{\cat{C} , \cat{D}})$-squares are parameters
corresponding to the variances of the homogeneous interchangers in $\cat{G}$.
When we need to be explicit about these we can add annotations, such as
``$(\tensor{\cat{C} , \cat{D}})_{\tuple*{\mathrm{lax},\mathrm{pseudo}}}$''
for a Gray ordered pair double category with lax generating $\tuple*{v , v}$ disks
and both oplax and lax generating $\tuple*{h , h}$ disks.

The composition structure and relations are precisely those of a locally cubical Gray category.
For example, by whiskering vertical functoriality \eqref{equation: whiskering vertical functoriality}
we can ``merge'' consecutive $\tuple*{1 , 0}$- or $\tuple*{0 , 1}$-type arrows
\begin{equation} \label{pair binary arrow composite}
	\comp{\tuple*{f , X} , \tuple*{g , X}}  =  \tuple*{\comp{f , g} , X}
	\quad \text{and} \quad
	\comp{\tuple*{A , i} , \tuple*{A , j}}  =  \tuple*{A , \comp{i , j}}
	,
\end{equation}
and an ordered pair containing an identity arrow is itself an identity arrow
\begin{equation} \label{pair nullary arrow composite}
	\tuple*{\comp{} \, A , X}  =  \comp{} \tuple*{A , X}  =  \tuple*{A , \comp{} \, X}
	.
\end{equation}
By whiskering horizontal functoriality \eqref{equation: whiskering horizontal functoriality}
we have the same results for proarrows
\begin{equation} \label{pair binary proarrow composite}
	\procomp{\tuple*{m , X} , \tuple*{n , X}}  =  \tuple*{\procomp{m , n} , X}
	\quad \text{and} \quad
	\procomp{\tuple*{A , p} , \tuple*{A , q}}  =  \tuple*{A , \procomp{p , q}}
	,
\end{equation}
with identity proarrow
\begin{equation} \label{pair nullary proarrow composite}
	\tuple*{\procomp{} \, A , X}  =  \procomp{} \tuple*{A , X}  =  \tuple*{A , \procomp{} \, X}
	.
\end{equation}

Compound morphisms are formed by composing along compatible boundaries.
For example, given arrows $f : \hom[\cat{C}]{A}{B}$ and $i : \hom[\cat{D}]{X}{Y}$
we can form the arrow $\comp{\tuple*{f , X} , \tuple*{B , i}} : \hom[(\tensor{\cat{C} , \cat{D}})]{\tuple*{A , X}}{\tuple*{B , Y}}$
by composing along $\tuple*{B , X}$.

We represent squares of the Gray ordered pair double category $\tensor{\cat{C} , \cat{D}}$
as surface diagrams containing two surfaces, one for each of $\cat{C}$ and $\cat{D}$.
For example, the following diagram represents the square
$
	\comp{(\procomp{\tuple*{α , X} , \tuple*{g , p}}) , (\procomp{\tuple*{n , i} , \tuple*{B′ , φ}})}
$,
which by \refer{middle-four exchange} is equal to
$
	\procomp{(\comp{\tuple*{α , X} , \tuple*{n , i}}) , (\comp{\tuple*{g , p} , \tuple*{B′ , φ}})}
$.
$$
	\begin{tikzpicture}[string diagram , x = {(16mm , -8mm)} , y = {(0mm , -16mm)} , z = {(24mm , 0mm)} , baseline=(align.base)]
		\coordinate (sheet) at (1/2 , 1/2 , 1) ;
		\draw [sheet]
			($ (sheet) + (-1/2 , -1/2 , 0) $) to coordinate [pos = 2/3] (top)
			($ (sheet) + (1/2 , -1/2 , 0) $) to coordinate [pos = 2/3] (right)
			($ (sheet) + (1/2 , 1/2 , 0) $) to coordinate [pos = 1/3] (bot)
			($ (sheet) + (-1/2 , 1/2 , 0) $) to coordinate [pos = 1/3] (left)
			cycle ;
		\draw [on sheet , name path = pro] (top) to node [fromabove] {p} node [tobelow] {q} (bot) ;
		\draw [on sheet , name path = arr] (left) to node [fromleft] {i} node [toright] {j} (right) ;
		\path [name intersections = {of = pro and arr , by = {cell}}] ;
		\node [bead] at (cell) {φ} ;
		\node at ($ (sheet) + (-1/6 , -1/6 , 0) $) {X} ;
		\coordinate (sheet) at ($ (sheet) + (0 , 0 , -1) $) ;
		\draw [sheet]
			($ (sheet) + (-1/2 , -1/2 , 0) $) to coordinate [pos = 1/3] (top)
			($ (sheet) + (1/2 , -1/2 , 0) $) to coordinate [pos = 1/3] (right)
			($ (sheet) + (1/2 , 1/2 , 0) $) to coordinate [pos = 2/3] (bot)
			($ (sheet) + (-1/2 , 1/2 , 0) $) to coordinate [pos = 2/3] (left)
			cycle ;
		\draw [on sheet , name path = pro] (top) to node [fromabove] {m} node [tobelow] {n} (bot) ;
		\draw [on sheet , name path = arr] (left) to node [fromleft] {f} node [toright] {g} (right) ;
		\path [name intersections = {of = pro and arr , by = {cell}}] ;
		\node [bead] at (cell) {α} ;
		\node at ($ (sheet) + (1/6 , 1/6 , 0) $) {B′} ;
		\node (align) at (sheet) {} ;
	\end{tikzpicture}
$$

A \refer{Gray-monoidal double category} $\cat{C}$ is then a double category with additional structure in the following sense.
We can define a functor
$\tensor{,}_{\cat{C}} : \hom{\cat{\tensor{\cat{C} , \cat{C}}}}{\cat{C}}$,
which sends an ordered pair to the corresponding whiskering or interchanger.
We can also define a functor
$\tensor{}_{\cat{C}} : \hom{\cat{\product{}}}{\cat{C}}$,
which picks out the unit for $\tensor{,}$.
A double category is Gray-monoidal when these functors determine a monoid.
Intuitively, this lets us combine a finite sequence of $\cat{C}$-cells into a single cell in this Gray, dimension-summing way.
We use the same surface diagram representation for a Gray-monoidal double category $\cat{C}$
as for the iterated Gray ordered pair double category $\tensor{\cat{C} , ⋯ , \cat{C}}$.

Using \refer[Gray ordered pair double category]{Gray ordered pair double categories} we can define functors
analogous to \refer[swap functor]{swap functors} for \refer[cartesian ordered pair double category]{cartesian ordered pair double categories}.

\begin{definition}[Gray swap functor]
	Let $a , p ∈ \set{\mathrm{oplax} , \mathrm{lax} , \mathrm{pseudo}}$,
	and let $\, \tilde{}\, $ be the involution on this set that transposes $\mathrm{oplax}$ and $\mathrm{lax}$, and fixes $\mathrm{pseudo}$.
	For double categories $\cat{C}$	and $\cat{D}$ the \define{Gray swap functor}
	$\swap{\cat{C}}{\cat{D}} : \hom{(\tensor{\cat{C} , \cat{D}})_{\tuple*{a , p}}}{(\tensor{\cat{D} , \cat{C}})}_{\tuple*{\tilde{a} , \tilde{p}}}$
	transposes the factors of each ordered pair, sending $\tuple*{x , y}$ to $\tuple*{y , x}$.
\end{definition}

For $i + j = k$, a swap functor sends $\tuple*{i , j}$-type $k$-cells to $\tuple*{j , i}$-type $k$-cells.
Moreover, it swaps the heterogeneous $\tuple*{1 , 1}$-type $2$-cells $\tuple*{v , h}$ and $\tuple*{h , v}$,
and sends oplax homogeneous $\tuple*{1 , 1}$-type $2$-cells $\tuple*{v , v}$ and $\tuple*{h , h}$ to lax ones, and vice-versa.
Swapping is an involution, in the sense that
\begin{equation} \label{swap functor involution}
\comp{\swap{\cat{C}}{\cat{D}} , \swap{\cat{D}}{\cat{C}}}
\; = \;
\comp{} ((\tensor{\cat{C} , \cat{D}})_{\tuple*{a , p}})
.
\end{equation}
We represent swapping in surface diagrams as permuting the order of the surfaces.
For example, the mapping of an oplax $\tuple*{h , h}$-type $2$-cell to a lax one by swapping is depicted as follows.
$$

$$

\begin{definition}[braiding]
	Let $\cat{C}$ be a \refer{Gray-monoidal double category} with invertible interchangers in both dimensions.
	A(n arrow-dimension) \define{braiding} on $\cat{C}$
	is an \refer{arrow-dimension pseudo transformation}
	$$
		\braiding : \hom[\hom[]{\tensor{\cat{C} , \cat{C}}}{\cat{C}}]
		{\tensor{,}_{\cat{C}}}
		{\comp{\swap{\cat{C}}{\cat{C}} , \tensor{,}_{\cat{C}}}}
	$$
	satisfying the following relations.
	\begin{description}
		\item[nullary tensor braiding coherence:]
			Braiding with the tensor unit is trivial in the sense that
			\begin{equation} \label{equation: nullary tensor braiding}
				\braid{A}{\tensor{}}
				\; = \;
				\comp{} \, A
				\; = \;
				\braid{\tensor{}}{A}
			\end{equation}
			and these equations constitute the object-components of \refer[identity modification]{identity modifications}
			\begin{equation} \label{equation: nullary tensor braiding naturality}
				\prohom{\braid{\arg}{\tensor{}}}{\comp{} \, \arg}
				\quad \text{and} \quad
				\prohom{\braid{\tensor{}}{\arg}}{\comp{} \, \arg}
			\end{equation}
		\item[binary tensor braiding coherence:]
			Braiding with a tensor product is given by consecutive braidings in the sense that
			\begin{equation} \label{equation: binary tensor braiding}
				\begin{array}{r c l}
					\braid{A}{\tensor{B , C}}
					& =  &
					\comp{(\tensor{\braid{A}{B} , C}) , (\tensor{B , \braid{A}{C}})}
					\\
					\braid{\tensor{A , B}}{C}
					& =  &
					\comp{(\tensor{A , \braid{B}{C}}) , (\tensor{\braid{A}{C} , B})}
					\\
				\end{array}
			\end{equation}
			and these equations constitute the object-components of \refer[identity modification]{identity modifications}
			\begin{equation} \label{equation: binary tensor braiding naturality}
				\begin{array}{l}
					\prohom
					{\braid{\arg[1]}{\tensor{\arg[2] , \arg[3]}}}
					{\comp{(\tensor{\braid{\arg[1]}{\arg[2]} , \arg[3]}) , (\tensor{\arg[2] , \braid{\arg[1]}{\arg[3]}})}}
					\quad \text{and}
					\\
					\prohom
					{\braid{\tensor{\arg[1] , \arg[2]}}{\arg[3]}}
					{\comp{(\tensor{\arg[1] , \braid{\arg[2]}{\arg[3]}}) , (\tensor{\braid{\arg[1]}{\arg[3]} , \arg[2]})}}
					\\
				\end{array}
			\end{equation}
		\item[Yang-Baxterator braiding coherence:]
			Reversing the order of three objects using three braidings is coherent in the sense that
			\begin{equation} \label{equation: Yang-Baxerator braiding coherence}
				\braid{\braid{A}{B}}{C}
				\; = \;
				\inv{\braid{A}{\braid{B}{C}}}
			\end{equation}
			where, recalling our convention, $\braid{\braid{A}{B}}{C}$ is an arrow-component disk with oplax variance
			and $\inv{\braid{A}{\braid{B}{C}}}$ is one with lax variance.
			We refer to these as ``Yang-Baxterators'', and describe them in detail below.
		\item[quaternary reversal braiding coherence:]
			Reversing the order of four objects using six braidings is coherent
			as described in equation \eqref{equation: quaternary reversal braiding coherence} below.
	\end{description}
	A Gray-monoidal double category equipped with a braiding is a \define{braided Gray-monoidal double category}.
\end{definition}

We now unfold this definition and introduce the corresponding graphical notation.
First we recall what it means for $\braiding$ to be a \refer{arrow-dimension oplax transformation}.

\paragraph{Braiding \refer[arrow transformation object-component arrow]{object-component arrows}}
For objects $A , X : \cat{C}$ we have a braiding component arrow $\braid{A}{X} : \hom[\cat{C}]{\tensor{A , X}}{\tensor{X , A}}$.
We represent this in a string diagram in the $\tuple*{\tensor{,} , \comp{,}}$-plane as a crossing of the lines representing $A$ and $X$.
$$
	\begin{tikzpicture}[string diagram , x = {(16mm , 0mm)} , y = {(0mm , -12mm)} , baseline=(align.base)]
		\coordinate (front in) at (1/4 , 0) ;
		\coordinate (front out) at (1/4 , 1) ;
		\coordinate (back in) at (3/4 , 0) ;
		\coordinate (back out) at (3/4 , 1) ;
		\draw [name path = A] (front in) to [out = south , in = north] node [fromabove] {A} node [tobelow] {A} (back out) ;
		\draw [name path = X] (back in) to [out = south , in = north] node [fromabove] {X} node [tobelow] {X} (front out) ;
		\path [name intersections = {of = A and X , by = {cell}}] ;
		\node (align) at (cell) {} ;
	\end{tikzpicture}
$$

\paragraph{Braiding \refer[arrow transformation proarrow-component square]{proarrow-component squares}}
For a proarrow $m : \prohom{A}{B}$ and object $X$
we have a $\tuple*{1 , 0}$-type proarrow $\tuple*{m , X} : \prohom{\tuple*{A , X}}{\tuple*{B ,  X}}$
and a $\tuple*{0 , 1}$-type proarrow $\tuple*{X , m} : \prohom{\tuple*{X , A}}{\tuple*{X ,  B}}$.
These give us the following component squares, which we depict as shown.
$$
	\braid{m}{X} : \doublehom{\braid{A}{X}}{\braid{B}{X}}{\tensor{m , X}}{\tensor{X , m}}
	\qquad \text{and} \qquad
	\braid{X}{m} : \doublehom{\braid{X}{A}}{\braid{X}{B}}{\tensor{X , m}}{\tensor{m , X}}
$$
\begin{equation} \label{diagram: braiding proarrow component square}

\end{equation}

\paragraph{Braiding \refer[arrow transformation arrow-component disk]{arrow-component disks}}
For an arrow $f : \hom{A}{B}$ and object $X$
we have a $\tuple*{1 , 0}$-type arrow $\tuple*{f , X} : \hom{\tuple*{A , X}}{\tuple*{B ,  X}}$
and a $\tuple*{0 , 1}$-type arrow $\tuple*{X , f} : \hom{\tuple*{X , A}}{\tuple*{X ,  B}}$.
These give us the following component disks, which we depict as shown.
$$
	\begin{array}{l}
		\braid{f}{X} : \prohom{\comp{(\tensor{f , X}) , \braid{B}{X}}}{\comp{\braid{A}{X} , (\tensor{X , f})}}
		\quad \text{and} \quad
		\\
		\braid{X}{f} : \prohom{\comp{(\tensor{X , f}) , \braid{X}{B}}}{\comp{\braid{X}{A} , (\tensor{f , X})}}
		\\
	\end{array}
$$
\begin{equation} \label{diagram: braiding arrow component disk}
	%
\end{equation}

\paragraph{Braiding composite $1$-cells}
The nullary clause of preservation of proarrow composition \eqref{arrow transformation proarrow composition preservation} and
the nullary clause of compatibility with arrow composition \eqref{arrow transformation arrow composition compatibility}
together with equation \eqref{pair nullary proarrow composite} say that
for objects $A$ and $X$ we have the following equations.
\begin{equation} \label{equation: braiding nullary composites}
	\braid{\procomp{} A}{X}
	=
	\braid{A}{\procomp{} X}
	=
	\procomp{} (\braid{A}{X})
	=
	\braid{\comp{} A}{X}
	=
	\braid{A}{\comp{} X}
\end{equation}
This provides a unique interpretation to the following surface diagram.
$$
	\begin{tikzpicture}[string diagram , x = {(16mm , -8mm)} , y = {(0mm , -24mm)} , z = {(12mm , 0mm)} , baseline=(align.base)]
		\coordinate (sheet) at (1/2 , 1/2 , 1) ;
		\coordinate (left top back) at ($ (sheet) + (-1/2 , -1/2 , 1/2) $) ;
		\coordinate (right top back) at ($ (sheet) + (1/2 , -1/2 , 1/2) $) ;
		\coordinate (left top front) at ($ (sheet) + (-1/2 , -1/2 , -1/2) $) ;
		\coordinate (right top front) at ($ (sheet) + (1/2 , -1/2 , -1/2) $) ;
		\coordinate (left seam) at ($ (sheet) + (-1/2 , 0 , 0) $) ;
		\coordinate (right seam) at ($ (sheet) + (1/2 , 0 , 0) $) ;
		\coordinate (left bot back) at ($ (sheet) + (-1/2 , 1/2 , 1/2) $) ;
		\coordinate (right bot back) at ($ (sheet) + (1/2 , 1/2 , 1/2) $) ;
		\coordinate (left bot front) at ($ (sheet) + (-1/2 , 1/2 , -1/2) $) ;
		\coordinate (right bot front) at ($ (sheet) + (1/2 , 1/2 , -1/2) $) ;
		\node [anchor = east] at (left seam) {\braid{A}{X}} ;
		\draw [sheet]
			(left seam) to
			(right seam) .. controls +(0 , 1/8 , 1/2) and +(0 , -1/8 , 0) ..
			(right bot back) to
			(left bot back) .. controls +(0 , -1/8 , 0) and +(0 , 1/8 , 1/2) ..
			cycle ;
		\draw [sheet]
			(left top back) to
			(right top back) .. controls +(0 , 1/8 , 0) and +(0 , -1/8 , 1/2) ..
			(right seam) to
			(left seam) .. controls +(0 , -1/8 , 1/2) and +(0 , 1/8 , 0) ..
			cycle ;
		\draw [sheet]
			(left seam) to
			(right seam) .. controls +(0 , 1/8 , -1/2) and +(0 , -1/8 , 0) ..
			(right bot front) to
			(left bot front) .. controls +(0 , -1/8 , 0) and +(0 , 1/8 , -1/2) ..
			cycle ;
		\draw [sheet]
			(left top front) to
			(right top front) .. controls +(0 , 1/8 , 0) and +(0 , -1/8 , -1/2) ..
			(right seam) to
			(left seam) .. controls +(0 , -1/8 , -1/2) and +(0 , 1/8 , 0) ..
			cycle ;
		\node [anchor = west] at (right seam) {\braid{A}{X}} ;
		\node (align) at (sheet) {} ;
	\end{tikzpicture}
$$
Similarly, the binary clauses
of \eqref{arrow transformation proarrow composition preservation}
and \eqref{arrow transformation arrow composition compatibility}
provide unique interpretations for all instances of this diagram
obtained by ``painting'' onto it two consecutive arrows or proarrows,
each independently of type $\tuple*{1 , 0}$ or $\tuple*{0 , 1}$,
with non-intersecting $\tensor{,}$-dimension projections.
For example, given arrows $f : \hom{A}{B}$ and $i : \hom{X}{Y}$
the equation \eqref{arrow transformation arrow composition compatibility} instance
$$
	\braiding (\comp{\tuple*{f , X} , \tuple*{B , i}})
	=
	\procomp{(\comp{\procomp{} (\tensor{f , X}) , \braid{B}{i}}) , (\comp{\braid{f}{X} , \procomp{} (\tensor{i , B})})}
$$
provides a unique interpretation to the following surface diagram on the left.
Likewise, for arrow $g : \hom{B}{C}$
equations \eqref{arrow transformation arrow composition compatibility} and \eqref{pair binary arrow composite} together imply
\begin{equation} \label{equation: braiding homogeneous composite arrow}
	\braid{\comp{f , g}}{X}
	=
	\procomp{(\comp{\procomp{} (\tensor{f , X}) , \braid{g}{X}}) , (\comp{\braid{f}{X} , \procomp{} (\tensor{X , g})})}
\end{equation}
providing a unique interpretation to the one on the right.
$$

$$

\paragraph{Braiding naturality for squares}
For square $φ : \doublehom{f}{g}{m}{n}$
we have a $\tuple*{2 , 0}$-type square $\tuple*{φ , X}$ and a $\tuple*{0 , 2}$-type square $\tuple*{X , φ}$.
In this case naturality for squares \eqref{arrow transformation square naturality} says
\begin{equation} \label{braiding naturality for squares (2,0)}
	\begin{array}{rcl}
		\procomp{(\comp{(\tensor{φ , X}) , \braid{n}{X}}) , \braid{g}{X}}
		& ≅ &
		\procomp{\braid{f}{X} , (\comp{\braid{m}{X} , (\tensor{X , φ})})}
		\\
		\procomp{(\comp{(\tensor{X , φ}) , \braid{X}{n}}) , \braid{X}{g}}
		& ≅ &
		\procomp{\braid{X}{f} , (\comp{\braid{X}{m} , (\tensor{φ , X})})}
		\\
	\end{array}
\end{equation}
This identifies the two boundary-preserving perturbations of each of the following surface diagrams
that move the square $φ$ up or down away from the braiding, as shown.
$$
	%
$$

For arrow $f : \hom{A}{B}$ and proarrow $p : \prohom{X}{Y}$
we have a $\tuple*{v , h}$-type square $\tuple*{f , p}$
and a  $\tuple*{h , v}$-type square $\tuple*{p , f}$.
In this case naturality for squares \eqref{arrow transformation square naturality} says
$$
	\begin{array}{rcl}
		\procomp{(\comp{\interchange{f , p} , \braid{B}{p}}) , \braid{f}{Y}}
		& ≅ &
		\procomp{\braid{f}{X} , (\comp{\braid{A}{p} , \interchange{p , f}})}
		\\
		\procomp{(\comp{\interchange{p , f} , \braid{p}{B}}) , \braid{Y}{f}}
		& ≅ &
		\procomp{\braid{X}{f} , (\comp{\braid{p}{A} , \interchange{f , p}})}
		\\
	\end{array}
$$
This identifies the two boundary-preserving perturbations of each of the following surface diagrams
that move the interchanger of $f$ and $p$ up or down away from the braiding.
$$
	%
$$

For proarrows $m : \prohom{A}{B}$ and $p : \prohom{X}{Y}$
we have an $\tuple*{h , h}$-type oplax disk $\tuple*{m , p}$
and an $\tuple*{h , h}$-type lax disk $\tuple*{p , m}$.
In this case the instance of naturality for squares in equation \eqref{arrow transformation naturality for proarrow disks} says
$$
	\comp{\interchange{m , p} , (\procomp{\braid{A}{p} , \braid{m}{Y}})}
	\quad = \quad
	\comp{(\procomp{\braid{m}{X} , \braid{B}{p}}) , \inv{\interchange{p , m}}}
$$
where, according to our convention on variance,
$\interchange{m , p}$ is an interchanger disk with oplax orientation and $\inv{\interchange{p , m}}$ is one with lax orientation.
This identifies the two boundary-preserving perturbations of the following surface diagram
that move the interchanger of $m$ and $p$ up or down away from the braiding.
$$

$$
Note that we need a \refer{braided Gray-monoidal double category} to have invertible interchangers for proarrows
because the \refer{Gray swap functor} inverts the interchanger variance below the \refer{braiding}.

For arrows $f : \hom{A}{B}$ and $i : \hom{X}{Y}$
we have a $\tuple*{v , v}$-type oplax disk $\tuple*{f , i}$
and a $\tuple*{v , v}$-type lax disk $\tuple*{i , f}$.
In this case the instance of naturality for squares in equation \eqref{arrow transformation naturality for arrow disks} says
$$
	\begin{array}{c}
		\procomp
		{
			[\comp{\interchange{f , i} , \procomp{} (\braid{B}{Y})}] ,
			[\comp{\procomp{} (\tensor{A , i}) , \braid{f}{Y}}] ,
			[\comp{\braid{A}{i} , \procomp{} (\tensor{Y , f})}]
		}
		\\ = \\
		\procomp
		{
			[\comp{\procomp{} (\tensor{f , X}) , \braid{B}{i}}] ,
			[\comp{\braid{f}{X} , \procomp{} (\tensor{i , B})}] ,
			[\comp{\procomp{} (\braid{A}{X}) , \inv{\interchange{i , f}}}]
		}
	\end{array}
$$
This identifies the two boundary-preserving perturbations of the following surface diagram
that move the interchanger of $f$ and $i$ up or down away from the braiding.
$$
	%
$$
Note that we need a \refer{braided Gray-monoidal double category} to have invertible interchangers for arrows
because the \refer{Gray swap functor} inverts the interchanger variance below the \refer{braiding}.

\paragraph{Braiding pseudo naturality}
The requirement that $\braiding$ be an \refer{arrow-dimension pseudo transformation}
means that we also have arrow-component disks with the lax variance
$$
	\begin{array}{l}
		\inv{\braid{f}{X}} : \prohom{\comp{\braid{A}{X} , (\tensor{X , f})}}{\comp{(\tensor{f , X}) , \braid{B}{X}}}
		\quad \text{and} \quad
		\\
		\inv{\braid{X}{f}} : \prohom{\comp{\braid{X}{A} , (\tensor{f , X})}}{\comp{(\tensor{X , f}) , \braid{X}{B}}}
		\\
	\end{array}%
$$
which are proarrow-dimension inverse \eqref{equation: pseudo interchanger laws}
to $\braid{f}{X}$ and $\braid{X}{f}$ respectively.
For example, the relation $\procomp{\braid{f}{X} , \inv{\braid{f}{X}}} = \procomp{} (\comp{(\tensor{f , X}) , \braid{B}{X}})$
is depicted below.
$$
	%
$$

\paragraph{$\tuple*{m, n}$-braiding coherence}
The nullary and binary tensor braiding laws \eqref{equation: nullary tensor braiding} and \eqref{equation: binary tensor braiding}
assert that each of the following string diagrams in the $\tuple*{\tensor{,} , \comp{,}}$-plane represents a unique arrow.
\begin{equation} \label{braiding coherence diagrams}
	%
\end{equation}
The braiding of an $m$-ary tensor of objects with an $n$-ary tensor of objects is well defined
up to interchangers of the intermediate (i.e., neither first nor last) braidings.
For example, when $m$ and $n$ are both $2$ we have an interchanger of braidings,
$$
	%
$$
Because a \refer{braided Gray-monoidal double category} has invertible interchangers
we will be opportunistically ambiguous about the precise order of the intermediate braidings
and refer to a braiding of an $m$-ary tensor with an $n$-ary tensor
as an ``$\tuple*{m , n}$-braiding''.

The requirement that equations \eqref{equation: nullary tensor braiding} and \eqref{equation: binary tensor braiding}
be the object-component disks of \refer[identity modification]{identity modifications}
\eqref{equation: nullary tensor braiding naturality} and \eqref{equation: binary tensor braiding naturality}
is a succinct way to make them natural
for arrows \eqref{arrow disk modification naturality for arrows} and
for proarrows \eqref{arrow disk modification naturality for proarrows} in all indices.
This implies that if we ``paint'' an arrow or proarrow onto any surface of one on the following diagrams,
the result has a unique interpretation.
$$

$$

For an arrow $f : \hom{A}{B}$ and objects $X$ and $Y$, 
in the case of a $\tuple*{2 , 1}$-braiding we have the following equations.
\begin{equation} \label{equation: (2,1)-braiding arrow coherence}
	\begin{array}{l}
		\braid{\tensor{f , X}}{Y}
		=
		\procomp
		{
			[\comp{\interchange{f , \braid{X}{Y}} , \procomp{} (\tensor{\braid{B}{Y} , X})}] ,
			[\comp{\procomp{} (\tensor{A , \braid{X}{Y}}) , (\tensor{\braid{f}{Y} , X})}]
		}
		\\
		\braid{\tensor{X , f}}{Y}
		=
		\procomp
		{
			[\comp{(\tensor{X , \braid{f}{Y}}) , \procomp{} (\tensor{\braid{X}{Y} , B})}] ,
			[\comp{\procomp{} (\tensor{X , \braid{A}{Y}}) , \inv{\interchange{\braid{X}{Y} , f}}}]
		}
		\\
		\braid{\tensor{X , Y}}{f}
		=
		\procomp
		{
			[\comp{(\tensor{X , \braid{Y}{f}}) , \procomp{} (\tensor{\braid{X}{B} , Y})}] ,
			[\comp{\procomp{} (\tensor{X , \braid{Y}{A}}) , (\tensor{\braid{X}{f} , Y})}]
		}
		\\
	\end{array}
\end{equation}

These provide unique interpretations for the following surface diagrams.
$$
	%
$$

For example, in the last case
by equation \eqref{equation: binary tensor braiding} we have the following identity disk in the $\tuple*{\procomp{,} , \comp{,}}$-plane.
\begin{equation} \label{temp1}
	%
\end{equation}
So for arrow $f : \hom{A}{B}$ we have
$$
	%
	\; \overset{\eqref{arrow disk modification naturality for arrows}}{=} \;
	%
$$
Of course, we have the corresponding $\tuple*{1 , 2}$-braiding laws for arrows as well.

For a proarrow $m : \prohom{A}{B}$ and objects $X$ and $Y$, 
in the case of a $\tuple*{1 , 2}$-braiding we have the following equations.
\begin{equation} \label{equation: (2,1)-braiding proarrow coherence}
	\begin{array}{l}
		\braid{m}{\tensor{X , Y}}
		=
		\comp
		{
			(\tensor{\braid{m}{X} , Y}) ,
			(\tensor{X , \braid{m}{Y}})
		}
		\\
		\braid{X}{\tensor{m , Y}}
		=
		\comp
		{
			(\tensor{\braid{X}{m} , Y}) ,
			\interchange{m , \braid{X}{Y}}
		}
		\\
		\braid{X}{\tensor{Y , m}}
		=
		\comp
		{
			\interchange{\braid{X}{Y} , m} ,
			(\tensor{Y , \braid{X}{m}})
		}
		\\
	\end{array}
\end{equation}
These provide unique interpretations for the following surface diagrams.
$$
	%
$$
Of course, we have the corresponding $\tuple*{2 , 1}$-braiding laws for proarrows as well.

\paragraph{Yang-Baxterator braiding coherence}
The braiding transformation $\braiding$ has object-component arrows,
and therefore arrow-component disks for these arrows,
$$
	\begin{array}{c}
		\braid{\braid{A}{B}}{C} :
		\prohom
		{\comp{(\tensor{\braid{A}{B} , C}) ,  \braid{\tensor{B , A}}{C}}}
		{\comp{\braid{\tensor{A , B}}{C} , (\tensor{C , \braid{A}{B}})}}
		\\ \text{and} \\
		\braid{A}{\braid{B}{C}} :
		\prohom
		{\comp{(\tensor{A , \braid{B}{C}}) , \braid{A}{\tensor{C , B}}}}
		{\comp{\braid{A}{\tensor{B , C}} , (\tensor{\braid{B}{C} , A})}}
	\end{array}%
$$

By equation \eqref{equation: binary tensor braiding} these both have the following proarrow-dimension boundaries,
left-to-right in the first case and right-to-left in the second.
$$
	%
\end{equation}
We can think of these arrow disks as \define[Yang-Baxterator]{Yang-Baxterators},
in the sense that they are \emph{structure} categorifying the \emph{property} of the Yang-Baxter equation.

Because $\braiding$ is a pseudo transformation
we also have their respective proarrow-dimension inverses,
the lax-oriented arrow-component disks
$$
	\begin{array}{c}
		\inv{\braid{\braid{A}{B}}{C}} :
		\prohom
		{\comp{\braid{\tensor{A , B}}{C} , (\tensor{C , \braid{A}{B}})}}
		{\comp{(\tensor{\braid{A}{B} , C}) ,  \braid{\tensor{B , A}}{C}}}
		\\ \text{and} \\
		\inv{\braid{A}{\braid{B}{C}}} :
		\prohom
		{\comp{\braid{A}{\tensor{B , C}} , (\tensor{\braid{B}{C} , A})}}
		{\comp{(\tensor{A , \braid{B}{C}}) , \braid{A}{\tensor{C , B}}}}
	\end{array}%
$$
By equation \eqref{equation: binary tensor braiding} these also have the proarrow-dimension boundary
depicted in diagram \eqref{yang-baxterator boundary},
right-to-left in the first case and left-to-right in the second.

Yang-Baxterator braiding coherence \eqref{equation: Yang-Baxerator braiding coherence} asserts that
$\inv{\braid{\braid{A}{B}}{C}}$ and $\inv{\braid{A}{\braid{B}{C}}}$
are respectively
$\braid{A}{\braid{B}{C}}$ and $\braid{\braid{A}{B}}{C}$;
that is, that
$\braid{\braid{A}{B}}{C}$ and $\braid{A}{\braid{B}{C}}$
are proarrow-dimension inverse to one another.
In light of this coherence, we will refer to the left-to-right orientation of \eqref{yang-baxterator boundary} as $\yangbaxt{A}{B}{C}$,
and to the right-to-left one as $\inv{\yangbaxt{A}{B}{C}}$.

\paragraph{Quaternary reversal braiding coherence}
For the sake of clarity and concision, in the following description we suppress interchangers of braiding object-component arrows.
Recall that it takes $n (n−1) / 2$-many transpositions to reverse a list of $n$ elements.
Thus, this is the minimum number of $\tuple*{1 , 1}$-braidings required to reverse the order of $n$ objects.
When $n$ is $4$ there are eight possible ways to do this;
that is, (quotiented by interchangers) there are eight terms $\hom{\tensor{A , B , C , D}}{\tensor{D , C , B , A}}$,
each composed of six whiskered $\tuple*{1 , 1}$-braidings.
These are related to one another by composites of whiskered Yang-Baxterators.
Quaternary reversal braiding coherence identifies parallel pairs of these.

Even suppressing the interchangers, an algebraic presentation of this law is quite verbose and imperspicuous.
The surface diagram representation, which is equivalent to figure 40 in \cite{carter-1997-knotted_surfaces},
is a bit better, but is still difficult to parse from a fixed perspective.
Instead, we present this law as an equation between sequences of \emph{rewrites}
of ``slice'' string diagrams in the $\tuple*{\tensor{,} , \comp{,}}$-plane.
As is customary in rewriting, we label the transitions for only the parts of the diagram that are changed.
This reduces notational clutter by suppressing identity morphisms.
Each of the following rewrite sequences represents an implied surface diagram, and the law identifies these.
$$

$$
\begin{equation} \label{equation: quaternary reversal braiding coherence}
=
\end{equation}
$$
	%
$$
This relation is sometimes called ``tetrahedron'' coherence
\cite{voevodsky-1994-2categories, baez-1996-braided_monoidal_2categories},
presumably because in each of the implied surface diagrams
the four Yang-Baxterators form the vertices of a tetrahedron\footnote
{I am grateful to John Baez for explaining this geometric interpretation to me.}.

\vspace{1em}

A natural question to ask is what happens if we have two consecutive braidings of the same pair of objects.
In a monoidal $1$-category the only thing we can do is impose a relation, known as \emph{symmetry}.
But as with the Yang-Baxterators, given more dimensions we can categorify this property to a structure.

\begin{definition}[syllepsis]
	A \define{syllepsis} for a \refer{braided Gray-monoidal double category} $\cat{C}$
	is an \refer[invertible modification]{invertible globular} \refer{modification}
	$$
		\syllepsis : \prohom[\hom[]{\tensor{,}_{\cat{C}}}{\tensor{,}_{\cat{C}}}]
		{\comp{} (\tensor{,}_{\cat{C}})}{\comp{\braiding , (\comp[2]{\swapper , \braiding})}}
	$$
	That is, for objects $A$ and $B$, a syllepsis $\syllepsis$ has \refer[modification object-component square]{object-component disk}
	$$
		\syllepsize{A}{B} : \prohom{\comp{} (\tensor{A , B})}{\comp{\braid{A}{B} , \braid{B}{A}}}
	$$
	It is required to satisfy the following relations.
	\begin{description}
		\item[nullary tensor syllepsis coherence:]
			Syllepsis with the tensor unit is trivial in the sense that
			\begin{equation} \label{equation: nullary tensor syllepsis}
				\syllepsize{A}{\tensor{}}
				\; = \;
				\comp[2]{} \, A
				\; = \;
				\syllepsize{\tensor{}}{A}
			\end{equation}
		\item[binary tensor syllepsis coherence:]
			Syllepsis with a tensor product is given by nested syllepses in the sense that
			\begin{equation} \label{equation: binary tensor syllepsis}
\begin{array}{r @{\hskip 4pt} c @{\hskip 4pt} l}
					\syllepsize{A}{\tensor{B , C}}
					& =  &
					\procomp
					{
						(\tensor{\syllepsize{A}{B} , C}) ,
						[\comp{\procomp{} (\tensor{\braid{A}{B} , C}) , (\tensor{B , \syllepsize{A}{C}}) , \procomp{} (\tensor{\braid{B}{A} , C})}]
					}
					\\
					\syllepsize{\tensor{A , B}}{C}
					& =  &
					\procomp
					{
						(\tensor{A , \syllepsize{B}{C}}) ,
						[\comp{\procomp{} (\tensor{A , \braid{B}{C}}) , (\tensor{\syllepsize{A}{C} , B}) , \procomp{} (\tensor{A , \braid{C}{B}})}]
					}
					\\
				\end{array}
			\end{equation}
		\item[Yang-Baxterator syllepsis coherence:]
			Syllepsis is compatible with Yang-Baxterators,
			as described in equation \eqref{equation: yang-baxterator syllepsis coherence} below.
	\end{description}
	A \refer{braided Gray-monoidal double category} with a syllepsis is a \define{sylleptic Gray-monoidal double category}.
\end{definition}

A syllepsis relates the ``un-braiding'' to a consecutive pair of braidings.
In surface diagrams we draw a syllepsis object-component disk $\syllepsize{A}{B}$ as shown.
$$
	%
$$

The nullary tensor syllepsis coherence laws \eqref{equation: nullary tensor syllepsis}
are dependent upon the nullary tensor braiding coherence laws \eqref{equation: nullary tensor braiding},
and ensure that the following surface diagram continues to have a unique interpretation in the presence of sylleptic structure.
$$
	%
$$
The binary tensor syllepsis coherence laws \eqref{equation: binary tensor syllepsis}
are dependent upon the binary tensor braiding coherence laws \eqref{equation: binary tensor braiding},
and provide unique interpretations to each of the following surface diagrams.
$$
	%
$$
To better understand these it may help to render the proarrow disks
on the right of equations \eqref{equation: binary tensor syllepsis} as rewrites.
$$
	%
$$

Yang-Baxterator syllepsis coherence says that the following cycle of invertible proarrow disks, along with its mirror image, is the identity.
\begin{equation} \label{equation: yang-baxterator syllepsis coherence}
	%
\end{equation}

Modification naturality for arrows \eqref{arrow disk modification naturality for arrows}
and for proarrows \eqref{arrow disk modification naturality for proarrows}
ensure that if we ``paint'' an arrow or proarrow onto any surface in a diagram involving a syllepsis
it will pass across the syllepsis as suggested by the topology.

\begin{definition}[symmetry]
	A syllepsis $\syllepsis$ is called a \define{symmetry}
	if $\tuple*{\braiding , \comp[2]{\swapper , \braiding} , \syllepsis , \comp[3]{\swapper , \inv{\syllepsis}}}$
	form an adjoint equivalence;
	that is, if for objects $A$ and $B$ we have an adjunction $\adjoint{\braid{A}{B} , \braid{B}{A}}$
	with unit $\syllepsize{A}{B}$ and counit $\inv{\syllepsize{B}{A}}$.
	
	A \refer{sylleptic Gray-monoidal double category} with a symmetry is a \define{symmetric Gray-monoidal double category}.
\end{definition}

We represent the adjunction laws in surface diagrams as follows.
$$

$$

In the globular literature symmetry is characterized by the following equation
\cite{street-1997-monoidal_bicategories}.
\begin{equation} \label{equation: syllepsis symmetry law}
	\comp{\syllepsize{A}{B} , \procomp{} (\braid{A}{B})}
	\quad = \quad
	\comp{\procomp{} (\braid{A}{B}) , \syllepsize{B}{A}}
\end{equation}
$$
	%
$$
Our definition is equivalent.

\begin{proposition}
	A syllepsis is a symmetry just in case equation \eqref{equation: syllepsis symmetry law} holds.
\end{proposition}

\begin{proof}
	Assuming the adjoint equivalence, we have
	$$
	%
	$$
	where in the first step we rewrite the identity arrow disk
	$\procomp{} (\comp{\braid{B}{A} , \braid{A}{B}})$
	to the composite of inverses
	$\procomp{\inv{\syllepsize{B}{A}} , \syllepsize{B}{A}}$,
	and in the second we apply the left adjunction law.
	
	Conversely, assuming equation \eqref{equation: syllepsis symmetry law}, we have
	$$
	%
	$$
	where in the first step we rewrite by equation \eqref{equation: syllepsis symmetry law} forwards,
	and in the second we rewrite the composite of inverses
	$\procomp{\syllepsize{B}{A} , \inv{\syllepsize{B}{A}}}$ to the identity square
	$\comp[2]{} (\tensor{B , A})$.
	Similarly, rewriting by equation \eqref{equation: syllepsis symmetry law} backwards
	we get the other adjunction law.
	$$
	%
	$$
\end{proof}

\paragraph{Related Constructions}
The development of symmetric Gray-monoidal double categories presented here is inspired
by Böhm's \emph{double category analogue of Gray monoids}
\cite{bohm-2019-monoidal_product_of_double_categories}, 
by Garner and Gurski's presentation of \emph{monoidal double categories}
\cite{garner-2009-tricategory_structures}, 
and by Shulman's development of \emph{symmetric monoidal double categories}
\cite{shulman-2010-symmetric_monoidal_bicategories}. 

Another thread of influence, coming from the globular perspective,
is the idea of a braided \define[Gray-monoidal 2-category]{Gray-monoidal $2$-category}
(or ``Gray monoid'' or ``semistrict monoidal $2$-category''),
which was proposed by Kapranov and Voevodsky \cite{voevodsky-1994-2categories}.
This was further developed by Baez and Neuchl \cite{baez-1996-braided_monoidal_2categories},
who added Yang-Baxterator braiding coherence,
by Day and Street \cite{street-1997-monoidal_bicategories},
who added syllepsis and symmetry,
and by Crans \cite{crans-1998-sylleptic_monoidal_2categories},
who added nullary tensor braiding coherence.

Our Gray-monoidal double categories are closest to Böhm's double categorical Gray monoids.
Our construction differs in being derived from locally cubical Gray categories,
and thus being parametric in the variance of the homogeneous interchangers.
We also go on to define braided, sylleptic, and symmetric structure.

Garner and Gurski define monoidal double categories
as one-object \refer{locally cubical bicategories} (section \ref{section: locally cubical gray categories}),
making their monoidal structure weak, but not Gray.
This is different from our construction
in that the monoidal structure on a double category $\cat{C}$
is given by a functor from the \refer{cartesian ordered pair double category}
$\hom{\cat{\product{\cat{C} , \cat{C}}}}{\cat{C}}$,
whereas for us it is given by a functor from the \refer{Gray ordered pair double category}
$\hom{\cat{\tensor{\cat{C} , \cat{C}}}}{\cat{C}}$.
Intuitively, in their setting you can do multiple (well-bracketed) things at once,
whereas in ours you can do only one thing at a time,
though you may change the order using suitably oriented interchangers if those things are independent.

Shulman gives a direct presentation by generators and relations of monoidal double categories,
and goes on to define symmetric braided structure for them. 
Shulman's braiding is a strict arrow-dimension transformation, 
while ours is pseudo.
We defend our choice on topological grounds:
the boundary arrows of a braiding arrow-component disk, depicted in diagram \eqref{diagram: braiding arrow component disk},
are topologically distinct, so we prefer not to identify them.

Our construction differs from those on Gray-monoidal $2$-categories mentioned above
most notably in that the two-dimensional categories involved are cubical rather than globular.
Another important difference is that in the definition of braiding
we impose the tensor braiding coherence equations \eqref{equation: binary tensor braiding}.
By doing so we are, in a sense, leaving food on the table
because these equations strictly identify structures whose dimension is not maximal.
In the globular literature cited above the modifications
corresponding to \eqref{equation: binary tensor braiding naturality}
are only invertible, not identities.
This results in an additional layer of tensor braiding coherence, which we avoid.

While the other approach is certainly more general,
we defend our choice on the following grounds.
First, maximal weakness need not be our goal.
If we don't consider them meaningful then compounding coherators only serves to add
bureaucracy to the constructions we are trying to perform.
The choice to use a strictly associative and unital Gray tensor product already precludes fully weak constructions,
so we are not the only ones leaving food on the table.

More affirmatively, our choice is guided by topology.
The diagrammatic representations of the terms in each of equations \eqref{equation: binary tensor braiding}
are topologically identical, as depicted in diagram \eqref{braiding coherence diagrams}.
This allows us to treat braidings compositionally.
Baez and Neuchl use a similar topological argument
to justify the addition of Yang-Baxterator coherence
to the braiding coherence laws originally proposed by Kapranov and Voevodsky,
although in that case the dimension of the structures involved is maximal so no food is wasted.

%% file: content/parts/cartesian_structure.tex

A presentation of cartesian products in $1$-categories
is typically given by the universal construction of a terminal cone over a functor from a discrete category.
Such a presentation is ``behavioral'' in the sense that it distinguishes product cones by a universal property in relation to all possible cones.
Of course, this presentation is equivalent to one
determining an isomorphism of $\Cat{Set}$-products of homs and homs to product objects.
Sadly, when we move away from the setting of $1$-categories,
such happy coincidences don't always obtain
\cite{moser-2022-2_limits}.

An alternative presentation is that of Fox,
who observed that in the case of symmetric monoidal $1$-categories, 
cartesian structure can be characterized by a pair of natural transformations
that are compatible with the monoidal structure and
whose components endow objects with the structure of cocommutative comonoids
\cite{fox-1976-cartesian_categories}.
The comonoid comultiplication acts as a \emph{duplicator} and its counit as a \emph{deletor}. 
This ``structural'' characterization of cartesian products has the benefit
that it can be given a presentation by generators and relations.

It is such a ``Fox-cartesian'' structure for Gray-monoidal double categories that we present here.
In the Gray-monoidal setting our duplicators must produce their copies in a given order,
which is permuted by the interchagers and braidings.
Of course, composition of arrows and of proarrows needs to be compatible with this ordering.
One consequence of this is that the map making two copies of a thing from one
cannot be a \refer[strict double category functor]{strict functor}
(definition \ref{definition: strict double category functor}).

A definition of lax functor for double categories
was proposed by Grandis and Paré \cite{grandis-1999-limits_in_double_categories},
but their maps are lax in only the proarrow dimension.
In order to copy composites of both arrows and proarrows,
we need morphisms of double categories that are lax in both dimensions.

\begin{definition}[lax functor of double categories]
	A (doubly) \define[lax double category functor]{lax functor} of double categories, $F : \hom{\cat{C}}{\cat{D}}$,
	consists of a collection of boundary-preserving maps of objects, arrows, proarrows, and squares,
	$$
		\begin{tikzpicture}[string diagram , x = {(12mm , 0mm)} , y = {(0mm , -12mm)} , baseline=(align.base)]
			\coordinate (c) at (1/2 , 1/2) ;
			\draw [name path = arr] (0 , 0 |- c) to node [fromleft] {f} node [toright] {g} (1 , 1 |- c) ;
			\draw [name path = pro] (c |- 0 , 0) to node [fromabove] {m} node [tobelow] {n} (c |- 1 , 1) ;
			\path [name intersections = {of = arr and pro , by = {cell}}] ;
			\node at ($ (cell) + (-1/3 , -1/3) $) {A} ;
			\node at ($ (cell) + (-1/3 , 1/3) $) {C} ;
			\node at ($ (cell) + (1/3 , -1/3) $) {B} ;
			\node at ($ (cell) + (1/3 , 1/3) $) {D} ;
			\node [bead] (align) at (cell) {α} ;
		\end{tikzpicture}
		\qquad \overset{F}{↦} \qquad
		\begin{tikzpicture}[string diagram , x = {(12mm , 0mm)} , y = {(0mm , -12mm)} , baseline=(align.base)]
			\coordinate (c) at (1/2 , 1/2) ;
			\draw [name path = arr] (0 , 0 |- c) to node [fromleft] {F f} node [toright] {F g} (1 , 1 |- c) ;
			\draw [name path = pro] (c |- 0 , 0) to node [fromabove] {F m} node [tobelow] {F n} (c |- 1 , 1) ;
			\path [name intersections = {of = arr and pro , by = {cell}}] ;
			\node at ($ (cell) + (-1/3 , -1/3) $) {F A} ;
			\node at ($ (cell) + (-1/3 , 1/3) $) {F C} ;
			\node at ($ (cell) + (1/3 , -1/3) $) {F B} ;
			\node at ($ (cell) + (1/3 , 1/3) $) {F D} ;
			\node [bead] (align) at (cell) {F α} ;
		\end{tikzpicture}
	$$
	together with a natural and coherent \define{comparison structure} in each dimension,
	$\comparitor^{\comp{,}}$ and $\comparitor^{\procomp{,}}$, as follows.
	We typically omit the dimension of a comparison structure
	as well as the arity of a comparitor disk
	when they are clear from the context.
	
	\paragraph{arrow comparison structure}
	For each object $A$
	there is a nullary comparitor arrow disk
	$\comparitor[0]^{\comp{,}} \tuple*{A} : \prohom{\comp{} (F A)}{F (\comp{} \, A)}$,
	and for consecutive arrows $f : \hom{A}{B}$ and $g : \hom{B}{C}$
	there is a binary comparitor arrow disk
	$\comparitor[2]^{\comp{,}} \tuple*{f , g} : \prohom{\comp{F f , F g}}{F (\comp{f , g})}$.
	\begin{equation} \label{lax arrow comparitor components}

	\end{equation}
	These are coherent in the sense that
	\begin{equation} \label{equation: lax arrow comparitor unit law}
		%
	\end{equation}
	and
	\begin{equation} \label{equation: lax arrow comparitor associative law}
		%
	\end{equation}
	and natural in the sense that
	\begin{equation} \label{equation: lax nullary comparitor arrow naturality}
		%
	\end{equation}
	and
	\begin{equation} \label{equation: lax binary comparitor arrow naturality}
		%
	\end{equation}
		
	\paragraph{proarrow comparison structure}
	For each object $A$
	there is a nullary comparitor proarrow disk
	$\comparitor[0]^{\procomp{,}} \tuple*{A} : \hom{\procomp{} (F A)}{F (\procomp{} \, A)}$,
	and for consecutive proarrows $m : \prohom{A}{B}$ and $n : \prohom{B}{C}$
	there is a binary comparitor proarrow disk
	$\comparitor[2]^{\procomp{,}} \tuple*{m , n} : \hom{\procomp{F m , F n}}{F (\procomp{m , n})}$.
	\begin{equation} \label{lax proarrow comparitor components}
		%
	\end{equation}
	These are coherent in the sense that
	\begin{equation} \label{equation: lax proarrow comparitor unit law}
		%
	\end{equation}
	and
	\begin{equation} \label{equation: lax proarrow comparitor associative law}
		%
	\end{equation}
	and natural in the sense that
	\begin{equation} \label{equation: lax nullary comparitor proarrow naturality}
		%
	\end{equation}
	and
	\begin{equation} \label{equation: lax binary comparitor proarrow naturality}
		%
	\end{equation}
\end{definition}

As usual, this presentation in terms of coherent nullary and binary generators is equivalent to
an \emph{unbiased} presentation in terms of coherent generators for all natural number arities \cite{leinster-2003-higher_categories}.
If the arrow- or proarrow-dimension nullary comparitors are invertible
then the lax functor is called \define[normal lax double category functor]{normal} in that dimension.
If the binary (and thus all positive arity) comparitors are invertible as well
then it is called \define[pseudo double category functor]{pseudo}.
If the nullary comparitors are identities
then it is \define[unitary lax double category functor]{unitary},
and if all comparitors in both dimensions are identities
then we are back to a \refer[strict double category functor]{strict functor of double categories}
(definition \ref{definition: strict double category functor}).

Our present interest is a particular class of unitary lax functors of double categories that will allow us to make copies.
In fact, we can define four such unitary lax functors, which differ only in the order in which the copies are  produced.
For simplicity, we discuss only the two that produce copies in the same order in both dimensions.
The following definition is the Gray-monoidal analogue
of a \refer{diagonal functor} to a \refer{cartesian ordered pair double category}
$\diag_{\cat{C}} : \hom{\cat{C}}{\product{\cat{C} , \cat{C}}}$ described in section \ref{section: double categories}.

\begin{definition}[Gray diagonal functor]
	Let $\tensor{\cat{C} , \cat{C}}$ be a \refer{Gray ordered pair double category} with oplax interchangers in both dimensions.
	We define the unitary lax \define{Gray diagonal functor}
	$\diag_{\cat{C}} : \hom{\cat{C}}{\tensor{\cat{C} , \cat{C}}}$
	so that it sends a square $α : \doublehom[\doublehom[]{\hom{A}{C}}{\hom{B}{D}}{\prohom{A}{B}}{\prohom{C}{D}}]{f}{g}{m}{n}$ to
	$
		\comp
		{
			(\procomp{\tuple*{A , α} , \tuple*{m , g}}) ,
			(\procomp{\tuple*{f , n} , \tuple*{α , D}})
		}
	$ as shown.
	\begin{equation} \label{gray diagonal functor components}

	\end{equation}
	For consecutive arrows $f : \hom{A}{B}$ and $g : \hom{B}{C}$
	it has comparitor disk
	$\comparitor[2] \tuple*{f , g} ≔ \comp{\procomp{} \tuple*{A , f} , \tuple*{f , g} , \procomp{} \tuple*{g , C}}$,
	and for consecutive proarrows $m : \prohom{A}{B}$ and $n : \prohom{B}{C}$
	it has comparitor disk
	$\comparitor[2] \tuple*{m , n} ≔ \procomp{\comp{} \tuple*{A , m} , \tuple*{m , n} , \comp{} \tuple*{n , C}}$
	(suppressing the associator terms),
	as shown below.
	These act to \emph{collate} copies of composite arrows and proarrows, respectively.
	\begin{equation} \label{gray diagonal functor comparitors}

	\end{equation}
\end{definition}

If $\tensor{\cat{C} , \cat{C}}$ has lax interchangers in both dimensions then we can define another Gray diagonal functor
$\diag′_{\cat{C}} : \hom{\cat{C}}{\tensor{\cat{C} , \cat{C}}} ≔ \comp{\diag_{\cat{C}} , \swap{\cat{C}}{\cat{C}}}$,
which \refer[Gray swap functor]{swaps} the order of the factors.
Once we have lax functors we need transformations for them as well.

\begin{definition}[arrow-dimension oplax transformation of lax functors]
	An \define{arrow-dimension oplax transformation of lax functors} of double categories
	$α : \hom[\hom[]{\cat{C}}{\cat{D}}]{F}{G}$ has
	\refer[arrow transformation object-component arrow]{object-component arrows},
	\refer[arrow transformation proarrow-component square]{proarrow-component squares}, and
	\refer[arrow transformation arrow-component disk]{arrow-component disks},
	just like those of an \refer{arrow-dimension oplax transformation} of strict functors
	(definition \ref{definition: arrow-dimension oplax transformation}).
	It satisfies the naturality condition for squares \eqref{arrow transformation square naturality} as well.
	However, the proarrow composition preservation condition \eqref{arrow transformation proarrow composition preservation}
	no longer make sense because
	$$
		α (\procomp{} \, A) : \doublehom{α A}{α A}{F (\procomp{} A)}{G (\procomp{} A)}
		\quad \text{and} \quad
		α (\procomp{m , n}) : \doublehom{α A}{α C}{F (\procomp{m , n})}{G (\procomp{m , n})}
		,
	$$
	whereas
	$$
		\procomp{} (α A) : \doublehom{α A}{α A}{\procomp{} (F A)}{\procomp{} (G A)}
		\quad \text{and} \quad
		\procomp{α m , α n} : \doublehom{α A}{α C}{\procomp{F m , F n}}{\procomp{G m , G n}}
		.
	$$
	Likewise, the arrow composition compatibility condition \eqref{arrow transformation arrow composition compatibility}
	no longer make sense because
	$$
		α (\comp{} \, A) : \prohom{\comp{F (\comp{} \, A) , α A}}{\comp{α A , G (\comp{} \, A)}}
		\quad \text{and} \quad
		α (\comp{f , g}) : \prohom{\comp{F (\comp{f , g}) , α C}}{\comp{α A , G (\comp{f , g})}}
		,
	$$
	whereas\footnote
	{
		Here we repeat the boundary of $\procomp{} (α A)$ using \refer{disk notation}
		to facilitate comparison with that of $α (\comp{} \, A)$.
	}
	$$
		\procomp{} (α A) : \prohom{α A}{α A}
		\quad \text{and} \quad
		\procomp{(\comp{\procomp{} (F f) , α g}) , (\comp{α f , \procomp{} (G g)})} : \prohom{\comp{F f , F g , α C}}{\comp{α A , G f , G g}}
		.
	$$
	Therefore, we replace these with the following comparitor compatibility conditions.
	
	\paragraph{proarrow comparitor compatibility:}
	$
		\comp{\comparitor[0]^F \tuple*{A} , α (\procomp{} A)}
		=
		\comp{\procomp{} (α A) , \comparitor[0]^G \tuple*{A}}
	$,
	\begin{equation} \label{equation: arrow transformation lax nullary proarrow comparitor compatibility}

	\end{equation}
	and
	$
		\comp{\comparitor[2]^F \tuple*{m , n} , α (\procomp{m , n})}
		=
		\comp{(\procomp{α m , α n}) , \comparitor[2]^G \tuple*{m , n}}
	$.
	\begin{equation} \label{equation: arrow transformation lax binary proarrow comparitor compatibility}
		%
	\end{equation}
	
	\paragraph{arrow comparitor compatibility:}
	$
		\procomp{(\comp{\comparitor[0]^F \tuple*{A} , \procomp{} (α A)}) , α (\comp{} \, A)}
		=
		\comp{\procomp{} (α A) , \comparitor[0]^G \tuple*{A}}
	$,
	\begin{equation} \label{equation: arrow transformation lax nullary arrow comparitor compatibility}
		%
	\end{equation}
	and
	$
		\procomp{(\comp{\comparitor[2]^F \tuple*{f , g} , \procomp{} (α C)}) , α (\comp{f , g})}
		=
		\procomp{(\comp{\procomp{} (F f) , α g}) , (\comp{α f , \procomp{} (G g)}) , (\comp{\procomp{} (α A) , \comparitor[2]^G \tuple*{f , g}})}
	$.
	\begin{equation} \label{equation: arrow transformation lax binary arrow comparitor compatibility}
		%
	\end{equation}
\end{definition}

Proarrow-dimension transformations of lax functors are defined similarly.
A cubical modification of transformations of lax functors
is defined just as \refer[modification]{those involving strict functors} in definition \ref{definition: modification},
without additional coherences involving the object-component squares of the modification
and the comparitor disks of the lax functors.

We now have the ingredients that we need to define duplication for symmetric Gray-monoidal double categories.
Note that because a \refer{braided Gray-monoidal double category} $\cat{C}$
has invertible interchangers in both dimensions, the unitary lax functors
$\comp{\diag_{\cat{C}} , \tensor{,}_{\cat{C}}} \, , \, \comp{\diag′_{\cat{C}} , \tensor{,}_{\cat{C}}} : \hom{\cat{C}}{\cat{C}}$
are necessarily pseudo (i.e., have invertible comparitors).

\begin{definition}[duplication structure]
	A(n arrow-dimension) \define{duplication structure}
	for a \refer{symmetric Gray-monoidal double category} $\cat{C}$
	consists of the following structure.
	\begin{description}
		\item[{\define*[duplicator]{duplicators}}:]
			\refer[arrow-dimension oplax transformation of lax functors]{arrow-dimension oplax transformations}
			of unitary pseudo functors
			$$
				δ : \hom[\hom[]{\cat{C}}{\cat{C}}]{\comp{} \, \cat{C}}{\comp{\diag_{\cat{C}} , \tensor{,}_{\cat{C}}}}
				\quad \text{and} \quad
				δ′ : \hom[\hom[]{\cat{C}}{\cat{C}}]{\comp{} \, \cat{C}}{\comp{\diag′_{\cat{C}} , \tensor{,}_{\cat{C}}}}
			$$
		\item[{\define*[coassociator]{coassociators}}:]
			\refer[invertible modification]{invertible} globular \refer[modification]{modifications}
			with object-component disks 
			$$
				s A : \prohom{\comp{δ A , (\tensor{δ A , A})}}{\comp{δ A , (\tensor{A , δ A})}}
				\quad \text{and} \quad
				s′ A : \prohom{\comp{δ′ A , (\tensor{δ′ A , A})}}{\comp{δ′ A , (\tensor{A , δ′ A})}}
			$$
		\item[{\define*[cocommutor]{cocommutors}}:]
			\refer[invertible modification]{invertible} globular \refer[modification]{modifications}
			with object-component disks
			$$
				c A : \prohom{δ A}{\comp{δ′ A , \braid{A}{A}}}
				\quad \text{and} \quad
				c′ A : \prohom{δ′ A}{\comp{δ A , \braid{A}{A}}}
			$$
	\end{description}
	These must satisfy the following relations.
	\begin{description}
		\item[nullary tensor duplicator coherence:]
			Duplicating the tensor unit is trivial in the sense that
			\begin{equation} \label{nullary tensor duplicator coherence}
				δ \tensor{} = \comp{} \, \tensor{} = δ′ \tensor{}
			\end{equation}
		\item[binary tensor duplicator coherence:]
			Duplicating a tensor product is given by duplicating the factors sequentially
			and using the braiding to permute the result into the required order.
			\begin{equation} \label{binary tensor duplicator coherence}
				\begin{array}{l}
					δ (\tensor{A , X}) = \comp{(\tensor{A , δ X}) , (\tensor{δ A , X , X}) , (\tensor{A , \braid{A}{X} , X})}
					\quad \text{and} \\
					δ′ (\tensor{A , X}) = \comp{(\tensor{δ′ A , X}) , (\tensor{A , A , δ′ X}) , (\tensor{A , \braid{A}{X} , X})}
					\\
				\end{array}
			\end{equation}
		\item[homogeneous coassociator coherence:]
			The two arrow disks with each of the following boundaries
			built using only coassociator components and lax interchangers for duplicators are identified,
			as described by equation \eqref{equation: duplication homogeneous coassociator coherence} below.
			$$
				\begin{array}{l}
					\prohom{\comp{δ A , (\tensor{δ A , A}) , (\tensor{δ A , A , A})}}{\comp{δ A , (\tensor{A , δ A}) , (\tensor{A , A , δ A})}}
					\quad \text{and} \\
					\prohom{\comp{δ′ A , (\tensor{δ′ A , A}) , (\tensor{δ′ A , A , A})}}{\comp{δ′ A , (\tensor{A , δ′ A}) , (\tensor{A , A , δ′ A})}}
					\\
				\end{array}
			$$
		\item[coassociator cocommutor coherence:]
			The coassociators are interdefinable using the cocommutors,
			as described by equation \eqref{equation: duplication coassociator cocommutor coherence} below.
		\item[cocommutor syllepsis coherence:]
			The cocommutors are interdefinable using the syllepsis,
			as described by equation \eqref{equation: duplication cocommutor syllepsis coherence} below.
	\end{description}
\end{definition}

We now unravel this definition and introduce the corresponding graphical syntax.

\paragraph{\refer[duplicator]{duplicators}:}
For each object $A$
we have a \refer[arrow transformation object-component arrow]{component arrow}
$δ A : \hom{A}{\tensor{A , A}}$,
for each proarrow $m : \prohom{A}{B}$
we have a \refer[arrow transformation proarrow-component square]{component square}
$δ m : \doublehom{δ A}{δ B}{m}{\procomp{(\tensor{A , m}) , (\tensor{m , B})}}$,
and for each arrow $f : \hom{A}{B}$
we have a \refer[arrow transformation arrow-component disk]{component disk}
$δ f : \prohom{\comp{f , δ B}}{\comp{δ A , (\tensor{A , f}) , (\tensor{f , B})}}$.

We represent the component arrow $δ A$ in a string diagram in the $\tuple*{\tensor{,} , \comp{,}}$-plane
as a splitting of the line representing $A$ into two lines.
We represent the component squares $δ m$ and $δ f$ in surface diagrams
as a splitting of the surface containing $m$ or $f$ into two surfaces along a ``seam'',
with the proarrow or arrow in the second surface preceding that in the first surface in the relevant composition order,
as prescribed by $\diag$ in diagram \eqref{gray diagonal functor components}.
\begin{equation} \label{diagram: duplicator components}

\end{equation}
The proarrow-component squares and arrow-component disks of $δ′$
produce their copies in the order opposite that for $δ$.
In string diagrams we
use a black dot $\tikz[string diagram]{\node[black dot] {} ;}$ to represent $δ$,
and a white dot $\tikz[string diagram]{\node[white dot] {} ;}$ to represent $δ′$.
When referring to either duplicator generically we omit the distinguishing dot.

We defined $δ$ to be an oplax transformation of unitary pseudo functors.
Its domain, the identity functor, is of course a strict functor; regarded as a lax functor it has identity comparitor disks.
But its codomain, $\comp{\diag_{\cat{C}} , \tensor{,}_{\cat{C}}}$, has nontrivial comparitor disks
that collate copies using interchangers as shown in diagram \eqref{gray diagonal functor comparitors}.
The proarrow comparitor compatibility relation for $δ$ \eqref{equation: arrow transformation lax binary proarrow comparitor compatibility}
says that for consecutive proarrows $m : \prohom{A}{B}$ and $n : \prohom{B}{C}$ we have (up to suppressed proarrow associators)
$$
	δ (\procomp{m , n})
	=
	\comp
	{
		(\procomp{δ m , δ n}) ,
		(\procomp{\comp{} (\tensor{A , m}) , \interchange{m , n} , \comp{} (\tensor{n , C})})
	}
	,
$$
and the arrow comparitor compatibility relation \eqref{equation: arrow transformation lax binary arrow comparitor compatibility}
says that for consecutive arrows $f : \hom{A}{B}$ and $g : \hom{B}{C}$ we have
$$
	δ (\comp{f , g})
	=
	\procomp
	{
		[\comp{\procomp{} f , δ g}] ,
		[\comp{δ f , \procomp{} (\tensor{B , g}) , \procomp{} (\tensor {g , C})}] ,
		[\comp{\procomp{} (δ A) , \procomp{} (\tensor{A , f}) , \interchange{f , g} , \procomp{} (\tensor{g , C})}]
	}
	.
$$
These are represented by the following surface diagrams,
$$

$$
which have the following projection string diagrams in the $\tuple*{\procomp{,} , \comp{,}}$-plane.
$$
	%
$$

The naturality relation for squares \eqref{arrow transformation square naturality} says that for a square
$α : \doublehom[\doublehom[]{\hom{A}{A′}}{\hom{B}{B′}}{\prohom{A}{B}}{\prohom{A′}{B′}}]{f}{g}{m}{n}$
we have
$$
	\procomp{(\comp{α , δ n}) , δ g}
	≅
	\procomp
	{
		δ f ,
		(\comp
		{
			δ m ,
			[\procomp{(\tensor{A , α}) , \interchange{m , g}}] ,
			[\procomp{\interchange{f , n} , (\tensor{α , B′})}]
		})
	}
	.
$$
$$
	%
$$
The story for $δ′$ is the same, except for the order of the copies.

\paragraph{\refer[coassociator]{coassociators}:}
The invertible object-component disks for coassociators $s$ and $s′$
relate the two possible ways to make three ordered copies using their respective duplicators.
$$
	%
$$

\begin{lemma}[coassociator braiding coherence]
	Coassociators are coherent for braiding
	in the sense that we have the following equation between invertible arrow disks.
	$$
		%
	$$
	We have analogous braiding coherences involving $δ′$ and $s′$,
	as well as coherences where we duplicate $X$ rather than $A$.
\end{lemma}

\begin{proof}
	The last three squares on the left compose to $\braid{\comp{δ A , (\tensor{A , δ A})}}{X}$
	and the first three squares on the right compose to $\braid{\comp{δ A , (\tensor{δ A , A})}}{X}$
	by \eqref{equation: braiding homogeneous composite arrow} and \eqref{equation: (2,1)-braiding arrow coherence}.
	The result then follows from braiding naturality for $\tuple*{2 , 0}$-type squares \eqref{braiding naturality for squares (2,0)}.
\end{proof}

\paragraph{\refer[cocommutor]{cocommutors}:}
The invertible object-component disks for the cocommutor $c$
relate duplication by $δ$ to duplication by $δ′$ followed by a braiding.
$$

$$
The components for $c′$ swap the roles of $δ$ and $δ′$.
The naturality requirement for an arrow $f : \hom{A}{B}$ \eqref{arrow disk modification naturality for arrows},
$$
	%
$$
and for a proarrow $m : \prohom{A}{B}$ \eqref{arrow disk modification naturality for proarrows},
$$
	%
$$
reveal the need for distinct duplicators that differ in the order of their copies.
If we had only a single duplicator then the naturality equations could not hold because the boundaries wouldn't even agree.

\begin{lemma}[cocommutor braiding coherence]
	Cocommutors are coherent for braiding
	in the sense that we have the following equation between invertible arrow disks.
	$$
		%
	$$
	We have analogous relations obtained by swapping $\tuple*{δ , c}$ and $\tuple*{δ′ , c′}$, and by braiding on the other side.
\end{lemma}

\begin{proof}
	The last two squares on the left compose to $\braid{\comp{δ′ A , \braid{A}{A}}}{X}$ by \eqref{equation: braiding homogeneous composite arrow}.
	The result then follows from braiding naturality for $\tuple*{2 , 0}$-type squares \eqref{braiding naturality for squares (2,0)}.
\end{proof}

\paragraph{tensor duplicator coherence:}
Equation \eqref{nullary tensor duplicator coherence} for nullary tensor duplicator coherence is well-bounded because
$\tensor{}$ is a strict unit for $\tensor{,}$ so $\tensor{\tensor{} , \tensor{}} = \tensor{}$.
This ensures that the empty diagram continues to have a unique interpretation in the presence of duplication structure.

Equation \eqref{binary tensor duplicator coherence} for binary tensor duplicator coherence
says what it means to duplicate a tensor product of objects.
In particular, we duplicate the objects independently
and use the braiding to collate the copies. 
$$

$$
As with equation \eqref{equation: binary tensor braiding} for tensor braiding coherence
we are leaving some food on the table here by asking for strict equality rather than just an invertible arrow disk.
This choice seems consistent with our decision to use strictly associative and unital Gray-monoidal structure,
but may be worth revisiting in the future.

Because we have object-component arrows for $δ$
we also have for each object $A$ an arrow-component disk
$
	δ (δ A) : \prohom{\comp{δ A , δ (\tensor{A , A})}}{\comp{δ A , (\tensor{A , δ A}) , (\tensor{δ A , A , A})}}
$
depicted on the left, with proarrow-dimension boundaries shown on the right.
$$

$$
Similarly, we have disks $δ′ (δ′ A)$, $δ (δ′ A)$, and $δ′ (δ A)$.

\paragraph{homogeneous coassociator coherence:}
This is the dual of Mac Lane's associator coherence ``pentagon equation'' \cite{mac_lane-1998-categories},
but made a hexagon by the explicit interchanger.
$$
	%
$$
\begin{equation} \label{equation: duplication homogeneous coassociator coherence}
=
\end{equation}
$$
	%
$$
Likewise for $δ′$.

\paragraph{coassociator cocommutor coherence:}
The following composite of invertible arrow disks is equal to the coassociator $s A$.
$$
	%
	\quad \overset{\mathclap{\inv{s′ A}}}{\proarrow}
$$
\begin{equation} \label{equation: duplication coassociator cocommutor coherence}
	%
\end{equation}
Dually, we can express $s′ A$ in terms of $s A$ by swapping the roles of $δ$ and $δ′$.

\paragraph{cocommutor syllepsis coherence:}
The following composite of invertible arrow disks is equal to the cocommutor $c A$.
\begin{equation} \label{equation: duplication cocommutor syllepsis coherence}
	%
\end{equation}
Dually, we can express $c′ A$ in terms of $c A$ by swapping the roles of $δ$ and $δ′$.

Using cocommutors we can define ``heterogeneous coassociators''.
We can do so either using $s$:
$$
	%
$$

\begin{lemma}
	The above heterogeneous coassociators are equal.
\end{lemma}

\begin{proof}
	Substitute the term for $\inv{s}$ from \eqref{equation: duplication coassociator cocommutor coherence}
	into the first term above and cancel the resulting consecutive inverses.
	The first two terms of the result are a pair of consecutive cocommutors;
	rewrite them to a syllepsis by \eqref{equation: duplication cocommutor syllepsis coherence}.
	The final $\inv{\syllepsis}$ term is independent of all the intermediate terms; permute it to directly follow the $\syllepsis$ term.
	Finally, cancel this inverse pair, and what remains is the second term above.
\end{proof}

A duplication structure for a symmetric Gray-monoidal double category
allows us to make more than one ordered copy of an object, arrow, proarrow, or square.
We would like to be able to make fewer than one copy as well.
For this we will need deletion.

\begin{definition}[deletion structure]
	A(n arrow-dimension) \define{deletion structure}
	for a \refer{symmetric Gray-monoidal double category} $\cat{C}$
	with a \refer{duplication structure}
	consists of the following structure.
	\begin{description}
		\item[{\define*{deletor}}:]
			an \refer{arrow-dimension oplax transformation} of strict functors
			$$
				ε : \hom[\hom[]{\cat{C}}{\cat{C}}]{\comp{} \, \cat{C}}{\comp{\tuple{}_{\cat{C}} , \tensor{}_{\cat{C}}}}
			$$
		\item[{\define*[counitor]{counitors}}:]
			\refer[invertible modification]{invertible} globular \refer[modification]{modifications}
			with object-component disks
			$$
				\begin{array}{l l}
					l A : \prohom{\comp{δ A , (\tensor{ε A , A})}}{\comp{} \; A},
					&
					r A : \prohom{\comp{δ A , (\tensor{A , ε A})}}{\comp{} \; A},
					\\
					l′ A : \prohom{\comp{δ′ A , (\tensor{ε A , A})}}{\comp{} \; A},
					&
					r′ A : \prohom{\comp{δ′ A , (\tensor{A , ε A})}}{\comp{} \; A}
					\\
				\end{array}
			$$
	\end{description}
	These must satisfy the following relations.
	\begin{description}
		\item[nullary tensor deletor coherence:]
			Deleting the tensor unit is trivial in the sense that
			\begin{equation} \label{nullary tensor deletor coherence}
				ε \tensor{}  =  \comp{} \, \tensor{}
			\end{equation}
		\item[binary tensor deletor coherence:]
			Deleting a tensor product is given by deleting the factors sequentially
			\begin{equation} \label{binary tensor deletor coherence}
				ε (\tensor{A , X})  =  \comp{(\tensor{ε A , X}) , ε X}
			\end{equation}
		\item[counitor coassociator coherence:]
			The two arrow disks with each of the following boundaries
			built using only counitor and coassociator components are identified,
			as described by equation \eqref{counitor coassociator coherence} below.
			$$
				\begin{array}{l}
					\prohom{\comp{δ A , (\tensor{δ A , A}) , (\tensor{A , ε A , A})}}{\comp{δ A , (\tensor{A , δ A}) , (\tensor{A , ε A , A})}}
					\quad \text{and} \\
					\prohom{\comp{δ′ A , (\tensor{δ′ A , A}) , (\tensor{A , ε A , A})}}{\comp{δ′ A , (\tensor{A , δ′ A}) , (\tensor{A , ε A , A})}}
					\\
				\end{array}
			$$
		\item[counitor  cocommutor coherence:]
			The counitors are interdefinable using the cocommutors,
			as described by equation \eqref{equation: duplication counitor cocommutor coherence} below.
	\end{description}
\end{definition}

We now unravel this definition and introduce the corresponding graphical syntax.

\paragraph{\refer{deletor}:}
For each object $A$
we have a \refer[arrow transformation object-component arrow]{component arrow}
$ε A : \hom{A}{\tensor{}}$,
for each proarrow $m : \prohom{A}{B}$
we have a \refer[arrow transformation proarrow-component square]{component square}
$ε m : \doublehom{ε A}{ε B}{m}{\procomp{} \, \tensor{}}$,
and for each arrow $f : \hom{A}{B}$
we have a \refer[arrow transformation arrow-component disk]{component disk}
$ε f : \prohom{\comp{f , ε B}}{ε A}$.

We represent the component arrow $ε A$ in a string diagram in the $\tuple*{\tensor{,} , \comp{,}}$-plane
as a termination of the line representing $A$.
We represent the component squares $ε m$ and $ε f$ in surface diagrams
as a termination of the surface containing $m$ or $f$ along a ``tear''.
\begin{equation} \label{diagram: deletor components}
	\begin{tikzpicture}[string diagram , x = {(16mm , 0mm)} , y = {(0mm , -12mm)} , baseline=(align.base)]
		\coordinate (cell) at (1/2 , 9/16) ; 
		\coordinate (in) at (cell |- 0 , 0) ;
		\draw (in) to [out = south , in = north] node [fromabove] {A} (cell) ;
		\node [gray dot] at (cell) {} ;
		\node (align) at (1/2 , 1/2) {} ;
	\end{tikzpicture}
	\qquad , \qquad
	\begin{tikzpicture}[string diagram , x = {(16mm , -8mm)} , y = {(0mm , -24mm)} , z = {(12mm , 0mm)} , baseline=(align.base)]
		\coordinate (sheet) at (1/2 , 1/2 , 0) ;
		\coordinate (left top) at ($ (sheet) + (-1/2 , -1/2 , 0) $) ;
		\coordinate (right top) at ($ (sheet) + (1/2 , -1/2 , 0) $) ;
		\coordinate (left seam) at ($ (sheet) + (-1/2 , 1/16 , 0) $) ;
		\coordinate (right seam) at ($ (sheet) + (1/2 , 1/16 , 0) $) ;
		\draw [sheet]
			(left top) to coordinate [pos = 1/2] (top)
			(right top) to
			(right seam) decorate [torn] {to coordinate [pos = 1/2] (cell)
			(left seam)} to
			cycle ;
		\draw [on sheet] (top) to node [fromabove] {m} coordinate [pos = 1/2] (square) (cell) ;
		\node (align) at (sheet) {} ;
	\end{tikzpicture}
	\qquad , \qquad
	\begin{tikzpicture}[string diagram , x = {(16mm , -8mm)} , y = {(0mm , -24mm)} , z = {(12mm , 0mm)} , baseline=(align.base)]
		\coordinate (sheet) at (1/2 , 1/2 , 0) ;
		\coordinate (left top) at ($ (sheet) + (-1/2 , -1/2 , 0) $) ;
		\coordinate (right top) at ($ (sheet) + (1/2 , -1/2 , 0) $) ;
		\coordinate (left seam) at ($ (sheet) + (-1/2 , 1/16 , 0) $) ;
		\coordinate (right seam) at ($ (sheet) + (1/2 , 1/16 , 0) $) ;
		\draw [sheet]
			(left top) to
			(right top) to
			(right seam) decorate [torn] {to coordinate [pos = 1/2] (cell)
			(left seam)} to coordinate [pos = 2/3] (left)
			cycle ;
		\draw [on sheet] (left) to [bend left = 30] node [fromleft] {f} (cell) ;
		\node (align) at (sheet) {} ;
	\end{tikzpicture}
\end{equation}

We defined $ε$ to be an oplax transformation of strict functors.
Intuitively, the reason that the functors for deletor are strict rather than pseudo, as in the case for duplicator,
is that the role of the comparitors is to collate copies;
but when making zero copies of something, collation is trivial.
The proarrow composition preservation condition \eqref{arrow transformation proarrow composition preservation}
says that
$
	ε (\procomp{m , n})
	=
	\procomp{ε m , ε n}
$
and the arrow composition compatibility condition \eqref{arrow transformation arrow composition compatibility}
says that
$
	ε (\comp{f , g})
	=
	\procomp{(\comp{\procomp{} f , ε g}) , ε f}
$.
$$

$$

The naturality relation for squares \eqref{arrow transformation square naturality} says that for a square
$α : \doublehom{f}{g}{m}{n}$
we have
$
	\procomp{(\comp{α , ε n}) , ε g}
	≅
	\procomp{ε f , ε m}
$.
$$
	%
$$

\paragraph{\refer[counitor]{counitors}:}
The invertible object-component disks for the counitors relate duplicating an object and then deleting one copy
to doing nothing at all.
$$
	%
$$

\begin{lemma}[counitor braiding coherence]
	Counitors are coherent for braiding
	in the sense that we have the following equation between invertible arrow disks.
	$$
		%
	$$
\end{lemma}

\begin{proof}
	The first three squares on the left compose to $\braid{\comp{δ A , (\tensor{ε A , A})}}{X}$
	by \eqref{equation: braiding homogeneous composite arrow} and \eqref{equation: (2,1)-braiding arrow coherence}.
	Braiding an identity is itself an identity by \eqref{equation: braiding nullary composites},
	so the result follows from braiding naturality for $\tuple*{2 , 0}$-type squares \eqref{braiding naturality for squares (2,0)}.
\end{proof}

\paragraph{tensor deletor coherence}
Equation \eqref{nullary tensor deletor coherence} for nullary tensor deletor coherence ensures that
the empty diagram continues to have a unique interpretation in the presence of deletion structure.

Equation \eqref{binary tensor deletor coherence} for binary tensor deletor coherence
says what it means to delete a tensor product of objects, as shown on the left.
The choice of ordering is arbitrary and is related to the other possibility
by the isomorphism $\interchange{ε A , ε X}$, as shown on the right.
$$

$$
As with equation \eqref{binary tensor duplicator coherence} for binary tensor duplicator coherence,
we could have chosen to make this an invertible arrow disk rather than a strict equality.

\paragraph{counitor  coassociator coherence:}
This is the dual of Mac Lane's middle unit coherence ``triangle equation'' \cite{mac_lane-1998-categories}.
It says that the following cycle of invertible modification components is the identity.
\begin{equation} \label{counitor coassociator coherence}
	%
\end{equation}
Likewise for $δ′$.

\paragraph{counitor  cocommutor coherence:}
The following composite of invertible arrow disks is equal to the counitor $l A$.
\begin{equation} \label{equation: duplication counitor cocommutor coherence}
	%
	\; \overset{\eqref{equation: nullary tensor braiding}}{=} \;
	%
\end{equation}
Note that this factorization of $l A$ through $r′ A$
is equivalent to a factorization of $r′ A$ through $l A$
because all of the arrow disks involved are invertible.
Dually, we get factorizations of the counitors $l′ A$ and $r A$ through one another
by swapping the roles of $δ$ and $δ′$.


We have now described duplication and deletion structure for symmetric Gray-monoidal double categories
that is natural for squares by virtue of being defined in terms of arrow-dimension oplax transformations.
We feel justified in calling this structure ``cartesian'',
at least in the ``Fox'' sense described above.

\begin{definition}[cartesian Gray-monoidal double category]
	A  \define{cartesian Gray-monoidal double category}
	is a symmetric Gray-monoidal double category
	equipped with \refer[duplication structure]{duplication} and \refer[deletion structure]{deletion} structure.
\end{definition}

We currently do not have a universal construction characterization
of our proposed notion of cartesian structure for Gray-monoidal double categories.
One obvious strategy is to define duplication structure
using some sort of right adjoints to the
\refer[Gray diagonal functor]{Gray diagonal functors}
$\diag _{\cat{C}} , \diag′_{\cat{C}} : \hom{\cat{C}}{\tensor{\cat{C} , \cat{C}}}$
and deletion structure
using a right adjoint to the
unique functor $\tuple{}_{\cat{C}} : \hom{\cat{C}}{\cat{\product{}}}$
so that the duplicators and deletor emerge as the units of these adjunctions.
We have not yet managed to do this.

\paragraph{Related Constructions}
There is substantial work on cartesian -- in the sense of finite product -- structure for $2$-dimensional categories
using the ordinary, as opposed to the Gray, monoidal product.
In the globular setting, Carboni, Kelly, Walters, and Wood
define finite products for bicategories, in the sense of bilimits as natural equivalences of hom categories, 
and prove that they give rise to a canonical symmetric monoidal structure,
establishing a $2$-dimensional Fox theorem for bicategories \cite{carboni-2007-cartesian_bicategories_2}. 
It is the cubical and Gray version of this finite product structure
(with the topologically-motivated strictifications previously discussed)
that we have sought to capture here.

Carboni et al. reserve the term ``cartesian bicategory'' for a bicategory
where only the locally full sub-bicategory of maps is required to have bicategorical finite products,
while all hom categories must have $1$-categorical finite products locally.
The term ``maps'' here refers to $1$-cells that have right adjoints.
These can be regarded as the arrows of a double category with companion and conjoint proarrows for all arrows
(also known as a ``fibrant double category'' or ``proarrow equipment'').
Shulman has observed that in many cases a cartesian bicategory
can be regarded as the \refer{proarrow bicategory} of such a double category
\cite{shulman-2010-symmetric_monoidal_bicategories}. 

This perspective is further developed by Aleiferi \cite{aleiferi-2018-cartesian_double_categories},
who defines a notion of \emph{cartesian double category} in which the finite product structure
is given by right adjoints to the
cartesian diagonal functor $\diag _{\cat{C}} : \hom{\cat{C}}{\product{\cat{C} , \cat{C}}}$
and the unique functor to the singleton double category $\tuple{}_{\cat{C}} : \hom{\cat{C}}{\cat{\product{}}}$
in the $2$-category of double categories, functors that are strict for arrows and pseudo for proarrows,
and strict arrow-dimension transformations.
A related  technique has recently been used by Patterson to give a characterization
of \emph{double categories with iso-strong finite products},
which have unbiased cartesian products as well as local products \cite{patterson-2024-double_category_products}.

Aleiferi's concise adjoint-theoretic presentation of cartesian structure for double categories is compelling.
However, it is too strict for our purposes, relying on the ordinary, rather than the Gray, monoidal product,
and using strict, rather than oplax, arrow-dimension transformations.
Still, it serves as inspiration for our goal of a characterization by universal construction.

%% file: content/parts/conclusion.tex

In this paper we have characterized the algebraic structure
comprising double categories together with their functors, transformations, and modifications
as a \emph{locally cubical Gray category},
have embedded the classical Gray categories into the locally cubical ones, 
have identified their one-object instances as \emph{Gray-monoidal double categories},
have added \emph{braided}, \emph{sylleptic}, and \emph{symmetric} structure to these,
and have proposed a notion of \emph{cartesian} structure inspired by Fox's theorem.
Each of these constructions has been accompanied by a graphical representation in a calculus of surface diagrams.

We have compared our constructions with those from the literature,
where the existing globular constructions tend to be weaker than ours,
while the cubical ones tend to be stricter.
While our design choices have a pleasing correspondence to ordered homotopy,
the absence of a universal construction characterization
leaves us uncertain whether we have managed to capture the ``right'' notion of cartesian Gray-monoidal double category,
and further investigation is warranted.

%% file: content/parts/version_information.tex

This version of this article was published by \emph{Applied Categorical Structures}
in the \emph{Collection honouring Robert Paré for his 80th birthday}
as \emph{Volume 34, article number 65 (2026)}
(DOI: \url{https://doi.org/10.1007/s10485-026-09883-x}).

Since the previous version, minor errors are corrected
and the proof of proposition \ref{theorem: locally cubical Gray category of double categories} is given in more detail
using properties of transformations and modifications described at the end of section \ref{section: double categories}.